\documentclass[a4paper,11pt]{amsart}
\usepackage{times}
\usepackage[margin=1.3in]{geometry}

\usepackage{amsmath,amsfonts,amssymb,amsthm}
\usepackage[colorlinks=true, pdfstartview=FitV, linkcolor=blue, citecolor=blue, urlcolor=blue,pagebackref=false]{hyperref}
\newcommand{\ignore}[1]{}

\newtheorem{theorem}{Theorem}[section]
\newtheorem{definition}[theorem]{Definition}
\newtheorem{proposition}[theorem]{Proposition}
\newtheorem{lemma}[theorem]{Lemma}
\newtheorem{corollary}[theorem]{Corollary}
\numberwithin{equation}{section}
\newcommand{\N}{\mathbb{N}}

\newcommand{\R}{\mathbb{R}}

\newcommand{\Z}{\mathbb{Z}}
\newcommand{\eps}{\varepsilon}
\newcommand{\ph}{\varphi}
\newcommand{\Ed}{E_{\mathrm{diag}}}

\newcommand{\diam}{\hspace{-0.2ex}\diamond\hspace{-0.2ex}}
\newcommand{\del}{\delta\hspace{-0.3ex}}

\usepackage{scalerel,stackengine}
\stackMath
\newcommand\reallywidehat[1]{%
\savestack{\tmpbox}{\stretchto{%
  \scaleto{%
    \scalerel*[\widthof{\ensuremath{#1}}]{\kern-.6pt\bigwedge\kern-.6pt}%
    {\rule[-\textheight/2]{1ex}{\textheight}}%WIDTH-LIMITED BIG WEDGE
  }{\textheight}% 
}{0.5ex}}%
\stackon[1pt]{#1}{\tmpbox}%
}

\begin{document}

\title{Quasilinear SPDEs via rough paths}
\author{Felix Otto and Hendrik Weber}

\begin{abstract}
We are interested in (uniformly) parabolic PDEs with a nonlinear
dependence of the leading-order coefficients, driven by a rough right 
hand side. For simplicity, we consider a space-time periodic setting 
with a single spatial variable:
\begin{equation*}
\partial_2u -P( a(u)\partial_1^2u + \sigma(u)f ) =0
    \end{equation*}
where $P$ is the projection on mean-zero functions, and  $f$ is a distribution which is only controlled in the low regularity norm of
$ C^{\alpha-2}$ for $\alpha > \frac{2}{3}$ on the parabolic H\"older scale. 
 The example we have in mind is a random forcing $f$ and our assumptions 
  allow, for example,  for an $f$ which is white in the time variable $x_2$ and only mildly coloured in the space variable $x_1$; any spatial covariance operator
  $(1 + |\partial_1|)^{-\lambda_1 }$ with $\lambda_1 > \frac13$ is admissible.
 
 On the deterministic side we obtain
a $C^\alpha$-estimate for $u$,
assuming  that we control products of the form 
 $v\partial_1^2v$ and $vf$ with $v$
solving the constant-coefficient equation $\partial_2 v-a_0\partial_1^2v=f$. 
As a consequence, we obtain  existence, uniqueness and stability with respect to $(f, vf, v \partial_1^2v)$
of small space-time periodic solutions for small data.  We then demonstrate how the required
 products can be bounded  in the case of a random forcing $f$
using stochastic arguments.

For this we extend the treatment of the singular product $\sigma(u)f$ 
via a space-time version of Gubinelli's notion of controlled rough paths 
to the product $a(u)\partial_1^2u$, which has the same
degree of singularity but is more nonlinear since the solution $u$ appears in both factors.
In fact, we develop a theory for the linear equation $\partial_t u - P(a\partial_1^2 u +\sigma f)=0$
with rough but given coefficient fields $a$ and $\sigma$ and then apply a 
fixed point argument.
 The PDE ingredient mimics the (kernel-free) Safonov approach to ordinary Schauder theory.
\end{abstract}

\maketitle

\tableofcontents

\section{Introduction}
\addtocontents{toc}{\protect\setcounter{tocdepth}{2}}

%$\|f\|_{\alpha-2}:=\sup_{T\le 1}(T^\frac{1}{4})^{2-\alpha}\|f_T\|$
%
%$\|[u,(\cdot)]\diam f\|_{2\alpha-2}:=\sup_{T\le 
%1}(T^\frac{1}{4})^{2-2\alpha}\|[u,(\cdot)_T]\diam f\|$
%
%$\|[u,(\cdot)\diam\|_{2\alpha-2,2,1}:=\sup_{a_0,a_0'}\sup_{T\le 
%1}(T^\frac{1}{4})^{2-2\alpha}\|\{1,\frac{\partial}{\partial 
%a_0},\frac{\partial^2}{\partial a_0^2}\}\{1,\frac{\partial}{\partial 
%a_0'}\}[u,(\cdot)]\diam f\|$

We are interested in the parabolic PDE 
\begin{equation}\label{i-1}
\partial_2u-P\big( a(u)\partial_1^2u + \sigma(u)f \big) =0
\end{equation}
 for a rough driver $f$. The non-linearities $a, \sigma$ are assumed to be regular and uniformly elliptic, see \eqref{2.15} below for precise assumptions. 
 In order to avoid difficulties related to initial and boundary values we adopt a more elliptic point of view and seek solutions which are periodic \emph{both} in 
 the space-like coordinate $x_1$ and in the time-like coordinate $x_2$. This is the reason for the non-standard labelling of coordinates and the presence
 of the operator $P$, the projection onto mean-zero functions.
% and $P$ is the projection on mean-zero  functions.
 For the right hand side $f$ we only assume control on the low regularity norm of $C^{\alpha-2}$ in the parabolic H\"older scale for $\alpha \in (\frac23,1)$ (see \eqref{c.1} for a precise statement).
The optimal control on $u$ one could aim to obtain under these assumption is in the $C^\alpha$ norm but in this regularity class there is no
classical functional analytic definition of the singular products $a(u) \partial_1^2 u$ and $\sigma(u) f$. In this article we assume that we have 
 an ``off-line'' interpretation for the  products  $v \partial_1^2 v$, $vf$ (see \eqref{x93}), where $v( \cdot, a_0)$ is the mean-free and space-time periodic 
 solution to the constant coefficient equation 
 \begin{equation}\label{e:v-eq}
 \partial_2 v(\cdot, a_0) -a_0\partial_1^2v(\cdot, a_0)= Pf \qquad\qquad\text{distributionally}
 \end{equation}
  and show that these bounds allow us to control $u$. We are ultimately interested in a stochastic forcing $f$ and 
 in this case the required control of products 
 can be obtained using explicit moment calculations to capture stochastic cancelations. 

\medskip

Our method is similar in spirit to Lyons' rough path theory \cite{Lyons1,LyonsQian,LyonsCaruanaLevy}. This theory is based on the observation that
the analysis of stochastic integrals  
\begin{align}\label{i1}
\int_0^t u(s) dv(s)
\end{align}
for irregular $v$, such as Brownian motion or even lower-regularity stochastic processes,
can be conducted efficiently by splitting it into a stochastic and a deterministic step. In the stochastic step
the integral \eqref{i1} is defined for a single well-chosen function $\bar{u}$, e.g. $v$ itself.
In the case where $ \bar{u}=v$ is a (multidimensional) Brownian motion  there is a one-parameter family of canonical definitions 
for these integrals, with the It\^o and the Stratonovich notions being the most prominent ones.
Information on this single integral  suffices to give a subordinate sense to integrals
for a whole class of functions $u$  with similar small-scale behaviour. This line of thought is expressed precisely in Gubinelli's notion of a controlled path \cite[Definition 1]{GubinelliControlling}.
There, a function $u$ in the usual H\"older space $C^\alpha$, $\alpha \in (\frac13, \frac12)$,  is said to be controlled by $\bar{u} \in C^\alpha$ if there exists a third function $\sigma \in C^\alpha$ 
such that for all $s,t \in \R$
\begin{equation}\label{i-3}
| u(t) - u(s)  - \sigma (s) (\bar{u}(t) - \bar{u}(s)) | \lesssim |t-s|^{2\alpha}.
\end{equation}
Loosely speaking, this means that the {\it increments} 
$u(t)-u(s)$ of the function $u$ can be approximated by
those of $\bar{u}$ , 
provided the latter are locally modulated by the
amplitudes $\sigma$. 
In \cite[Theorem 1]{GubinelliControlling} it is then shown that this assumption, together with an ``off-line'' bound of the form
\begin{align}\label{i2}
\Big| \int_s^t  \bar{u}(r) dv(r) - \bar{u}(s) (v(t) - v(s)) \Big| \lesssim |t-s|^{2\alpha}, 
\end{align}
suffices to define the integral $\int u(r) dv(r)$ and to obtain the bound
\begin{align}
\notag
&\Big| \int_s^t u(r) dv(r) - u(s) (v(t) - v(s)) - \sigma(s) \int_s^t (\bar{u}(r) - \bar{u}(s)) dv(r) \Big|  \\
\label{i-17}
& \qquad \lesssim |t-s|^{3\alpha}.
\end{align}
The construction of the integrals \eqref{i2} for the specific function $\bar{u}$   can be accomplished under a less restrictive set of assumptions than required for the classical It\^o theory. In many applications
  this construction can be carried out using Gaussian calculus without making reference to an underlying martingale structure. 
  The construction makes no use of the linear ordering of time and lends itself well to extensions to higher-dimensional index sets.
  
\medskip
This last point was the starting point for Hairer's work on singular stochastic PDE -- the observation that the variable $t$ in the rough path theory
could represent ``space'' rather than ``time'' was the key insight that allowed him to define stochastic PDEs with non-linearities of Burgers type 
\cite{HairerBurgers} and the KPZ equation \cite{HairerKPZ}. The notion of controlled path was also the starting point for his definition of 
regularity structures \cite{HairerRegularity} which permits treatment of semilinear stochastic PDE with an extremely irregular right hand side, possibly involving a renormalization procedure. Parallel to that, Gubinelli, Imkeller and Perkowski put forward a notion of paracontrolled distributions
\cite{GubinelliImkellerPerkowski}, a Fourier-analytic variant of \eqref{i-3} which has also been used to treat singular stochastic PDE.

\medskip
In this article we propose yet another higher-dimensional generalization   of the notion of controlled path, see Definition~\ref{D1} below,
and  use it to provide a solution and stability theory for \eqref{i-1}.
This definition is an immediate generalization of Gubinelli's definition \eqref{i-3}  and also closely related to Hairer's notion \cite[Definition 3.1]{HairerRegularity} of a modelled 
distribution in a certain regularity structure. 
 However, the definition comes with a twist because the quasilinear nature of \eqref{i-1} forces us to allow the realization of the 
 \emph{model},
  $v( \cdot, a_0)$ in our notation, to depend on a parameter $a_0$, which (ultimately) corresponds to the variable diffusion coefficient $a(u)$. 
 In our theory the  ``off-line products'' $v f$ and $v \partial_1^2 v$ play the role of the ``off-line integral''$\int \bar{u} dv$ above and the regularity assumption 
 \eqref{i2} is translated into a control on the  commutators 
 \begin{align*}
 [v,(\cdot)_T] \diam \{ \partial_1^2 v, f\}:= v ( \{ \partial_1^2 v, f\})_T-(v \diam   \{ \partial_1^2 v, f\} )_T ,
  \end{align*}
 where $( \cdot )_T$ denotes the convolution with a smooth kernel at scale $T$ (see \eqref{1.10} and the discussion that follows it) and where we use the notation
 $\diam$ to indicate that products are not classically defined and that their interpretations have to be specified\footnote{In the literature $\diam$ is sometimes used to denote the Wick product of two random variables. Our 
 products need not be Wick products.}. Furthermore, here and below we use the abbreviated notation $[v,(\cdot)_T] \diam \{ \partial_1^2 v, f\}$ when we speak about $[v,(\cdot)_T] \diam \partial_1^2 v$ 
 and $[v,(\cdot)_T] \diam f$ simultaneously.
  Based on these assumptions we derive bounds in the spirit of \eqref{i-17} on the singular products $a(u) \diam \partial_1^2 u$ and $\sigma(u) \diam f$ (see Lemma~\ref{L2a} and \ref{L1}) which can also be seen as a  variant of Hairer's Reconstruction Theorem \cite[Theorem 3.10]{HairerRegularity}
   in a simpler situation. 
 We want to point out that our method completely avoids  the use of wavelet analysis which features prominently in Hairer's proof of the Reconstruction Theorem.
On the PDE side, in Lemma~\ref{L3},  we obtain an optimal regularity result on solutions $u$ of \eqref{i-1} based on a control of the commutators
$[a, (\cdot)_T] \diam \partial_1^2 u$ and $[\sigma, (\cdot)_T] \diam f$. This result is similar in spirit to Hairer's Integration Theorem~\cite[Theorem 5.12]{HairerRegularity}.
Our proof mimics Safonov's 
approach to Schauder theory (as popularized in the monograph \cite{Krylov})  and therefore does not make reference to a parabolic heat kernel. These
ingredients are combined in Proposition \ref{P}, to obtain a robust existence and uniqueness theory for the linear version of \eqref{i-1} (i.e. $a$ and $\sigma$ do not depend on 
$u$) including stability in the input data, and in Theorem~\ref{Theo} these results are used to develop a small data theory for 
the non-linear problem \eqref{i-1}.
We want to point out that the deterministic analysis does not depend on the 
assumption of a $1+1$ dimensional space and would go through completely unchanged if $\partial_2 - a(u) \partial_1$ were replaced by a 
uniformly parabolic operator $\partial_{n+1} - \sum_{i,j=1}^n  a^{ij}(u) \partial_i \partial_j  $ over $\R^n \times \R$.

\medskip

On the stochastic side, we consider a class of stationary Gaussian distributions $f$ of class $C^{\alpha-2}$. This class  includes, for example, the case where 
$f$ is ``white'' in the time-like variable $x_2$ and has covariance operator $(1+ |\partial_{1}|)^{-\lambda_1}$ for $\lambda_1 >\frac13$ in the $x_1$ variable, or the case where the noise is constant 
in the time-like variable $x_2$ and has covariance operator $(1+|\partial_1|)^{-\lambda_1}$ for $\lambda_1> -\frac{5}{3}$ for the $x_1$ variable (see the end of Section~\ref{s:stochastic}
for a more detailed discussion of admissible $f$).
%(this includes in particular a noise 
%term with regularity akin to the noise term in the two-dimensional Parabolic Anderson model). 
For such $f$ we construct the generalized products $v  \diam \partial_1^2 v$ and $v \diam f$ 
 as limits 
of renormalized smooth approximations: More precisely, let $\varphi$ be an arbitrary Schwartz function with $\int \varphi =1$ and for $\eps\in (0,1]$  set 
\begin{align}
\varphi_\eps(x_1,x_2):=\frac{1}{\eps^\frac{3}{4}}\varphi \Big( \frac{x_1}{\eps^\frac{1}{4}},\frac{x_2}{\eps^\frac{1}{2}}\Big),
\qquad 
f_\eps := f \ast \varphi_\eps, \label{i-def-fepsveps}
\qquad  v_\eps(\cdot, a_0) :=v(\cdot, a_0) \ast \varphi_\eps 
\end{align}
  and construct the $C^{\alpha-2}$ distributions  $v \diam f$ and $v  \diam \partial_1^2 v$  as %
\begin{align}
\notag
v( \cdot, a_0) \diam f
:=& 
\lim_{\eps \to 0} 
	\big( 
		v_\eps( \cdot, a_0) f_\eps 
		- \big\langle 
			v_\eps( \cdot, a_0) f_\eps  
		\big\rangle 
	\big),
\\ %\notag 
v( \cdot, a_0) \diam \partial_1^2 v(\cdot, a_0') 
:=& 
\lim_{\eps \to 0} 
	\big( 	
		v_\eps( \cdot, a_0) \partial_1^2 v_\eps(\cdot, a_0') 
%\\
\label{Into5} %& 
%\qquad \qquad \qquad \qquad  - 
		-\big\langle 
			v_\eps( \cdot, a_0) \partial_1^2 v_\eps(\cdot, a_0') 
		\big\rangle 
	\big), 
\end{align} 
where we use angled brackets $\langle \cdot \rangle$ for the expectation of a random variable, see Proposition~\ref{P2} below. 
(We use the non-standard scaling in $\eps$ for consistency  spatial scaling given by the convolution with $\Psi_T$, see Section~\ref{s:setup} below).
%
%In many of the examples we consider,  the expectations of the regularized  products 
%
%
%\begin{align}\notag
%g_1(\eps, a_0)& = \langle v_\eps( \cdot, a_0) f_\eps  \big\rangle \\
%g_2(\eps, a_0,a_0') &= \langle v_\eps( \cdot, a_0) \partial_1^2 v_\eps(\cdot, a_0')\rangle 
%\label{Intro77}
%\end{align}

\medskip

The construction of these  renormalized products  and the deterministic well-posedness theory can be combined 
to the following theorem:

\begin{theorem}\label{Theo:Introduction1}
Let the non-linearities $a, \sigma$ be smooth and uniformly elliptic in the sense that 
\begin{align}\label{2.15}
\begin{array}{c}
a\in[\lambda,\frac{1}{\lambda}],\quad\|a'\|,\|a''\|,\|a'''\|\le\frac{1}{\lambda},\\[1ex]
\sigma\in[-1,1],\quad\|\sigma'\|,\|\sigma''\|,\|\sigma'''\|\le\frac{1}{\lambda},
\end{array}
\end{align}
where $\lambda>0$ is some fixed constant and $\| \cdot \|$ denotes the supremum norm.
Let $f$ be a  space-time periodic random Schwartz distribution, which is stationary, centered and Gaussian, and which satisfies the regularity assumption \eqref{A1} for
$\frac23<\alpha'<1$ and let $\alpha$ satisfy  $\frac23< \alpha <\alpha'$. 
Let $f_\eps$ be as in \eqref{i-def-fepsveps}.

\medskip

For any noise amplitude $\eta>0$ we consider the following regularized and renormalized version of \eqref{i-1}
\begin{align}
\notag
&\partial_2 u_\eps - P\big( a(u_\eps) \partial_1^2 u_\eps  - a'(u_\eps) \sigma^2(u_\eps) \eta^2 g_2(\eps, a(u_\eps), a(u_\eps))  \\
\label{glatteGleichung}
&\qquad \qquad \qquad+ \sigma(u_\eps) \eta f_\eps 
- \sigma'(u_\eps) \sigma(u_\eps) \eta^2 g_1( \eps, a(u_\eps)) \big) =0,
\end{align}
where 
\begin{align}\notag
g_1(\eps, a_0)& := \langle v_\eps( \cdot, a_0) f_\eps  \big\rangle, \\
g_2(\eps, a_0,a_0') &:= \langle v_\eps( \cdot, a_0) \partial_1^2 v_\eps(\cdot, a_0')\rangle ,
\label{Intro77}
\end{align}
and where $v_\eps( \cdot , a_0) $ is defined in \eqref{i-def-fepsveps}.
%= v (\cdot, a_0) \ast \varphi_\eps$ and $v(\cdot, a_0)$ is the unique zero-mean and space-time periodic solution of \eqref{e:v-eq}. 
\medskip

%\medskip
There exists a random  $\eta_0>0$ and a deterministic constant $\delta = \delta(\lambda,\alpha) \in (0,1]$ such that almost surely for any $\eta\leq \eta_0$ and for any 
$0<\eps\leq1$ there exists a unique space-time periodic smooth random function  $u_\eps$ which satisfies \eqref{glatteGleichung} and which is small
in the sense $[u_\eps]_{\alpha} \leq \delta$, where $[u_\eps]_{\alpha}$ refers to the parabolic H\"older semi-norm, defined in \eqref{w12}. Furthermore
$\eta$ is not too small in the sense that
\begin{align}\label{eta-moment-bound}
\langle \eta_0^{-p} \rangle^{\frac{1}{p}} < \infty\qquad \qquad \text{for all $p<\infty$.}
\end{align}
\medskip

Almost surely, for any fixed $\eta\leq\eta_0$ the solutions $u_\eps$ converge as $\eps \downarrow 0$ to a limit $u$. This convergence takes place 
uniformly and with respect to $[ \cdot ]_{\alpha}$. The limit $u$ does not depend on the choice of mollifying kernel $\varphi$ although  $g_1$ and $g_2$ do.
\end{theorem}
The small amplitude  $\eta$ appears here because of our choice to work with space-time  periodic solutions rather than 
treating the initial value problem (space-time periodic here means that functions/distributions are periodic of fixed period which 
without loss of generality we set  to be $1$, 
both in the space-like coordinate $x_1$ and the time-like coordinate $x_2$).
 In initial value problems it is common to show ``local'' existence and uniqueness of 
solutions, i.e. existence and uniqueness on some small time interval (the length of which is random if there are random 
terms in the equation). The small amplitude $\eta$ plays the role of this small time interval here. The smallness assumption 
 $[u_\eps]_\alpha\leq\delta$ also appears because of the periodic space-time boundary conditions and is needed to ensure 
 uniqueness of solutions.
%The precise form of  the functions $g_1$ and $g_2$, which are defined in \eqref{Intro77},
%is given in Lemma~\ref{lem6} below. 
%In general under our Assumption \eqref{A1} on the law of the noise $f$, they both 
%diverge as $\eps \downarrow 0$. 
%If  $f$ satisfies the additional stronger regularity assumption \eqref{A2},
%$g_1$ and $g_2$ converge to a finite limit as $\eps$ to $0$, and the statement of Theorem~\ref{Theo:Introduction1}
% remains true without renormalisation, i.e. for $g_1= g_2=0$. This is stronger assumption holds, for example, if $f$ is ``white'' in $x_1$ and ``trace-class'' in $x_2$. 
%
%\medskip
%
The following theorem gives a characterization of the limit $u$ obtained in Theorem~\ref{Theo:Introduction1}.

\begin{theorem}\label{Theo:Introduction2}
\medskip
Under the assumptions of Theorem~\ref{Theo:Introduction1}, $u$ is almost surely the unique mean-free space-time periodic
function with the properties
\begin{align}
%\notag 
u\;\mbox{is modelled after}\;v\;\mbox{according to}\;a(u)\;\mbox{and}\;\sigma(u) 
\mbox{ (in the sense of Definition \ref{D1})},\label{v95A}\\
\partial_2u-P\big[ a(u)\diam\partial_1^2u+\sigma(u)\diam \eta f \big]=0\quad\mbox{distributionally},\label{v94A}
\end{align}
satisfying 
\begin{align}
[u]_\alpha \leq \delta \label{v96bisA}.
\end{align}
\end{theorem}
We stress that the definition  of the non-standard products  $a(u)\diam\partial_1^2u$ and $\sigma(u)\diam \eta f$ in \eqref{v94A}
(see Corollary~\ref{C1} and Lemma~\ref{L1})
relies on the "modelledness" of $u$ as well as the definition of the renormalized products \eqref{Into5}. 

\medskip

We finally mention that briefly before posting the second version of our result,  the article 
\cite{GubinelliFurlan} was posted on the arXiv. In this article Furlan and Gubinelli study the equation
\begin{equation}
\partial_t u - a(u) \Delta u = \xi ,
\end{equation}
where $u = u(t,x)$ for $x$ taking values in the two-dimensional torus, and  $\xi = \xi(x)$ is 
a white noise over the two-dimensional torus, which is constant in the time variable $t$. This
noise term $\xi$ is of class $C^{-1-}$ and therefore essentially behaves like our term $f$.
They also define a notion of solution and prove short time existence and uniqueness of solutions 
for the initial value problem, as well as convergence for renormalized approximations  similar to  \eqref{glatteGleichung}. 
Following the approach we present here, they locally approximate the 
solutions $u$ by a family of solutions to constant coefficient problems. Their approach then proceeds in
 the framework of paracontrolled distributions. 
 Yet another approach    by Bailleul, Debussche and Hofmanova \cite{BailleulDebusscheHofmanova}
 was put forward shortly after posting our second version.
 They deal with the system
 \[
 \partial_t u - a(u)\Delta u = g(u) \xi,
 \]
 where $\xi$ is again a two dimensional white noise and they also obtain a short-time existence and 
 stability result for renormalized solutions. Their method is easier than ours or Furlan and Gubinelli's 
 as they only need a single random function, namely $X = (-\Delta)^{-1} \xi$ to locally describe $u$. 
 However, this makes strong use of the fact that the noise only depends on the space variable 
 and it would also not work if the operator $a(u) \Delta$ were replaced by the more 
 general uniformly elliptic operator $a_{ij}(u) \partial_{i} \partial_j u$.

\bigskip

\section{Setup}
\label{s:setup}

%\subsection*{Setup}

%\subsubsection*{Metric}
%{\sc Metric}. 
The parabolic operator $\partial_2-a_0\partial_1^2$ and its
mapping properties on the scale of H\"older spaces (i.e. Schauder theory) imposes its intrinsic 
(Carnot-Carath\'eodory) metric, which is given by
\begin{align}\label{1.11}
d(x,y):=|x_1-y_1|+\sqrt{|x_2-y_2|},
\end{align}
see for instance \cite[Section 8.5]{Krylov}.
The H\"older semi-norm $[\cdot]_\alpha$ is defined based on (\ref{1.11}):
\begin{align}\label{w12}
[u]_\alpha:=\sup_{x\not= y}\frac{|u(x)-u(y)|}{d^\alpha(x,y)}.
\end{align}

\medskip
%\subsubsection*{Convolution} 
In order to define negative norms of distributions in an intrinsic way, cf. \eqref{c.1} below, 
 it is convenient to have a family $\{(\cdot)_T\}_{T>0}$ of mollification operators $(\cdot)_T$ consistent 
with the relative scaling $(x_1,x_2)=(\ell\hat x_1,\ell^2\hat x_2)$ of the two variables
dictated by (\ref{1.11}).  
It will turn out to be extremely convenient to have in addition the semi-group property
\begin{align}\label{1.10}
(\cdot)_T \circ (\cdot)_t=(\cdot)_{T+t}.
\end{align}
All is achieved by convolution with the semi-group $\exp(-T(\partial_1^4-\partial_2^2))$
of the elliptic operator $\partial_1^4-\partial_2^2$, which is the simplest positive operator
displaying the same relative scaling between the variables as $\partial_2-\partial_1^2$ and being symmetric
in $x_2$ and $x_1$. We note that the corresponding convolution kernel $\psi_T$ is easily
characterized by its Fourier transform $\hat\psi_T(k)=\exp(-T(k_1^4+k_2^2))$; since
the latter is a Schwartz function, also $\psi_T$ is a Schwartz function. The only two (minor) inconveniences
are that 1) the $x_1$-scale is played by $T^\frac{1}{4}$ (in line with (\ref{1.11}) the $x_2$-scale
is played by $T^\frac{1}{2}$) since we have
$\psi_T(x_1,x_2)=\frac{1}{T^\frac{3}{4}}\psi_1(\frac{x_1}{T^\frac{1}{4}},\frac{x_2}{T^\frac{1}{2}})$ 
and that 2) $\psi_1$ (and thus $\psi_T$) does not have a sign. The only properties of the kernel we need
are moments of derivatives: 
\begin{align}\label{1.13}
%\begin{array}{rcl}
\int dy|\partial_1^k \partial_2^\ell \psi_T(x-y)|d^\alpha(x,y) \leq C(k,\ell,\alpha) (T^\frac{1}{4})^{-k-2\ell+\alpha}
%\int dy|\partial_2^k\psi_T(x-y)|d^\alpha(x,y)&\lesssim&(T^\frac{1}{4})^{-2k+\alpha}
%\end{array}
\end{align}
for all orders of derivative $k,\ell=0,1,\cdots$ and moment exponents $\alpha\ge 0$, as well as the fact that 
$\int \psi(x) x_1 dx =0$.  
 Estimates (\ref{1.13}) follow immediately
from the scaling and the fact that $\psi_1$ is a Schwartz function.
In Lemma  \ref{L:AMMM1} we show however that our main regularity assumption  \eqref{c.1} on $f$ as well as the bounds on 
the commutators do not depend on the specific  choice of Schwartz kernel $\psi$. In particular, the statements 
ultimately  do not depend on the semi-group property although this property plays an important part in the proofs. 
%
%
%\medskip
%
%\subsubsection*{Finite domain} 
%We mimic a finite domain by imposing periodicity in both directions;
%w.l.o.g. we may set this scale equal to one. 
We will typically measure the size of the
distribution $f$ by the expression 
\begin{equation} \label{c.1}
\| f\|_{\alpha-2}:= \sup_{T\le 1}(T^\frac{1}{4})^{2-\alpha}\|f_T\|,
\end{equation}
%%%
%%% Here we could introduce the constant $N_0$ &&&&

where the restriction $T\le 1$ reflects the period unity. By Lemma \ref{LA1}, cf. Step \ref{LA1S1}, 
 this expression agrees with  the standard definition of the norm of $C^{\alpha-2}$.  
% Here and below 
% we use the convention to write 
%% 
% \begin{align}\label{conven}
% \| A(T)\|_{T,\beta} = \sup_{T\leq 1} (T^{\frac14})^{-\beta}\| A(T)\|,
% \end{align}
%%

%\medskip
%
%\subsubsection*{Standing assumptions on the nonlinearities} We assume throughout 
%that the coefficients $a$ and $\sigma$ satisfy the regularity and uniform ellipticity 
%assumption \eqref{2.15} for some constant  $\lambda>0$.
%%%
%%\begin{align}\label{2.15}
%%\begin{array}{c}
%%a\in[\lambda,\frac{1}{\lambda}],\quad\|a'\|,\|a''\|,\|a'''\|\le\frac{1}{\lambda},\\[1ex]
%%\sigma\in[-1,1],\quad\|\sigma'\|,\|\sigma''\|,\|\sigma'''\|\le\frac{1}{\lambda}.
%%\end{array}
%%\end{align}
%%%
%We express the bound on the various norms of $a$ and $\sigma$ by the ellipticity contrast $\lambda$
%in order to have a single constant that measures the quality of the data. Note that the assumption
%$\sigma\in[-1,1]$ is only seemingly stronger than $\|\sigma\|\lesssim\frac{1}{\lambda}$ since
%that constant can always be absorbed into the rhs $f$ in the equation.
%These fairly high regularity assumptions intervene in the proof of Lemma \ref{L0},
%they could be slightly weakened in the sense of \cite[Proposition 4]{GubinelliControlling}, at the expense
%of a more complicated notation. 
%
\medskip

Here and throughout the entire deterministic analysis presented in Sections~\ref{s:deterministic},~\ref{sec:DetProofs}
and Appendix~\ref{s:Appendix}  $\lesssim$ means
$\le C$ with a constant $C$ only depending on $\lambda$ and the exponent $\alpha$.
In the derivation of the stochastic bounds in Sections~\ref{s:stochastic} and~\ref{s:stoch-proofs}
the implicit constant may depend on additional parameters which are specified there.
Similarly, we write $\ll 1$ for $\leq \delta$ for $\delta= \delta(\lambda, \alpha)$ small enough.

\section{Deterministic Analysis}
\label{s:deterministic}

We start with the following central definition which is a straightforward generalization of
Gubinelli's definition \cite[Definition 1]{GubinelliControlling}
of a ``controlled path'', a generalization from the time variable $x_2$ to multiple variables $x$,
and to a ``model'' $(v_1,\cdots,v_I)$ (in the language of Hairer \cite{HairerRegularity}) 
that here may depend on an additional parameter $a_0$.  It
states that the {\it increments} 
$u(y)-u(x)$ of the function $u$ can be approximated by
those of several functions $v_i$, 
if  the latter are locally modulated by the
amplitudes $\sigma_i$ and the functions $a_i$ that locally determine the value of the parameter $a_0$.
The functions $\sigma_i$ can therefore be interpreted as ``derivatives'' of 
$u$ with respect to $v_i$.
The increments of the  linear function $x_1$ also have to be included
because of $\alpha>\frac{1}{2}$.
In fact, since $2\alpha>1$, given the model $(v_1,\cdots,v_I)$ (as modulated by the functions $a_i$),
the ``derivatives'' $(\sigma_1,\cdots,\sigma_I)$ and $\nu$ determine $u$ up to a constant.
In our situation, we expect $u$ and $(v_1,\cdots,v_I)$ 
to be H\"older continuous with exponent not (much) larger than $\alpha$,
so that imposing closeness of the increments to order $2\alpha$ contains valuable additional information.
%However, since $2\alpha<2$, the functions $(a_1,\cdots,a_I)$ and $(\sigma_1,\cdots,\sigma_I)$ 
%determine $(u,\nu)$ only up to adding $(\del u,\partial_1\del u)$ with $\del u$ a function of class $C^{2\alpha}$.

\begin{definition}\label{D1}
Let $\frac{1}{2}<\alpha<1$ and $I\in\mathbb{N}$.
We say that a function $u$ is modelled after the functions $(v_1,\cdots,v_I)$ of $(x,a_0)$
according to the functions $(a_1,\cdots,a_I)$ and $(\sigma_1,\cdots,\sigma_I)$ provided there exists a function $\nu$ 
(which because of $2\alpha>1$ is easily seen to be unique) such that
\begin{align}\label{1.3}
M:=&\sup_{x\not=y}\frac{1}{d^{2\alpha}(y,x)} \nonumber\\
&
|u(y)-u(x)-\sigma_i(x)(v_i(y,a_i(x))-v_i(x,a_i(x)))-\nu(x)(y-x)_1|
\end{align}
is finite. Here and in the sequel we use Einstein's convention of summation over repeated indices.
\end{definition}

Note that imposing (\ref{1.3}) also for distant points $x$ and $y$ is consistent with periodicity 
despite the non-periodic term $(y-x)_1$ since by $\alpha\ge\frac{1}{2}$ the latter is dominated by
$d^{2\alpha}(x,y)$ for $d(x,y)\ge 1$. 
Note also that (\ref{1.3}) is reminiscent of a H\"older norm: In case of $(\sigma_1,\cdots,\sigma_I)=0$,
the finiteness of (\ref{1.3}) implies that $u$ is continuously differentiable in $x_1$ and that
$\nu(x)=\partial_1u(x)$ so that $M$ turns into the parabolic $C^{2\alpha}$-norm of $u$.
In this spirit, Step \ref{L2aS0} in the proof of Lemma \ref{L2a} shows that the modelledness constant $M$ in (\ref{1.3}) controls
the $(2\alpha-1)$-H\"older norm of $\nu$, provided $x\mapsto\sigma_i(x)v_i(\cdot,a_i(x))$ is $\alpha$-H\"older continuous
with values in $C^\alpha$.
In addition, in the presence of periodicity, $M$ also controls the $\alpha$-H\"older norm of $u$ and the
supremum norm of $\nu$, which are of lower order, cf. Step \ref{L2aS-1} in the proof of Lemma \ref{L2a}.

%%%%%%%%%%%%%%%%%%%%%%%%%%%%%%%%%%%%%%%%%%%%%%%%%%%%%%%%%%%%%%%%%%%%%%%%%%%%%%%%%%%%%%%%%%%%%%%%%%%%%%%
\medskip

The following lemma shows that the notion of modelledness in Definition \ref{D1}
is well-behaved under sufficiently smooth nonlinear {\it pointwise} transformation; it will
be used in the proof of Theorem \ref{Theo}. It is essentially
identical to \cite[Proposition 4]{GubinelliControlling}, which in turn is a consequence of Taylor's formula; 
and we omit the proof.
%because of the minor modifications due to the presence of a more general model, we reproduce the proof.

\begin{lemma}\label{L0}
Let $\frac12<\alpha<1$.\\ 
i) Suppose that $u\in C^\alpha$ is modelled after $v $ according to $a$ and $\sigma$ with constant $M$.
Let the function $b$ be twice differentiable. Then $b(u)$ is modelled after $v$ according to
$a$ and $\mu:=b'(u)\sigma$ with constant $\tilde M$ estimated by
\begin{align}
\tilde M+[b(u)]_\alpha&\le(\|b'\|+\|b''\|[u]_\alpha)(M+[u]_\alpha),\label{wj14}\\
[\mu]_\alpha+\|\mu\|&\le(\|b'\|+\|b''\|[u]_\alpha)([\sigma]_\alpha+\|\sigma\|).\label{wj13}
\end{align}
ii) Suppose that for $i=0,1$, $u_i \in C^\alpha$ is modelled after $v_i $ according to $a_i$ and $\sigma_i$ with constant $M_i$.
Suppose further that $u_1-u_0$ is modelled after
$(v_1,v_0)$ according to $(a_1,a_0)$ and $(\sigma_1,-\sigma_0)$ with constant $\del M$.
Let the function $b$ be three times differentiable. Then $b(u_1)-b(u_0)$
is modelled after $(v_1,v_0)$ according to 
$(a_1,a_0)$ and $(\mu_1:=b'(u_1)\sigma_1, -\mu_0:= -b'(u_0)\sigma_0)$ with constant $\del\tilde M$
estimated by
\begin{align}
\del\tilde M&+[b(u_1)-b(u_0)]_\alpha+\|b(u_1)-b(u_0)\| \nonumber\\
&\le \Big(\|b'\|+\|b''\|(\max_iM_i+\max_i[u_i]_\alpha) 
+\|b'''\|(\max_i[u_i]_\alpha)^2\Big)
\nonumber\\
&  \qquad \times (\del M+[u_1-u_0]_\alpha+\|u_1-u_0\|),\label{wj11}\\
[\mu_1&-\mu_0]_\alpha+\|\mu_1-\mu_0\| \nonumber\\
&\le(\|b'\|+\|b''\|\max_i[u_i]_\alpha)([\sigma_1-\sigma_0]_\alpha+\|\sigma_1-\sigma_0\|)\nonumber\\
&+(\|b''\|+\|b'''\|\max_i[u_i]_\alpha)\max_i([\sigma_i]_\alpha+\|\sigma_i\|)
%\nonumber\\
 ([u_1-u_0]_\alpha+\|u_1-u_0\|).\label{wj12}
\end{align}

\end{lemma}

\medskip
%%%%%%%%%%%%%%%%%%%%%%%%%%%%%%%%%%%%%%%%%%%%%%%%%%%%%%%%%%%%%%%%%%%%%%%%%%%%%%%%%%%%%%%%%%%%%%%%%%%%%%%

As discussed in the introduction, the main  challenge in solving stochastic ordinary differential equations  is to give a sense
to integrals of the form \eqref{i1}. In the spirit of Hairer  \cite{HairerRegularity}  we interpret this problem as giving a meaning 
to the {\it product} $u \partial_t v$, which does not have a canonical functional analytic definition because both $u$ and $v$
 are only H\"older continuous 
in the time variable $t$ of exponent less than $\frac{1}{2}$, because they behave like Brownian motion. 
In view of the parabolic scaling, we encounter
the same difficulty when giving a distributional sense to $b\diam\partial_1^2u$ when
$b$ and $u$ are only H\"older continuous of exponent $\alpha<1$ (from now we use the non-standard notation 
$b \diam \partial_1^2 u$ instead of $b\, \partial_1^2 u$ to indicate that the definition of this product is non-standard). 

\medskip
As discussed in the introduction a main insight of Lyons' theory of rough paths, 
was the observation that such products can be defined provided $u$ is controlled by $\bar u$ and
the off-line product $\bar{u} \partial_t v$ 
satisfies the bound  \eqref{i2}, which can be rewritten as
 $\int_{s}^{t}(\bar{u}(r)-\bar{u}(s))\diam\partial_r v(r)$ $=-\bar u(s)\int_{s}^{t} \partial_r v(r)  -\int_{s}^{t} \bar u\diam \partial_r v $
$=:-([\bar u,\int^{t}]\diam \partial_r v)(s)$, that is, the expression on both sides of (\ref{i2}) amount 
to a commutator $[\bar{u},\int^{t}]$ of multiplication with 
$\bar{u}$ and integration, applied to a distribution $\partial_r v$. In our multi-dimensional framework, we replace
integration $\frac{1}{t-s}\int_{s}^{t}$ by (smooth) averaging:
\begin{align}\label{v38}
[v,(\cdot)_T]\diam f:=vf_T-(v\diam f)_T.
\end{align}
It is (only the control of) $[v,(\cdot)_T]\diam f$ that relates the distribution $v\diam f$ to
the function $v$ and the distribution $f$. 
In our set up, the role of the crucial ``algebraic relationship'' \cite[(24)]{GubinelliControlling} from 
rough path theory is played by 
the following straightforward consequence of the semi-group property (\ref{1.10})
\begin{align}\label{v36}
[v,(\cdot)_{t+T}]\diam f-([v,(\cdot)_T]\diam f)_t=[v,(\cdot)_t]f_T,
\end{align}
cf. (\ref{wi16}) in the proof of Lemma \ref{L2a}.
We also stress that the bound of order $(T^{\frac14})^{2\alpha-2}$ on the commutator \eqref{v38}
we impose below, is equivalent to the condition on the ``model'' imposed in \cite{HairerPardoux} in the framework 
of regularity structures. In fact, there in  \cite[Equation (3.9)]{HairerPardoux}  the condition  (in their notation)
\begin{equation*}
|(\Pi_z)(\tau) (\varphi^{\lambda}_z)| \lesssim \lambda^{|\tau|},
\end{equation*}
is assumed for all ``stochastic basis elements'' $\tau$. Specialized to $\tau = \mathcal{I}(\Xi) \Xi$  (still in their notation)
and following the definition  of the ``canonical admissible model'' (see \cite[Section 3.3]{HairerPardoux}) this condition
translates to our notation as
\begin{equation*}
\| [ v , \ast \varphi^{\lambda}] \diam f \| =
\sup_{x_0}\Big| \int (v \diam f(x) - v(x_0) f(x))  \varphi^{\lambda}(x-x_0) dx \Big| \lesssim \lambda^{|\tau|}, 
\end{equation*}
where $\varphi^\lambda$ is a (parabolically) scaled test-function. This only differs from our assumption in our specific 
choice of regularising kernel $\psi$.

\medskip

For our quasilinear SPDE, we need to give a sense to the {\it two} singular products $\sigma(u)\diam f$
and $a(u)\diam\partial_1^2u$,
so in particular to products of the form $u\diam f$ and $b\diam\partial_1^2u$,
where $u$ and $b$ behave $v$ defined by \eqref{e:v-eq}. Hence we will
need the two off-line products $v\diam f$ and $v\diam\partial_1^2v$. For simplicity, we split
the argument into Lemma \ref{L2a} and Corollary~\ref{C2} dealing with the first 
and Lemma \ref{L1} with the second factor in the singular
products. We will use Corollary \ref{C2}, 
in order to pass from the definition of $v \diam  f$ and $v\diam\partial_1^2v$
to the definition of $u\diam f$ and $b\diam\partial_1^2v$, respectively (since the distribution
$\partial_1^2v$ plays a role very similar
to $f$, the lemma and the corollary are formulated in the notation of the former case).
We will then use Lemma \ref{L1} to pass from $b\diam\partial_1^2v$ to $b\diam\partial_1^2u$.

\medskip

%Lemmas \ref{L2a}, \ref{L2}, and \ref{L1}, 
These upcoming statements reveal a clear hierarchy of norms and measures of size: 
\begin{itemize}
\item Functions $u$ are measured 
in terms of the H\"older semi-norm $[u]_\alpha$ (the supremum norm $\|\sigma\|$ of a function
$\sigma$ only intervenes in scaling-wise suboptimal estimates like (\ref{x44}) that rely on the
periodicity or the constraint $T\le 1$ providing a large-scale cut-off, otherwise just
as part of the product $\|\sigma\|[a]_\alpha$ with the H\"older norm of $a$),
\item 
distributions are measured in the $C^{\alpha-2}$-norm %$\sup_{T\le 1}$ $(T^\frac{1}{4})^{2-\alpha}$ 
$\|f\|_{\alpha-2}$ (defined in \eqref{c.1}),
%see Step \ref{LA1S1} in the proof of Lemma \ref{LA1} for this equivalence of norms,
\item 
commutators $[u,(\cdot)_T]\diam f$ 
are measured on level $2\alpha-2<0$ via 
\begin{equation}\label{e:2alphadef}
\|[u,(\cdot)]\diam f\|_{2\alpha-2}:= \sup_{T\le 1}(T^\frac{1}{4})^{2-2\alpha}\|[u,(\cdot)_T]\diam f\|,
\end{equation}
and 
\item differences $[u,(\cdot)_T]\diam f-[v,(\cdot)_T]\diam f$ of commutators, like in case of the
rough path expression (\ref{i-17}) divided by $(t-s)$, are measured on level $3\alpha-2>0$ via
\begin{equation}\label{e:difference3alpha}
\|[u,(\cdot)]\diam f-[v,(\cdot)]\diam f\|_{3\alpha-2} := \sup_{T\le 1}(T^\frac{1}{4})^{2-3\alpha}\|[u,(\cdot)_T]\diam f-[v,(\cdot)_T]\diam f\|,
\end{equation}
see (\ref{wi84}) of Lemma \ref{L2a}.
\end{itemize}

\medskip

Equipped with this dictionary, Corollary~\ref{C2} and Lemma \ref{L1} can be seen to be very close to
\cite[Theorem 1]{GubinelliControlling}; in particular, (\ref{wi84}) in Lemma \ref{L2a} 
is very close to (28) in
\cite[Corollary 3]{GubinelliControlling}. The major difference is the multi-dimensional extension through
(\ref{v38}). A minor difference coming from the parabolic nature is the appearance of the
commutator $[x_1,(\cdot)_T]f$, which however is regular, cf. Lemma \ref{LA2}. 
A further minor difference arises from the $a_0$-dependence of the model $v$ and the related appearance
of the function $a$, which necessitates control of $\frac{\partial}{\partial a_0}$-derivatives
of the functions and the commutators and manifests itself via the evaluation operator $E$. 
However, these minor differences can be embedded into the more general form of the upcoming Lemma \ref{L2a}.

\medskip

\begin{lemma}\label{L2a}
Let $\frac{2}{3}<\alpha<1$. Suppose we have a family of functions $\{v(\cdot,x)\}_x$ of class $C^\alpha$,
parameterized by points $x$, a distribution $f$, and a family of distributions
$\{v(\cdot,x)\diam f\}_x$, both of class $C^{\alpha-2}$, satisfying
%%%%
%%%% Here I could refer to \eqref{c.1} mit N_0 anstelle von N_1
%%%%
%
\begin{align}
[v(\cdot,x)-v(\cdot,x')]_\alpha&\le Nd^\alpha(x,x'),\label{wi80}\\
%\sup_{T\le 1}(T^\frac{1}{4})^{2-\alpha}\|f_T\|&\le N_1,\label{wi78}\\
\|f\|_{\alpha-2}&\le N_1,\label{wi78}\\
%\sup_{T\le 1}(T^\frac{1}{4})^{2-2\alpha}
%\|[v(\cdot,x),(\cdot)_T]\diam f
%-[v(\cdot,x'),(\cdot)_T]\diam f\|&\le NN_1d^\alpha(x,x')\label{wi79}
\|[v(\cdot,x),(\cdot)]\diam f
-[v(\cdot,x'),(\cdot)]\diam f\|_{2\alpha-2}&\le NN_1d^\alpha(x,x')\label{wi79}
\end{align}
for all pairs of points $x,x'$ and for some constants $N,N_1$\footnote{in \eqref{wi79} the $2\alpha-2$ semi-norm
of the difference of commutators is defined as \eqref{e:difference3alpha} with $3\alpha-2$ replaced by $2\alpha-2$.}.
Suppose we are given a function $u$ such
that
\begin{align}\label{wi81}
|(u(y)-u(x))-(v(y,x)-v(x,x))
-\nu(x)(y-x)_1|\le M d^{2\alpha}(y,x)
\end{align}
for all pairs of points $y,x$ for some constant $M$ and some function $\nu$.
Then there exists a unique distribution $u\diam f$ such that
\begin{align}\label{wi84}
% \sup_{T\le 1}(T^\frac{1}{4})^{2-3\alpha}
 \|[u,(\cdot)]\diam f 
-\Ed [v,(\cdot)]\diam f-\nu[x_1,(\cdot)]f\|_{3\alpha-2}\lesssim (M+N)N_1,
\end{align}
where $\Ed$ stands for the evaluation of the continuous function 
$(x,y)\mapsto([v(\cdot,x),(\cdot)_T]\diam f)(y)$ on the diagonal $y=x$.
\medskip

If moreover all functions and distributions are space-time periodic and we use the constant $N$ to also estimate
the lower-order expressions
\begin{align}
[v(\cdot,x)]_\alpha&\le N,\label{wi75}\\
%\sup_{T\le 1}(T^\frac{1}{4})^{2-2\alpha}
\|[v(\cdot,x),(\cdot)]\diam f\|_{2\alpha-2}&\le NN_1\label{wi85}
\end{align}
for all points $x$ then also
\begin{align}\label{wi83}
%\sup_{T\le 1}(T^\frac{1}{4})^{2-2\alpha}
\|[u,(\cdot)]\diam f\|_{2\alpha-2}&\le (M+N)N_1.
\end{align}
\end{lemma}

\medskip
%%%%%%%%%%%%%%%%%%%%%%%%%%%%%%%%%%%%%%%%%%%%%%%%%%%%%%%%%%%%%%%%%%%%%%%%%%%%%%%%%%%%%%%%%%%%%%%%%%%%%%%%%%5

Equipped with Lemma \ref{L2a}, the upcoming corollary specifies the form of the model.
The general form of Lemma \ref{L2a} is in
particular convenient for part iii), where the Lipschitz continuity
of the product $\sigma\diam f$ in terms of the off-line product $v\diam f$ and
the modulating property (both constant and modulating functions) is established.

To shorten some of the formulas, from now on we add some more indices to the (semi-) norms
 referring to parameter derivatives
with respect to $a_0$ and $a_0'$: If $| \cdot |$ is a semi-norm and 
if $u$ depends on a parameter $a_0$ we write
\begin{equation}\label{conv1}
|u|_{n} := \sup_{a_0 \in [\lambda, \frac{1}{\lambda}]} \max_{i=0, \ldots, n} \Big| \Big(\frac{\partial}{\partial a_0}\Big)^i  u(a_0)\Big|
\end{equation}
and if $u$ depends on two parameters $a_0$ and $a_0'$ we write
\begin{equation}\label{conv2}
|u|_{n,m}:= \sup_{a_0, a_0' \in [\lambda, \frac{1}{\lambda}]} \max_{i=0, \ldots, n} \max_{j=0,\ldots,m}  \Big| \Big(\frac{\partial}{\partial a_0}\Big)^i  \Big(\frac{\partial}{\partial a_0}\Big)^j  u(a_0, a_0') \Big|.
\end{equation}

%
%, i.e. we will write
%$[  v  ]_{\alpha,1} := \sup_{a_0} \big( [  v  ]_{\alpha} + [ \frac{\partial}{\partial a_0} v  ]_{\alpha}\big)$,
%$\|  f \|_{\alpha-2,2} :=\sup_{a_0'}\big( \|   f \|_{\alpha-2} + \| \frac{\partial}{\partial a_0'} f \|_{\alpha-2} + \| \frac{\partial^2}{\partial a_0'^2} f \|_{\alpha-2}$,
%%
%$\|  [v,(\cdot)]  \diam f \Big\|_{ 2\alpha-2,1,2}  := \sup_{a_0,a_0'}\sum_{i=0,1\; j=0,1,2} \big\| 
%	 \big(\frac{\partial}{\partial a_0}\big)^i 
%	\big( \frac{\partial}{\partial a_0'} \big)^j
%	 [v,(\cdot)]\diam f \big\|_{ 2\alpha-2} $.

\begin{corollary}\label{C2}
i) Let $\{v(\cdot,a_0)\}_{a_0}$ be a family of functions and let
$\{f(\cdot,a_0')\}_{a_0'}$,  $\{v(\cdot,a_0)\diam f(\cdot,a_0')\}_{a_0,a_0'}$
be two families of distributions satisfying 
%%%%%%%%%%
%%%%%%%%%% Hier gibt es auch eine f Annahme, aber die ist nicht genau wie vorher
%
\begin{align}
[  v  ]_{\alpha,2}  &\le N_0,\label{L0.1}\\
\|  f \|_{\alpha-2,2}&\le N_1,\label{wj36}\\
\|  [v,(\cdot)]  \diam f \|_{ 2\alpha-2,1,2} &\le N_1N_0
\label{wj41}
\end{align}
for some constants $N_0$ and $N_1$. 
If $u$ is modelled after $v$ according to the $\alpha$-H\"older functions
$a$ and $\sigma$ with constant $M$ and $\nu$ as in (\ref{1.3}), 
then there exists a unique family of distributions $\{  u\diam f\}_{a_0, a_0'} $ such that
\begin{align}\label{w14}
\lim_{T\downarrow0}  \big\|  [u,(\cdot)_T]\diam f
-\sigma E[v,(\cdot)_T]\diam f-\nu [x_1,(\cdot)_T]f \big\|=0,
\end{align}
where $E$ evaluates a function of $(x,a_0)$ at $(x,a(x))$. Furthermore, in case of 
\begin{align}\label{wk58}
[\sigma]_\alpha\le 1,\;[a]_\alpha\le 1\quad\mbox{and}\quad\|\sigma\|\le 1
\end{align}
and when all functions are space-time periodic we have the sub-optimal estimate
\begin{align}\label{wj40}
%\sup_{T\le 1}(T^\frac{1}{4})^{2-2\alpha}
%\sup_{a_0'}
\big\| 
	[u,(\cdot)]\diam f \big\|_{2\alpha-2,2}
\lesssim N_1(M+N_0).
\end{align}
ii) Let $\{v(\cdot,a_0)\}_{a_0}$, $\{f_j(\cdot,a_0')\}_{a_0'}$, and 
$\{v(\cdot,a_0)\diam f_j(\cdot,a_0')\}_{a_0,a_0'}$, $j=0,1$, be as in i) 
and suppose in addition
\begin{align}
\|
	f_{1}-f_{0}
	\|_{\alpha-2,1} 
&\le \del N_1,\label{wj52}\\
	\| 
	[v,(\cdot)]\diam f_1 -[v,(\cdot)]\diam f_0
	\|_{2\alpha-2,1,1}
&\le \del N_1 N_0\label{wj48}
\end{align}
for some constant $\del N_1$. Then for $u$ and $u \diam f_i$ as in i) we have
\begin{align}\label{wk35}
%\sup_{T\le 1}&(T^\frac{1}{4})^{2-2\alpha}
\|
	[u,(\cdot)]\diam f_1 -[u,(\cdot)]\diam f_0
	\|_{2\alpha-2,1}\lesssim \del N_1(M+N_0).
\end{align}
iii) Let the two families of functions $\{v_i(\cdot,a_0)\}_{a_0}$, $i=0,1$, 
and the three families of distributions $\{f(\cdot,a_0')\}_{a_0'}$, 
$\{v_i(\cdot,a_0)\diam f(\cdot,a_0')\}_{a_0,a_0'}$
satisfy \eqref{wj36} and in addition  
%
%in addition (\ref{wi20}), (\ref{wi21}).
%Suppose we have in addition
%
\begin{align}
[
	v_i
	]_{\alpha,2} 
	&\le N_0,\label{wi20}\\
[
	v_1-v_0
]_{\alpha,1}&\le \del N_0,\label{wi21}\\
\| 
%	\Big\{1,\frac{\partial}{\partial a_0},\frac{\partial^2}{\partial a_0^2} \Big\}
%	\Big\{1,\frac{\partial}{\partial a_0'} \Big\} 
[v_i,(\cdot)]\diam f 
\|_{2\alpha-2,2,1}
&\le N_1N_0,\label{wj50}\\
%
%\sup_{T\le 1}(T^\frac{1}{4})^{2-2\alpha}
\|
%	\Big\{ 1,\frac{\partial}{\partial a_0} \Big\}
%	\Big\{ 1,\frac{\partial}{\partial a_0'} \Big\} 
	[v_1,(\cdot)]\diam f -[v_0,(\cdot)]\diam f
\|_{2\alpha-2,1,1}
&\le N_1\del N_0.\label{wj49}
\end{align}
Let $u_i$ be two functions like in part i) and let $u_i \diam f$ be as constructed there.  Suppose that $u_1-u_0$ is modelled after $(v_1,v_0)$
according to $(a_1,a_0)$ and $(\sigma_1,-\sigma_0)$ with constant $\del M$.
Then we have
\begin{align}\label{wk36}
\notag
%&\sup_{T\le 1}(T^\frac{1}{4})^{2-2\alpha}\sup_{a_0'}
& \|
%\{1,\frac{\partial}{\partial a_0'}\}(
[u_1,(\cdot)]\diam f-[u_0,(\cdot)]\diam f
%) 
\|_{2\alpha-2,1} \\
& \lesssim N_1\big(\del M 
+N_0([\sigma_1-\sigma_0]_\alpha
+\|\sigma_1-\sigma_0\|
+[a_1-a_0]_\alpha+\|a_1-a_0\|)+\del N_0\big).
\end{align}
\end{corollary}

\medskip
%%%%%%%%%%%%%%%%%%%%%%%%%%%%%%%%%%%%%%%%%%%%%%%%%%%%%%%%%%%%%%%%%%%%%%%%%%%%%%%%%

We now turn to Lemma \ref{L1} that deals with the second factor in $a\diam\partial_1^2u$.
The reason why we consider several functions $v_1,\cdots,v_I$ in Lemma \ref{L1}
instead of a single one for our scalar PDE 
is that this seems necessary when establishing the contraction property for Proposition \ref{P};
because of the $a_0$-dependence, it turns out that we need not just $I=2$ but in fact $I=3$,
cf. Corollary \ref{C1}.

\begin{lemma}\label{L1}
Let $\frac{2}{3}<\alpha<1$ and $I \in \N$. We are given a function $b$, $I$ families of
functions $\{v_1(\cdot,a_0),$ $\cdots,v_I(\cdot,a_0)\}_{a_0}$, and $I$ families of 
distributions $\{b\diam\partial_1^2v_1(\cdot,a_0),$ 
$\cdots,b\diam\partial_1^2v_I(\cdot,a_0)\}_{a_0}$ with
\begin{align}
[v_i]_{\alpha,1}&\le N_i,\label{wi5}\\
\| [b,(\cdot)]\diam\partial_1^2 v_i\|_{2\alpha-2,1} &\le N_0N_i\label{1.5}
\end{align}
for some constants $N_0,\cdots,N_I$.
Let the function $u$ be modelled after $(v_1,\cdots,v_I)$ according to the $\alpha$-H\"older functions
$a$ and $(\sigma_1,\cdots,\sigma_I)$ with constant $M$, cf. Definition \ref{D1}.
Then there exists a unique distribution $b\diam\partial_1^2u$ such that 
\begin{align}
\lim_{T\downarrow0}\|[b,(\cdot)_T]\diam\partial_1^2u-\sigma_i E[b,(\cdot)_T]\diam\partial_1^2v_i\|
=0,\label{1.9}
\end{align}
where $E$ denotes the operator that evaluates a function in two variables $(x,a_0)$ at
$(x,a(x))$.
Moreover, provided $[a]_\alpha\le 1$, we have the sub-optimal estimate
\begin{align}
%\sup_{T\le 1}(T^\frac{1}{4})^{2-2\alpha}
\|[b,(\cdot)]\diam\partial_1^2u\|_{2\alpha-2}
\lesssim [b]_\alpha M+N_0N_i([\sigma_i]_{\alpha}+\|\sigma_i\|).
\label{x44}
\end{align}
\end{lemma}

\bigskip

%%%%%%%%%%%%%%%%%%%%%%%%%%%%%%%%%%%%%%%%%%%%%%%%%%%%%%%%%%%%%%%%%%%%%%%%%%%%%%%%%%%
The following lemma is the only place where we use the PDE. It might be seen as an extension
of Schauder theory in the sense that it compares, on the level of $C^{2\alpha}$, the solution $u$ of
a variable-coefficient equation $\partial_2u-a\diam\partial_1^2u=\sigma\diam f$ to
the solutions of the corresponding constant-coefficient equation (\ref{v53}), by saying
that $u$ is modelled after $v$ according to $a$ and $\sigma$. To this purpose we apply $(\cdot)_T$
to the equation and rearrange to
\begin{align*}
\partial_2u_T-P(a\partial_1^2u_T+\sigma f_T)=-P\big([a,(\cdot)_T]\diam\partial_1^2u+[\sigma,(\cdot)_T]\diam f\big).
\end{align*}
Since the previous lemmas estimate the commutators on the right hand side,
we will right away assume that the left hand side is estimated accordingly, cf. (\ref{y39}).
Working with the commutator of multiplication with a coefficient $a$ and convolution 
is reminiscent of the DiPerna-Lions theory, which however deals with a transport
instead of a parabolic equation with a rough coefficient, that is $\partial_2u-a\partial_1u$
instead of $\partial_2u-a\partial_1^2u$.
In our proof, we follow the approach to classical Schauder theory of  Safonov, 
 \cite{Krylov}, in particular Section 8.6. This approach avoids the use of kernels.

\begin{lemma}\label{L3}
Let $\frac{1}{2}<\alpha<1$ and suppose all functions and distributions are periodic. 
We are given $I$ families of distributions $\{f_1(\cdot,a_0),\cdots,$ $f_I(\cdot,a_0)\}_{a_0}$ with
%
%
% f Annahme etwas anders 
%
\begin{align}\label{3.31}
 \|f_{i}\|_{\alpha-2,1}  \le N_i
\end{align}
for some constants $N_1,\cdots,N_I$.
For $a_0\in[\lambda,\frac{1}{\lambda}]$ we denote by $v_i(\cdot,a_0)$ the function
of vanishing mean solving
\begin{align}\label{v53}
(\partial_2-a_0\partial_1^2)v_i(\cdot,a_0)=Pf_i(\cdot,a_0)\quad\mbox{distributionally}.
\end{align}
We are also given a function $u$, modelled after $(v_1,\cdots,v_I)$ 
according to some functions $a\in[\lambda,\frac{1}{\lambda}]$ and 
$(\sigma_1,\cdots,\sigma_I)$.
We assume that 
%$u$ approximately satisfies the PDE 
%$\partial_2u-Pa\partial_1^2u=P\sigma_i Ef_i$ in the sense of 
%
\begin{align}\label{y39}
\sup_{T\le 1}(T^\frac{1}{4})^{2-2\alpha}
\|\partial_2u_T-P(a\partial_1^2u_T+\sigma_i Ef_{iT})\|\le N^2
\end{align}
for some constant $N$, where $E$ is defined  in Lemma \ref{L1}.
Then we have for the modelling and the H\"older constant of $u$
\begin{align}
M&\lesssim N^2+[a]_\alpha M+N_i([\sigma_i]_{\alpha}+\|\sigma_i\|[a]_\alpha),\label{3.15}\\
[u]_\alpha&\lesssim M+N_i\|\sigma_i\|.\label{wk01}
\end{align}
\end{lemma}

\medskip

%%%%%%%%%%%%%%%%%%%%%%%%%%%%%%%%%%%%%%%%%%%%%%%%%%%%%%%%%%%%%%%%%%%%%%%%%%%%%%%%%%%%%%%%%%%%%%%%%%%%%%%
%%%%%%%%%%%%%%%%%%%%%%%%%%%%%%%%%%%%%%%%%%%%%%%%%%%%%%%%%%%%%%%%%%%%%%%%%%%%%%%%%%%%%%%%%%%%%%%%%%%%%%%

In the upcoming Corollary \ref{C1}, we  combine Corollary~\ref{C2} on the product $\sigma \diam f$, Lemma \ref{L1} on the product $a\diam\partial_1^2u$ and Lemma \ref{L3}
to obtain an {\it a priori} estimate on the modelling and H\"older constants. 
The use of the ``infinitesimal'' part ii) of this corollary will
be explained in the discussion of Proposition \ref{P}.

\begin{corollary}\label{C1}
Let $\frac{2}{3}<\alpha<1$. 
\medskip

i) Suppose we are given two functions $\sigma$ and $a$, two distributions $f$ and $\sigma\diam f$,
and a family of distributions $\{a\diam\partial_1^2v(\cdot,a_0)\}_{a_0}$ with
\begin{align}
[\sigma]_\alpha+[a]_\alpha&\le N,\label{wk59}\\
%\sup_{T\le 1}(T^\frac{1}{4})^{2-\alpha}
\|f\|_{\alpha-2}&\le N_0,\label{wj58}\\
%\sup_{T\le 1}(T^\frac{1}{4})^{2-2\alpha}
\|[\sigma,(\cdot)]\diam f\|_{2\alpha-2}&\le N N_0,\label{wj56}\\
%\sup_{T\le 1}(T^\frac{1}{4})^{2-2\alpha}\sup_{a_0}
\|
%\{1,\frac{\partial}{\partial a_0},
%\frac{\partial^2}{\partial a_0^2}\}
[a,(\cdot)]\diam\partial_1^2v\|_{2\alpha-2,2}&\le N N_0,\label{wj57}
\end{align}
for some constants $N_0$ and $N$, where $v(\cdot,a_0)$ denotes the mean-free solution of \eqref{e:v-eq},
%%
%\begin{align}\label{wj64}
%(\partial_2-a_0\partial_1^2)v(\cdot,a_0)=Pf,
%\end{align}
%
and satisfying the constraints
\begin{align}\label{wj59}
\sigma\in[-1,1],\;a\in[\lambda,\frac{1}{\lambda}],\;
[\sigma]_\alpha\le 1,\;[a]_\alpha\ll 1.
\end{align}
Then if a function $u$ is modelled after $v$ according to $a$ and $\sigma$ with
\begin{align}\label{wj55}
\partial_2u-P(a\diam\partial_1^2u+\sigma\diam f)=0 \qquad \text{distributionally}
\end{align}
we have for the modelling and H\"older constants
\begin{align}
M&\lesssim N_0N,\label{wj63}\\
[u]_\alpha&\lesssim N_0(N+1).\label{wk02}
\end{align}
ii) In addition, suppose we are given two functions $\del\sigma$ and $\del a$, three distributions $\del f$,
$\sigma\diam \del f$, and $\del\sigma\diam f$,
and two families of distributions $\{a\diam\partial_1^2\del v(\cdot,a_0)\}_{a_0}$ and
$\{\del a\diam\partial_1^2v(\cdot,a_0)\}_{a_0}$ with
\begin{align}
[\del\sigma]_\alpha+\|\del\sigma\|+[\del a]_\alpha+\|\del a\|&\le \del N,\label{wk60}\\
%\sup_{T\le 1}(T^\frac{1}{4})^{2-\alpha}
\|\del f\|_{\alpha-2} &\le \del N_0,\label{wj78}\\
%\sup_{T\le 1}(T^\frac{1}{4})^{2-2\alpha}
\|[\sigma,(\cdot)]\diam \del f\|_{2\alpha-2}
&\le N \del N_0,\label{wj67}\\
%\sup_{T\le 1}(T^\frac{1}{4})^{2-2\alpha}
\|[\del\sigma,(\cdot)]\diam f\|_{2\alpha-2}&\le \del N N_0,\label{wj68}\\
%\sup_{T\le 1}(T^\frac{1}{4})^{2-2\alpha}\sup_{a_0}\|\{1,\frac{\partial}{\partial a_0}\}
\| [a,(\cdot)]\diam\partial_1^2 \del v\|_{2\alpha-2,1}&\le N \del N_0,
\label{wj69}\\
%\sup_{T\le 1}(T^\frac{1}{4})^{2-2\alpha}\sup_{a_0}
\|
%\{1,\frac{\partial}{\partial a_0}\}
[\del a,(\cdot)]\diam\partial_1^2v\|_{2\alpha-2,1}&\le \del N N_0\label{wj92}
\end{align}
for some constants $\del N_0,\del N$ and where $\del v(\cdot, a_0)$ is the mean-free solution of
\begin{align}\label{wj65}
(\partial_2-a_0\partial_1^2)\del v(\cdot, a_0) =P\del f \qquad \text{distributionally}.
\end{align}
Then if a function $\del u$ is modelled after $(v,\frac{\partial v}{\partial a_0},\del v)$ according
to $a$ and $(\del\sigma,\sigma\del a,\sigma)$ with
\begin{align}\label{wj66}
\partial_2\del u-P(a\diam\partial_1^2\del u+\del a\diam\partial_1^2u+\sigma\diam\del f
+\del\sigma\diam f)=0
\end{align}
then we have for the modelling and H\"older constants
\begin{align}
\del M         &\lesssim N_0\del N+\del N_0N\quad\mbox{provided}\;N\le 1,\label{wj80}\\
[\del u]_\alpha&\lesssim N_0\del N+\del N_0 \quad\mbox{provided}\;N\le 1.\label{wk03}
\end{align}
\end{corollary}

\medskip

%%%%%%%%%%%%%%%%%%%%%%%%%%%%%%%%%%%%%%%%%%%%%%%%%%%%%%%%%%%%%%%%%%%%%%%%%%%%%%%%%%%%%%%%%%%%%%%%%%%%%%
%%%%%%%%%%%%%%%%%%%%%%%%%%%%%%%%%%%%%%%%%%%%%%%%%%%%%%%%%%%%%%%%%%%%%%%%%%%%%%%%%%%%%%%%%%%%%%%%%%%%%%%
%%%%%%%%%%%%%%%%%%%%%%%%%%%%%%%%%%%%%%%%%%%%%%%%%%%%%%%%%%%%%%%%%%%%%%%%%%%%%%%%%%%%%%%%%%%%%%%%%%%%%%%

The following Proposition \ref{P} may be seen as the main contribution of this paper.
It establishes a solution theory for the linear equation
$\partial_2u-P(a\diam\partial_1^2u+\sigma\diam f)=0$ for given driver $f$ (a distribution) and
given coefficients $\sigma$ and $a$. Because of the roughness
of $f$, it does not only require a definition of $\sigma\diam f$ but also of $a\diam\partial_1^2v$,
%where $v(\cdot,a_0)$ solves $\partial_2v-a_0\partial_1^2v=Pf$, 
so that when $u$ is modelled after $v$ according to $a$ and $\sigma$,
also $a\diam\partial_1^2u$ may be given a sense by Lemma \ref{L1}. The most subtle point is to establish
Lipschitz continuity of $u$ in the data $(a,a\diam\partial_1^2v)$.
This involves considering differences of solutions and quantifying
\begin{align}\label{wi11bis}
u_1-u_0\;\mbox{is modelled after}\;(v_1,v_0)\nonumber\\
\mbox{according to}\;(a_1,a_0)\;\mbox{and}\;(\sigma_1,-\sigma_0).
\end{align}
When quantifying differences of solutions, variable coefficients require a somewhat different strategy
compared to constant coefficients, as we shall explain now.
The modelledness \eqref{wi11bis} has to come from the PDE, that is, Lemma \ref{L3}. The naive approach is to consider
the difference of the PDE for two given pairs of data $(\sigma_i,a_i,f_i)$, $i=0,1$, (plus the products),
and to rearrange as follows
\begin{align}\label{wi10}
\lefteqn{\partial_2(u_1-u_0)-P(a_0\diam\partial_1^2u_1-a_0\diam\partial_1^2u_0)}\nonumber\\
&=P\big(\sigma_1\diam f_1-\sigma_0\diam f_0+(a_1\diam\partial_1^2u_1-a_0\diam\partial_1^2u_1)\big),
\end{align}
which already means breaking the permutation symmetry in $i=0,1$ and therefore does not bode well.
By the modelledness of $u_1$ we expect that for the purpose
of Lemma \ref{L3}, we may replace $u_1$ by $v_1$ on the right hand side of (\ref{wi10}), leading to
\begin{align}\nonumber
\lefteqn{\partial_2(u_1-u_0)-P(a_0\diam\partial_1^2u_1-a_0\diam\partial_1^2u_0)}\nonumber\\
&\approx P\big(\sigma_1\diam f_1-\sigma_0\diam f_0+\sigma_1(E_1 a_1\diam\partial_1^2v_1- E_1 a_0\diam\partial_1^2v_1)\big).
\end{align}
In view of Lemma \ref{L3} and the discussion preceding it, this suggests that we obtain
\begin{align}\label{wi11}
u_1-u_0\;\mbox{is modelled after}\;(v_1,v_0,(\partial_2-a_0\partial_1^2)^{-1}PE_1\partial_1^2v)\nonumber\\
\mbox{according to}\;a_0\;\mbox{and}\;(\sigma_1,-\sigma_0,\sigma_1(a_1-a_0)),
\end{align}
which is {\it not} the desired (\ref{wi11bis}) unless $a_1=a_0$. Instead, our strategy will be to construct a {\it curve} $\{u_s\}_{s\in[0,1]}$
interpolating between $u_0$ and $u_1$. For this, we interpolate the data linearly, that is,
$f_s:=sf_1+(1-s)f_0$, $\sigma_s:=s\sigma_1+(1-s)\sigma_0$, and $a_s:=sa_1+(1-s)a_0$, and solve
\begin{align}\label{wi14}
\partial_2u_s-P(a_s\diam\partial_1^2u_s+\sigma_s\diam f_s)=0.
\end{align}
Provided we interpolate the products bi-linearly, that is,
\begin{align}\label{wi13}
\sigma_s\diam f_s:=s^2\sigma_1\diam f_1+s(1-s)(\sigma_1\diam f_0+\sigma_0\diam f_1)+(1-s)^2\sigma_0\diam f_0
\end{align}
and the same definition for $a_s\diam\partial_1^2v_s$,
Leibniz's rule for $\sigma_s\diam f_s$ holds, and we expect it to hold for $a_s\diam\partial_1^2u_s$ so that
differentiation of (\ref{wi14}) gives
\begin{align*}
\partial_2\partial_su-P(a_s\diam\partial_1^2\partial_su)=P(\partial_sa\diam\partial_1^2u_s+\partial_s\sigma\diam f_s
+\sigma_s\diam\partial_sf),
\end{align*}
where we write $\partial_su$ as short hand for $\partial_s u_s$ and the same for $a$, $\sigma$ and $f$. 
In view of (\ref{wi14}) we approximate the right hand side by 
\begin{align*}
\partial_2\partial_su-P(a_s\diam\partial_1^2\partial_su)\approx P(\sigma_sE_s\partial_sa\diam\partial_1^2v_s
+\partial_s\sigma\diam f_s
+\sigma_s\diam\partial_sf),
\end{align*}
with $v_s=sv_1+(1-s)v_0$. It is this form that motivates part ii) of Corollary \ref{C1}.
Noting that $(\partial_2-a_0\partial_1^2)\frac{\partial v_s}{\partial a_0}=\partial_1^2v_s$ 
\footnote{Here we use the symbol $a_0$ with two different meanings: as the concrete coefficient field $a_0$
and as an abstract parameter in $v_s$. It will always be clear from the context which of these interpretations is meant.  }
 we obtain
\begin{align}\label{wi15}
\partial_su\;\mbox{is modelled after}\;(v_s,\frac{\partial v_s}{\partial a_0},\partial_sv)\nonumber\\
\mbox{according to}\;a_s\;\mbox{and}\;(\partial_s\sigma,\sigma_s\partial_sa,\sigma_s),
\end{align}
which compares favorably to (\ref{wi11}). Using Leibniz's rule once more, but this time in the classical form of
\begin{align*}
\lefteqn{\frac{\partial}{\partial s}\big(\sigma_s(x)v_s(y,a_s(x))\big)
=(\partial_s\sigma)(x) v_s(y,a_s(x))}\nonumber\\
&+(\sigma_s\partial_sa)(x)\frac{\partial v_s}{\partial a_0}(y,a_s(x))
+\sigma_s(x)\partial_sv(y,a_s(x)),
\end{align*}
and integrating (\ref{wi15}) in $s\in[0,1]$ yields  the desired (\ref{wi11bis}).
We note that this strategy differs from \cite{GubinelliControlling} even in case when $a$ is constant:
When passing from the modelledness of $u_1-u_0$ to the modelledness of $\sigma(u_1)-\sigma(u_0)$,
the argument in \cite[Proposition 4]{GubinelliControlling} uses the linear interpolation $u_s=su_1+(1-s)u_0$
(as we do in Lemma \ref{L0}),
which implicitly amounts to the interpolation $\sigma_s\diam f_s=s\sigma_1\diam f_1+(1-s)\sigma_0\diam f_0$,
as opposed to (\ref{wi13}).

\begin{proposition}\label{P}
Let $\frac{2}{3}<\alpha<1$. 

\medskip
i) Suppose we are given two functions $\sigma$ and $a$, two distributions $f$ and $\sigma\diam f$,
and a family of distributions $\{a\diam\partial_1^2v(\cdot,a_0)\}_{a_0}$ satisfying \eqref{wk59} --\eqref{wj59}.
%%
%%%
%%\begin{align}
%%[\sigma]_\alpha+[a]_\alpha&\le N,\label{wk51}\\
%%%\sup_{T\le 1}(T^\frac{1}{4})^{2-\alpha}
%%\|f\|_{\alpha-2}&\le N_0,\label{wj95}\\
%%%\sup_{T\le 1}(T^\frac{1}{4})^{2-2\alpha}
%%\|[\sigma,(\cdot)]\diam f\|_{2\alpha-2}&\le N_0 N,\label{wj97}\\
%%%\sup_{T\le 1}(T^\frac{1}{4})^{2-2\alpha}
%%%\sup_{a_0}
%%\|\
%%%{1,\frac{\partial}{\partial a_0},\frac{\partial^2}{\partial a_0^2}\}
%%[a,(\cdot)]
%%\diam\partial_1^2v\|_{2\alpha-2,2}
%%&\le N_0 N,\label{wj98}
%%\end{align}
%%%
%%%for some constants $N_0,\;N$, where $v(\cdot,a_0)$ denotes the mean-free solution of $(\partial_2-a_0\partial_1^2)v(\cdot,a_0)=Pf$.
%%%Let the functions $\sigma$ and $a$ satisfy
%
%\begin{align}\label{wj96}
%%\sigma\in[-1,1],\;a\in[\lambda,\frac{1}{\lambda}],\;[\sigma]_\alpha\le 1,\;
%[a]_\alpha\ll 1.
%\end{align}
%
Then there exists a unique mean-free function $u$ modelled after $v$ according to
$a$ and $\sigma$ and such that
\begin{align}\label{wk07}
\partial_2u-P(a\diam\partial_1^2u+\sigma\diam f)=0.
\end{align}
The modelling and H\"older constants are estimated as follows
\begin{align}
M&\lesssim N_0N,\label{wk17}\\
[u]_\alpha&\lesssim N_0(N+1).\label{wk18}
\end{align}

\medskip
ii) Suppose we are given  functions $\sigma_i$ and $a_i$, $i=0,1$, 
distributions $f_i$ and $\sigma_i\diam f_j$, $j=0,1$, and  families
of distributions 
$\{a_i\diam\partial_1^2v_j(\cdot,a_0)\}_{a_0}$, where $v_i(\cdot,a_0)$ is the mean-free solution
of \eqref{e:v-eq} corresponding to $f_i$,
% $(\partial_2-a_0\partial_1^2)v_i$ $=Pf_i$
% 
  satisfying the assumption 
(\ref{wk59})--
%(\ref{wj58}), (\ref{wj56}), and 
(\ref{wj57})
%the two latter 
with cross terms, that is,
\begin{align}
%\sup_{T\le 1}(T^\frac{1}{4})^{2-2\alpha}
\|[\sigma_i,(\cdot)]\diam f_j\|_{2\alpha-2} &\le N_0 N,\label{wj97bis}\\
%\sup_{T\le 1}(T^\frac{1}{4})^{2-2\alpha}\sup_{a_0}
\|
%\{1,\frac{\partial}{\partial a_0},\frac{\partial^2}{\partial a_0^2}\}
[a_i,(\cdot)]
\diam\partial_1^2v_j\|_{2\alpha-2,2}
&\le N_0 N\label{wj98bis}
\end{align}
and \eqref{wj59}.
% and (\ref{wj96}). 
We measure the distance of $(f_1,\sigma_1, a_1)$ to $(f_0, \sigma_0, a_0)$ in terms of the constants $\del N_0$
and $\del N$ with
\begin{align}
[\sigma_1-\sigma_0]_\alpha+\|\sigma_1-\sigma_0\|+[a_1-a_0]_\alpha+\|a_1-a_0\|&\le\del N,\label{wk50}\\
%\sup_{T\le 1}(T^\frac{1}{4})^{2-\alpha}
\|f_{1}-f_{0}\|_{\alpha-2} &\le\del N_0,\label{wk24}\\
%\sup_{T\le 1}(T^\frac{1}{4})^{2-2\alpha}
\|[\sigma_i,(\cdot)]\diam f_1
-[\sigma_i,(\cdot)]\diam f_0\|_{2\alpha-2}
&\le N \del N_0,\label{wk13}\\
%
%\sup_{T\le 1}(T^\frac{1}{4})^{2-2\alpha}
\|[\sigma_1,(\cdot)]\diam f_j
-[\sigma_0,(\cdot)]\diam f_j\|_{2\alpha-2}
&\le\del N N_0,\label{wk25}\\
%
%\sup_{T\le 1}(T^\frac{1}{4})^{2-2\alpha}
%\sup_{a_0}
\|
%\{1,\frac{\partial}{\partial a_0}\}
[a_i,(\cdot)]\diam\partial_1^2v_1
-[a_i,(\cdot)]\diam\partial_1^2v_0\|_{2\alpha-2,1}
&\le N \del N_0,
\label{wk14}\\
%
%\sup_{T\le 1}(T^\frac{1}{4})^{2-2\alpha}\sup_{a_0}
\|
%\{1,\frac{\partial}{\partial a_0}\}
%\big(
[a_1,(\cdot)]\diam\partial_1^2v_j
-[a_0,(\cdot)]\diam\partial_1^2v_j\|_{2\alpha-2,1}&\le\del N N_0.\label{w26}
\end{align}
 Let $u_i$ denote the corresponding solutions ensured by part i).
Then $u_1-u_0$ is modelled after $(v_1,v_0)$ according to $(a_1,a_0)$ and
$(\sigma_1,-\sigma_0)$ with modelling constant and H\"older norm estimated as follows
\begin{align}
\del M&\lesssim N_0\del N+\del N_0N,\label{wk15}\\
[u_1&-u_0]_\alpha+\|u_1-u_0\|\lesssim N_0\del N+\del N_0\label{wk16} 
\qquad   \mbox{both provided}\;N\leq 1.  
\end{align}
\end{proposition}

\medskip
%%%%%%%%%%%%%%%%%%%%%%%%%%%%%%%%%%%%%%%%%%%%%%%%%%%%%%%%%%%%%%%%%%%%%%%%%%%%%%%%%%%%%%%%%%%%%
%%%%%%%%%%%%%%%%%%%%%%%%%%%%%%%%%%%%%%%%%%%%%%%%%%%%%%%%%%%%%%%%%%%%%%%%%%%%%%%%%%%%%%%%%%%%%%%%%%%%%%%
%%%%%%%%%%%%%%%%%%%%%%%%%%%%%%%%%%%%%%%%%%%%%%%%%%%%%%%%%%%%%%%%%%%%%%%%%%%%%%%%%%%%%%%%%%%%%%%%%%%%%%%

We now proceed to Theorem \ref{Theo}, the main deterministic result of this paper.
It can be seen as a PDE version of the ODE result in \cite[Section 5]{GubinelliControlling}.
Part i) of the theorem provides existence and uniqueness by
a contraction mapping argument, corresponding to 
\cite[Proposition 7]{GubinelliControlling}; 
part ii) provides continuity of the fixed point in the model, the analogue
of the Lyons' sense of continuity for the It\^o map and corresponding to \cite[Proposition 8]{GubinelliControlling}.

\begin{theorem}\label{Theo}
Let $\frac{2}{3}<\alpha<1$ and
let the non-linearities satisfy \eqref{2.15}.

i)  Suppose we are given a distribution $f$ satisfying
\begin{align}\label{x92}
%\sup_{T\le 1}(T^\frac{1}{4})^{2-\alpha}
\|f\|_{\alpha-2}\le N_0
\end{align}
for some constant $N_0\ll 1$;
denote by $v(\cdot,a_0)$ the space-time periodic and mean-free  solution of \eqref{e:v-eq}.
%$(\partial_2-a_0\partial_1^2)v=Pf$.
Suppose further that we are given a one-parameter family of distributions
$\{ v(\cdot,a_0')\diam f\}_{a_0'}$ and a two-parameter family of distributions $\{v(\cdot,a_0')\diam\partial_1^2v(\cdot,a_0)\}_{a_0,a_0'}$
satisfying
\begin{align}\label{x93}
\|
[v,(\cdot)]\diam f \|_{2\alpha-2,2}, \;
\|
[v,(\cdot)]\diam \partial_1^2v \|_{2\alpha-2,2,2}
\le N_0^2.
\end{align}
(In fact, we do not need the highest cross-derivative  $\frac{\partial^2}{\partial a_0'^2}\frac{\partial^2}{\partial a_0^2} [v,(\cdot)_T]\diam \partial_1^2v$).
Then there exists a unique mean-free function $u$ with the properties
\begin{align}
u\;\mbox{is modelled after}\;v\;\mbox{according to}\;a(u)\;\mbox{and}\;\sigma(u),\label{v95}\\
\partial_2u-P(a(u)\diam\partial_1^2u+\sigma(u)\diam f)=0\quad\mbox{distributionally},\label{v94}
\end{align}
under the smallness condition
\begin{align}\label{v96bis}
[u]_\alpha\ll 1.
\end{align}
This unique $u$ satisfies the estimate
\begin{align}\label{v76}
[u]_\alpha+\|u\|\lesssim N_0\quad\mbox{and}\quad M\lesssim N_0^2,
\end{align}
where $M$ denotes the modelling constant in (\ref{v95}).

\medskip

ii) 
Now suppose we have two distributions $f_j$, $j=0,1$, with
\begin{align}\label{v97}
%\sup_{T\le 1}(T^\frac{1}{4})^{2-\alpha}
\|f_{j}\|_{\alpha-2}\le N_0;
\end{align}
and let $v_j(\cdot,a_0)$ be the corresponding solutions of \eqref{e:v-eq}.
Suppose further that for $i=0,1$ we are given four one-parameter families of distributions
$\{ v_i(\cdot,a_0')\diam f_j \}_{a_0'}$ and four two-parameter families of distributions
$\{ v_i(\cdot,a_0')\diam\partial_1^2v_j(\cdot,a_0)\}_{a_0,a_0'}$
satisfying the analogue of (\ref{x93}) including the cross-terms
\begin{align}\label{v14}
\|
[v_i,(\cdot)]\diam f_j \|_{2\alpha-2,2}, \;
\|
[v_i,(\cdot)]\diam \partial_1^2v_j \|_{2\alpha-2,2,2}
\le N_0^2.
\end{align}
We measure the distance of $f_1$ to $f_0$ in terms of a constant $\del N_0$ with
\begin{align}
%\sup_{T\le 1}(T^\frac{1}{4})^{2-\alpha}
\| f_1-f_0\|_{\alpha-2}&\le\del N_0,\label{y56}\\
%
%\sup_{T\le 1}(T^\frac{1}{4})^{2-2\alpha}
\| [v_i,(\cdot)]\diam\{f_1,\partial_1^2v_1\}
-[v_i,(\cdot)]\diam\{f_0,\partial_1^2v_0\}
\|_{2\alpha-2,1,1}
&\le N_0\del N_0,\label{y57}\\
\|
%\{1,\frac{\partial}{\partial a_0'}\}\{1,\frac{\partial}{\partial a_0}\}
  [v_1,(\cdot)]\diam\{f_j,\partial_1^2v_j\}
-[v_0,(\cdot)]\diam\{f_j,\partial_1^2v_j\} \|_{2\alpha-2,1,1}
&\le N_0\del N_0.\label{y58}
\end{align}
If $u_i$, $i=0,1$, denote the corresponding solutions of \eqref{v95}--(\ref{v96bis})
%\&(\ref{v94})
we have
\begin{align}\label{wk54}
[u_1-u_0]_\alpha+\|u_1-u_0\|\lesssim \del N_0.
\end{align}
Moreover, $u_1-u_0$ is modelled after $(v_1,v_0)$ according to $(a(u_1),a(u_0))$ and $(\sigma(u_1),-\sigma(u_0))$
with modelling constant $\del M$ estimated by
\begin{align}\label{wk55}
\del M\lesssim N_0\del N_0.
\end{align}

\end{theorem}

It remains to establish a link between the solution theory presented in Theorem~\ref{Theo} and the classical solution 
theory in the case where $f$ is smooth, e.g.  $f \in C^\beta$ for any $0<\beta<1$.  In this case 
by classical Schauder theory $\sup_{a_0} [\{ \partial_{1}^2, \partial_2 \} v(\cdot, a_0) ]_\beta \lesssim [f]_\beta$ and in particular 
there is the classical choice for the products 
%\begin{align}\label{EM1}
%
$v( \cdot, a_0') \diam \{ f,  \partial_1^2 v(\cdot, a_0) \} = v( \cdot, a_0') \{ f, \partial_1^2 v(\cdot, a_0)\}$.  
%
%\end{align}
In the language of Hairer \cite[Sec. 8.2]{HairerRegularity}, this corresponds to the \emph{canonical model} built from a smooth 
noise term. The only assumption on the products $v( \cdot, a_0') \diam \{ f,  \partial_1^2 v(\cdot, a_0) \} $ entering the definition 
of the singular products
are the regularity bounds  \eqref{x93} expressed in terms of  commutators and they are easily seen to be satisfied in this case. 
For example we have
\begin{align}
\| [v,(\cdot)_T] f \| &= \sup_x \Big| \int \psi_T(x-y) (v(x) - v(y)) f(y) dy \Big| 
\label{EM1A}
\overset{\eqref{1.13}}{ \lesssim} T^{\frac14} \| \{ \partial_1, \partial_2\} v\| \|f\|,
\end{align}
which is much more than needed. However, the canonical definition %\eqref{EM1} 
is by no means the only possible choice of 
product. In fact, as \eqref{x93} is the only requirement on $v( \cdot, a_0') \diam \{ f,  \partial_1^2 v(\cdot, a_0) \}$ we can set 
for example
\begin{align}
%\notag
v( \cdot, a_0') \diam  \{f,  \partial_1^2 v(\cdot, a_0)  \}  :=& v( \cdot, a_0') \{ f ,   \partial_1^2 v(\cdot, a_0)  \} + \{g_1,g_2\}  
\label{EM2}
%v( \cdot, a_0') \diam  :=& v( \cdot, a_0') \partial_1^2 v(\cdot, a_0) +g_2
\end{align}
for a one-parameter family of distributions $g_1$ indexed by $a_0'$ and a two-parameter family $g_2$ 
  indexed by $a_0$, $a_0'$.
For this choice of ``products'' $\diam$ the commutators turn into
\begin{align*}
[v, (\cdot)_T] \diam \{f, \partial_1^2 v\}  = [v, (\cdot)_T]  \{f, \partial_1^2 v\}  - (\{g_1, g_2 \})_T
\end{align*}
so that \eqref{x93} reduces to the regularity assumption 
\begin{align}
 \|g_1 \|_{2\alpha-2,2} , 
 \label{EM3}
  \| g_2 \|_{2\alpha-2,2,2} 
   <\infty.
\end{align}
%This mild assumption leaves a lot of freedom to choose $g_{i}$ 
%%(any 
%%distribution of order $2\alpha -2$ that is smooth in the parameter would do) 
%but  we are mostly  interested in the case where they are constant in $x$  depending only on $a_0$ and $a_0'$. 
The following corollary 
provides a link between  solutions of \eqref{v94}  and classical solutions  in the case where   the products $\diam$ are defined by \eqref{EM2}.
\begin{corollary}\label{Coro}
Let $f$ be a  space-time periodic function in $C^\beta$ for some $0<\beta<1$  and let the products $v( \cdot, a_0') \diam \{ f,  \partial_1^2 v(\cdot, a_0) \} $ be defined by \eqref{EM2} 
for $g_1$, $g_2$
 which satisfy \eqref{EM3}. Then for a periodic mean-free function $u$ the following are equivalent:
\begin{itemize}
\item[i)] $u$ is modelled after $v$ according to $a(u)$ and $\sigma(u)$ and solves 
$\partial_2u-P(a(u)\diam \partial_1^2u+\sigma(u)\diam f)=0$ distributionally.
\item[ii)] $u$ is of class $C^{\beta+2}$ and  a classical solution of
\begin{align*}
\partial_2 u - P\big(a(u) \partial_1^2 u  + a'(u) \sigma^2(u) g_2(\cdot, a(u), a(u))  \\
+ \sigma(u) f 
+ \sigma'(u) \sigma(u) g_1( \cdot, a(u)) \big) =0.
\end{align*}
\end{itemize}
\end{corollary}

%%%%%%%%%%%%%%%%%%%%%%%%%%%%%%%%%%%%%%%%%%%%%%%%%%%%%%%%%%%%%%%%%%%%%%%%%%%%%%%%%%%%%%%%%%%%%%%

%%%%%%%%%%%%%%%%%%%%%%%%%%%%%%%%%%%%%%%%%%
%%%%%%%%%%%%%%%%%%%%%%%%%%%%%%%%%%%%%%%%%%
%%%%%%%%%%%%%%%%%%%%%%%%%%%%%%%%%%%%%%%%%%
%%%%%%%%%%%%%%%%%%%%%%%%%%%%%%%%%%%%%%%%%%
%%%%%%%%%%%%%%%%%%%%%%%%%%%%%%%%%%%%%%%%%%

\section{Stochastic bounds}
\label{s:stochastic}

We now present the stochastic bounds which are necessary as input into our deterministic 
theory. We consider a random distribution $f$, construct (renormalized) commutators, and 
show that the bounds \eqref{x92} and \eqref{x93} hold for these objects. 
The calculations in this section are inspired by a similar reasoning (in a more complicated 
situation)  in \cite[Sec. 5]{HairerPardoux}, \cite[Sec. 10]{HairerRegularity}; 
for the reader's convenience we provide self-contained proofs.

\medskip

Let $f$ be a stationary centered Gaussian distribution which is periodic in both the $x_1$ and the 
$x_2$ direction. Such a distribution is  most conveniently represented in terms of its Fourier series 
\begin{align}\label{SB1}
f(x) = \sum_{k \in (2 \pi \Z)^2} \sqrt{\hat{C}(k)} e^{i k \cdot x}Z_k,
\end{align}
which converges in a suitable topology on distributions. The $Z_k$ are complex-valued centered Gaussians which are independent except for the symmetry constraint
$Z_k = \bar{Z}_{-k}$ and satisfy $ \langle Z_k Z_{-\ell}\rangle = \delta_{k,\ell}$, where as in the introduction
 we use angled brackets $\langle  \cdot \rangle$ to denote the expectation of a random variable.
 The coefficients $\sqrt{\hat{C}}$ are assumed to be real-valued, non-negative, and symmetric
$\sqrt{\hat{C}(k)} = \sqrt{\hat{C}(-k)}$. This notation is chosen because in the case 
where  realizations from $f$ are (say smooth) functions the coefficients in \eqref{SB1} 
coincide with the square root of the  Fourier series of the covariance function.

\medskip

Throughout this section we assume  that $\hat{C}(0)=0$, i.e.  $f$ has vanishing average. 
Our quantitative assumptions on the regularity of $f$ are expressed in terms of $\hat{C}$:
We assume that there exist   $\lambda_1, \lambda_2 \in \R$ and $\alpha' \in (\frac14,1)$  such that 
\begin{align}\label{A1}
 &\hat{C}(k)  \leq \frac{1}{(1+|k_1|)^{\lambda_1} (\sqrt{1+|k_2|})^{\lambda_2}},  \qquad k =(k_1,k_2) \in (2\pi \Z)^2,\\
\notag
&\lambda_1 +  \lambda_2 = -1 + 2 \alpha'    \qquad \lambda_1, \frac{\lambda_2}{2} < 1.   
\end{align}
The second condition may be confusing, because larger values of $\lambda_i$, corresponding to more smoothness for 
$f$, should help our theory. The point here is that decay in one of the directions beyond summability cannot compensate for a lack of decay in the other direction. 
The upcoming Lemma~\ref{lem5} shows that assumption \eqref{A1} corresponds to the regularity assumption \eqref{x92} on $f$.
In order to use the bounds presented in Lemma~\ref{lem5} and Proposition~\ref{P2} as input for the deterministic 
theory in Section~\ref{s:deterministic} we only need the case where $\alpha' >\frac23$ but  the construction 
presented in this section works under the weaker assumption $\alpha'>\frac14$ without additional difficulty.
%
%
%% but this condition does not play a role in the proof of these
%%bounds.
%

\medskip

As in the introduction, we fix an arbitrary Schwartz function  $\varphi$  with $\int_{\R^2} \varphi =1$ and define the rescaling $\varphi_\eps$
and the regularized noise $f_\eps$ as in \eqref{i-def-fepsveps}.
%
%
%in \eqref{i-def-fepsveps} we define the rescaling $\varphi_\eps(x_1,x_2)=\frac{1}{\eps^\frac{3}{4}}\varphi(\frac{x_1}{\eps^\frac{1}{4}},\frac{x_2}{\eps^\frac{1}{2}})$ and 
%\begin{align}\label{def:feps}
%f_\eps = f \ast \varphi_\eps.
%\end{align}
Of course,  $\varphi = \psi_1$ for $\psi_1$ as in
the deterministic analysis is an admissible choice, but in the following 
analysis of stochastic moments the semi-group property  for $\varphi$ is not needed, 
and we therefore do not need to restrict ourselves to this particular choice.
\begin{lemma}\label{lem5}
Let $f$ be  given by \eqref{SB1} satisfying \eqref{A1} for some $\alpha' < 1$
and let $f_\eps$ be as in \eqref{i-def-fepsveps}.
Then we have for any $p<\infty$ and $\alpha< \alpha'$
 \begin{align}\label{SB2}
 \Big\langle  \sup_{\eps \in [0,1]}  \| f_\eps \|_{\alpha-2}^p \Big\rangle^{\frac{1}{p}} \lesssim 1,
 \end{align}
where we use the convention $f_0:=f$. If additionally $0\leq \kappa \leq 1$, then  
 \begin{align}\label{SB2A}
 \big\langle  \big( \sup_{\eps \in (0, 1]}(\eps^{\frac14})^{ -\kappa} \| f_\eps - f \|_{\alpha-2-\kappa}\big) ^p \big\rangle^{\frac{1}{p}} \lesssim 1.
 \end{align} 
 Here and in the proof the implicit constant in $\lesssim$ depends only on  the $\lambda_i$, $p$, $\alpha$  as well as our choice of regularising 
 kernel $\varphi$. 
\end{lemma}
%
%

%\medskip

\medskip 

As before let $v( \cdot, a_0)$ denote the space-time periodic and mean-free
solution to \eqref{e:v-eq}.  We aim at giving a meaning to  the products $v(\cdot, a_0)\diam f$ and $v(\cdot, a_0)\diam \partial_1^2 v(\cdot, a_0')$, 
 and obtaining bounds for the families of 
commutators 
%$[v ( \cdot,a_0), (\cdot)_T ]\diam  f$, $[v (\cdot,a_0), ( \cdot)_T]\diam \partial_1^2 v( \cdot, a_0')$ 
derived from them. 
 The regularities of $v(\cdot, a_0)$, $f$ and $\partial_1^2 v(\cdot, a_0)$ are 
not sufficient to give a deterministic interpretation to these products, and we therefore  seek a probabilistic 
argument to show the convergence of regularized products: We define   $v_\eps(\cdot, a_0)$  as in \eqref{i-def-fepsveps}
% := v(\cdot, a_0) \ast \varphi_\eps$ 
and study the convergence of $v_\eps(\cdot, a_0)f_\eps$, $v_\eps(\cdot, a_0) \partial_1^2 v_\eps(\cdot, a_0')$ as $\eps$ goes to zero by bounding stochastic moments. 
In general under assumption~\eqref{A1} these regularized products do not converge as the regularization is 
removed, but convergence can be enforced by subtracting their expectation. Therefore, we define the renormalized products
%
%
%
%%
%As a first step, in the following lemma we calculate the expectations of $  v_\eps(\cdot, a_0)f_\eps$ and $v_\eps(\cdot, a_0)\partial_1^2 v_\eps(\cdot, a_0')$, which by stationarity do not depend on the point $x \in [0,1)^2$ they are evaluated at.
%%
%%
%%
%
%The regularity assumption \eqref{A1} does not imply that the constants $g_1(\eps, a_0)$ and $g_2(\eps, a_0,a_0')$ converge to a finite limit as $\eps$ tends to zero, although there are interesting 
%cases in which they do converge. This is discussed below, but for the moment we study the convergence of the renormalized products 
%
\begin{align}
\notag
 v_\eps( \cdot, a_0) \diam f_\eps :=& v_\eps( \cdot, a_0) f_\eps - g_1(\eps, a_0) , \\
 \label{e:RenProducts}
 v_\eps( \cdot, a_0) \diam \partial_1^2 v_\eps(\cdot, a_0')   :=& v_\eps( \cdot, a_0) \partial_1^2 v_\eps(\cdot, a_0')  - g_2(\eps, a_0,a_0') , 
\end{align}  
where as in \eqref{Intro77} we  set  $g_1(\eps, a_0) = $   $ \langle v_\eps(\cdot, a_0)f_\eps \rangle$ and $g_2(\eps, a_0, a_0')=
\langle v_\eps(\cdot, a_0)\partial_1^2 v_\eps(\cdot, a_0') \rangle$.
%The corresponding commutators are given as
%%
%%
%%as well as the corresponding commutators,
%%
%\begin{align}
%\notag
%[v_\eps( \cdot, a_0),(\cdot)_T ] \diam f_\eps =&  v_\eps( \cdot, a_0)  (f_\eps )_T - (v_\eps( \cdot, a_0) \diam f_\eps )_T   ,  \\
%\notag
% [v_\eps( \cdot, a_0),(\cdot)_T ]\diam  \partial_1^2 &v_\eps(\cdot, a_0')\\ 
% \label{commutator1}
% =  &   v_\eps( \cdot, a_0)  (\partial_1^2 v_\eps(\cdot, a_0'))_T -  (v_\eps( \cdot, a_0)\diam  \partial_1^2 v_\eps(\cdot, a_0'))_T  .
%\end{align}
%%  
%Observe that while the singular products appearing in this expression, $v_\eps( \cdot, a_0) f_\eps$ and $ v_\eps( \cdot, a_0)$ $ \partial_1^2 v_\eps(\cdot, a_0')$, are renormalized by subtracting the expectation, the products $ v_\eps( \cdot, a_0)  (f_\eps )_T$ and $v_\eps( \cdot, a_0)  (\partial_1^2 v_\eps(\cdot, a_0'))_T$ are not changed. In particular, unlike the renormalized products in 
%\eqref{e:RenProducts} the renormalized commutators in \eqref{commutator1} do not have vanishing expectation. 

\medskip

The key result of this section is the following proposition which shows the convergence of 
the renormalized products and provides a control for  stochastic moments of the renormalized commutators   as well as their derivatives with respect to $a_0, a_0'$. 
\begin{proposition}\label{P2}
Let $f$ be a stationary centered Gaussian distribution given by \eqref{SB1} satisfying \eqref{A1} for some $\frac14<\alpha'<1$,
let $v(\cdot, a_0')$ be the space-time periodic mean-free solution of \eqref{e:v-eq} 
and let $f_\eps$ and $v_\eps(\cdot, a_0')$ be as in \eqref{i-def-fepsveps}.
%
%and $\hat{v}_\eps( \cdot, a_0') = \hat{G}( \cdot, a_0') \hat{f_\eps}$ for $G( \cdot, a_0')$
%as in \eqref{e:Greens}.

\medskip
i) For any $n, m \geq 0$ the random distributions $\big(\frac{\partial}{\partial a_0}\big)^n $ $\big(\frac{\partial}{\partial {a_0'}}\big)^m $  $ v_\eps(\cdot, a_0) \diam \{ f_\eps, \partial_1^2 v_\eps( \cdot, a_0') \}$  converge as $\eps \to 0$.  
This convergence takes place almost surely uniformly over $a_0, a_0'$ and 
with respect to any $C^{\alpha-2}$ norm for $\alpha < \alpha'$. 
We denote the limits by $\big(\frac{\partial}{\partial a_0}\big)^n $ $\big(\frac{\partial}{\partial {a_0'}}\big)^m $   $v(\cdot, a_0)$  $\diam$  $\{ f, \partial_1^2 v( \cdot, a_0') \}$.
\medskip

ii) For all $p<\infty$ we have the estimates 
\begin{align}
\notag
  \Big\langle
% 	 \Big(
%	 	 \sup_{a_0, a_0' \in [\lambda,\frac{1}{\lambda}]} 
		 \sup_{\eps_0, \eps_1 \in [0,1]}  
%		 \sup_{T \leq 1}  (T^{\frac14})^{2-2\alpha'} 
		 	\| 
%				\frac{\partial^n}{\partial a_0^n} 
%				\frac{\partial^{m}}{\partial (a_0')^{m}} 
				[ v_{\eps_0},(\cdot) ]
  				\diam f_{\eps_1}
			\|_{2\alpha-2,n} ^p
%	\Big)^p 
\Big\rangle^{\frac{1}{p}} 
&\lesssim 1,\\
  \Big\langle
% 	 \Big(
%	 	 \sup_{a_0, a_0' \in [\lambda,\frac{1}{\lambda}]} 
		 \sup_{\eps_0,\eps_1 \in [0,1]}  
%		 \sup_{T \leq 1}  (T^{\frac14})^{2-2\alpha'} 
		 	\| 
%				\frac{\partial^n}{\partial a_0^n} 
%				\frac{\partial^{m}}{\partial (a_0')^{m}} 
				[ v_{\eps_0},(\cdot) ]
 \label{FinalStoch1}
  				\diam \partial_1^2 v_{\eps_1}  
			\|_{2\alpha-2,n,m} ^p
%	\Big)^p 
\Big\rangle^{\frac{1}{p}} 
&\lesssim 1,
\end{align}
%
%,
as well as for $0< \kappa \ll 1$ (where $\ll$ depends only on $\lambda_1,\lambda_2$)
\begin{align}
\notag
\Big\langle
	\Big( 
%		\sup_{a_0, a_0' \in [\lambda,\frac{1}{\lambda}]}   
		 \sup_{\eps \in (0,1]}  (\eps^\frac14)^{-\kappa} 
%		\sup_{T \leq 1} (T^{\frac14})^{2-2\alpha'+\kappa}  \\
%\notag
%&\qquad		
		\big\|   
%			\frac{\partial^n}{\partial a_0^n}  
%			\frac{\partial^{m}}{\partial (a_0')^{m}}
%			\Big(
				[v_\eps,(\cdot) ] \diam					
				 f_\eps 
				-[v,(\cdot) ] 
				\diam  f  
%				\Big)  
			\big\|_{2\alpha-2-\kappa,n}
	\Big)^p 
\Big\rangle^{\frac{1}{p}} 
&\lesssim  1,\\
\Big\langle
	\Big( 
%		\sup_{a_0, a_0' \in [\lambda,\frac{1}{\lambda}]}   
		 \sup_{\eps \in (0,1]}  (\eps^\frac14)^{-\kappa} 
%		\sup_{T \leq 1} (T^{\frac14})^{2-2\alpha'+\kappa}  \\
%\notag
%&\qquad		
		\big\|   
%			\frac{\partial^n}{\partial a_0^n}  
%			\frac{\partial^{m}}{\partial (a_0')^{m}}
				[v_\eps,(\cdot) ] 					
				\diam  \partial_1^2 v_\eps 
 \label{FinalStoch3}
%&\qquad \qquad  
				-[v,(\cdot) ] 
				\diam  \partial_1^2 v  
			\big\|_{2\alpha-2-\kappa,n,m}
	\Big)^p 
\Big\rangle^{\frac{1}{p}} 
&\lesssim  1,
\end{align}
where here and in the proof $\lesssim $ means up to a constant depending only on $n$, $m$, the $\lambda_i$, $\alpha$, $\kappa$, $p$, the ellipticity contrast $\lambda$ as well as the specific choice of regularising kernel $\varphi$.
In both estimates the subscripts $n,m$ in the norms refer to parameter derivatives with respect to $a_0, a_0'$
as in \eqref{conv1} and \eqref{conv2}.
\end{proposition}
Proposition~\ref{P2} is a consequence of the following estimate on the second moments of  commutators. 
\begin{lemma}\label{cch1}
Let $f$ and $v(\cdot, a_0)$ be as in  Proposition~\ref{P2}.
%be defined by \eqref{SB1} for $\hat{C}$ satisfying \eqref{A1} for some $0<\alpha'<1$ and
%let $v(\cdot, a_0')$ be the space-time periodic mean-free solution of \eqref{e:v-eq}.
% and let $v(\cdot,a_0)$ solve 
% and let $G$ be the Greens function defined in \eqref{e:Greens} for some $a_0 \in [\lambda,\frac{1}{\lambda}]$. 
%
Let $\hat{M}_1, \hat{M}_2$ be Fourier multipliers satisfying  
\begin{align}\label{mn00}
\hat{M}_i(k) = \overline{\hat{M}_i(-k)} \quad \text{and} \quad  |\hat{M}_i(k) | \leq (k_1^4 + k_2^2)^{\frac{\kappa_i}{4}},  \quad k \in (2 \pi\Z)^2, i = 1,2
\end{align}
for $0 \leq \kappa_1, \kappa_2 \ll 1$ (where $\ll$ depends only on $\lambda_1, \lambda_2$).
 Let $f'$ and $v'(\cdot, a_0)$ be defined through their Fourier series
\[
 \hat{f'} = \hat{M}_1 \hat{f} \qquad  \text{and} \qquad \hat{v'}(\cdot, a_0) = \hat{M}_2 \hat{v}(\cdot, a_0).
 \]
 We  make the qualitative assumption that $f'$ and $v'(\cdot, a_0)$ are smooth and set 
% 4
$$v'(\cdot, a_0) \diam f': = v'(\cdot, a_0) f' - \langle v'(\cdot, a_0) f' \rangle.$$ 
%and 
% $$[v'(\cdot, a_0), (\cdot)_T] \diam f' := v' (\cdot, a_0)f'_T - (v'(\cdot, a_0) \diam f')_T+  \langle v'(\cdot, a_0) f' \rangle.$$ 
 Then 
for all $a_0 \in [\lambda, \frac{1}{\lambda}]$
\begin{align}\label{mn01}
\langle ( [v'(\cdot, a_0), ( \cdot)_T] \diam f' )^2 \rangle^{\frac12} \lesssim (T^{\frac14})^{2 \alpha' -2 - \kappa_1 - \kappa_2}.
\end{align} 
Here and in the proof the implicit constant depends on $\lambda_1, \lambda_2$, $\kappa_1, \kappa_2$ as well as the ellipticity contrast $\lambda$ (but not on the qualitative smoothness assumption on $f', v'$).
\end{lemma}
In the proof of Proposition~\ref{P2} this lemma is used in the form of the following immediate corollary:
\begin{corollary}\label{cch2}
%
%
%For $\eps>0$ and $a_0, a_0' \in [\lambda,\frac{1}{\lambda}]$  let $f_\eps = \ph_\eps \ast f$  and $\hat{v}_\eps( \cdot, a_0') = \hat{G}( \cdot, a_0') \hat{f_\eps}$ for $G( \cdot, a_0')$ as in \eqref{e:Greens}. Furthermore, let $[v_\eps(\cdot, a_0), (\cdot)_T]$  $\diam \{f_\eps, \partial_1^2 v_\eps(\cdot, a_0') \}$ be defined as in \eqref{commutator1}. 
%
Let $f$, $f_\eps$, $v$ and $v_\eps$ be as in Proposition~\ref{P2}.
Then for $n,m \geq 0$ we have
\begin{align}
%\notag
&\Big\langle
	 \Big(
	 	 \Big[ 
		 	\Big(\frac{\partial}{\partial a_0}\Big)^n v_{\eps_0}( \cdot, a_0),(\cdot)_T 
		\Big] 
		\diam \Big\{
			  f_{\eps_1}, \Big( \frac{\partial}{\partial a_0'}\Big)^{m} \partial_1^2 v_{\eps_1}( \cdot, a_0')
			 \Big \} 
	\Big)^2 
\Big\rangle^{\frac12}
%\\
\label{mn02}
%& \qquad \qquad \qquad \qquad \qquad \qquad \qquad \qquad \qquad 
 \lesssim (T^{\frac14})^{2 \alpha' - 2}.
\end{align}  
Furthermore, we have for $0 \leq \kappa \ll 1$ 
($\ll$ depends only on $\lambda_1, \lambda_2$) and for $i=0,1$
\begin{align}
%\notag
&\Big\langle 
	\Big(
	 	\eps_i \frac{\partial}{\partial \eps_i}
	 	\Big(
			 \Big[ 
			\Big( \frac{\partial}{\partial a_0}\Big)^n v_{\eps_0}( \cdot, a_0),(\cdot)_T 
			 \Big] 
			 \diam \Big\{  f_{\eps_1}, \Big( \frac{\partial}{\partial {a_0'}}\Big)^{m} \partial_1^2 v_{\eps_1} \Big\}
		\Big) 
	\Big)^2 
\Big\rangle^{\frac12} 
%\\
\label{mn03}
%& \qquad \qquad \qquad \qquad \qquad \qquad \qquad \qquad \qquad
\lesssim  
(T^{\frac14})^{2 \alpha' - 2 - \kappa}(\eps_i^{\frac14})^{ \kappa}.
\end{align} 
Here and in the proof the implicit constant 
depends on $\lambda_1, \lambda_2$, $\kappa$  the ellipticity contrast $\lambda$, $n,m$
as well as the specific choice of regularising kernel $\varphi$.
\end{corollary}

\medskip

Finally, the following lemma deals with the behaviour 
of the expectations $g_1,g_2$ as the regularization is removed.

%  come back to the products and commutators without renormalization.
%
%
% According to Lemma~\ref{lem6} the constants 
%$g_1(\eps, a_0)$ converge to a  limit if and only if 
%%
%%
%
%

%
\begin{lemma}
\label{lem6}
i)
For $\eps>0$ we have
\begin{align}
\label{const1}
g_1(\eps, a_0) 
= &
%\big\langle 
%	v_\eps(\cdot, a_0)f_\eps 
%\big\rangle 
%= 
\sum_{k \in (2 \pi \Z)^2 \setminus \{0\}}  
	\frac{a_0k_1^2}{a_0^2 k_1^4 + k_2^2} 
	\hat{C}(k) 
	|\hat{\ph}_{\eps}(k)|^2 ,   \\
%\notag
g_2(\eps, a_0, a_0') 
=& 
%\big\langle 
%	v_\eps(\cdot, a_0)\partial_1^2 v_\eps(\cdot, a_0') 
%\big\rangle \\
\label{const2}
%=& 
\sum_{k \in (2 \pi \Z)^2 \setminus \{0\}} 
	\frac{(-a_0 a_0'k_1^4 + k_2^2)k_1^2 }{(a_0^2 k_1^4 + k_2^2)((a_0')^2 k_1^4 + k_2^2)} 		\hat{C}(k) 
	|\hat{\ph}_{\eps}(k)|^2  . 
\end{align}

\medskip
ii) The expectation $g_1(\eps, a_0)$ converges to a finite limit as $\eps \to 0$  if and only if 
\begin{align}\label{A2}
\sum_{k \in (2 \pi \Z)^2 \setminus \{0\}}  \frac{k_1^2}{ k_1^4 + k_2^2} \hat{C}(k)  < \infty.
\end{align}

 If \eqref{A2} holds, then $g_2(\eps, a_0,a_0')$ as well as all parameter derivatives $(\frac{\partial}{ \partial a_0})^n g_1(\eps, a_0)$ 
and $(\frac{\partial}{ \partial a_0})^n $   $(\frac{\partial}{ \partial {a_0'}})^m $ $ g_2(\eps, a_0, a_0')$  for $n,m \geq 0$ converge as well.
\end{lemma}
 In particular we immediately get the following:

\begin{corollary}\label{Cc5}
Assume that both \eqref{A1}  and \eqref{A2} hold. Then the statements of Proposition~\ref{P2} remain true if all of the
renormalized products are replaced by products without renormalization.
\end{corollary}
%
%
%The limits which exist under the assumptions of this corollary will be denoted by  $[v( \cdot, a_0),(\cdot)_T ] f $, $ [v( \cdot, a_0),(\cdot)_T ] \partial_1^2 v(\cdot, a_0)$ etc..

%

\medskip

We finish this section by 
 discussing the assumptions \eqref{A1} and  \eqref{A2} in  particular cases. First consider the case
\begin{align}
\label{S1}
\hat{C}(k) 
= &
\frac{1}{(1+|k_1|)^{\lambda_1} (\sqrt{1+|k_2|})^{\lambda_2}}.
\end{align}
For this choice of $\hat{C}$ the regularity assumption \eqref{A1} is equivalent to 
\begin{align}\label{A11}
\lambda_1 + \lambda_2 \geq -1 +2 \alpha',  \qquad \lambda_1 > -3 +2\alpha', \quad \text{and} \quad  \lambda_2 > -2 +2\alpha'.
\end{align}
Note that equality is not necessary in the first condition, because in the case of strict inequality, one can find $\lambda_1' \leq \lambda_1$ and $\lambda_2' \leq \lambda_2$ that satisfy \eqref{A1} with equality. 
%However, $\lambda_1 \leq -3 +2\alpha$ or 
%$\lambda_2 \leq -2 +2\alpha$ can never be compensated without violating the second condition in \eqref{A1}.
The condition \eqref{A2} on the other hand is equivalent to 
\begin{equation}\label{A22}
\lambda_1 +  \lambda_2 > 1  
\qquad 
\lambda_1 >-1, 
\quad \text{and}  \quad  
\lambda_2 >- 2.
\end{equation}
%
%For any $\alpha' \in (0, 1)$ the first requirement in \eqref{A11} is weaker than the corresponding assumption in \eqref{A22}. 
An interesting case in which both assumptions 
are satisfied and for which our theory can therefore be applied without renormalization is the case where $\lambda_1 >1$ and 
$\lambda_2 =0$; this corresponds to the case of noise which is white in the time-like variable $x_2$ but ``trace-class'' in $x_1$. 
However, if we are willing to accept renormalization, the regularity requirement in the $x_1$ direction reduces to $\lambda_1 >\frac13$ (recall that the deterministic analysis is applicable if $\alpha >\frac23$). Another interesting case is the 
covariance
\begin{align*}
\hat{C}(k) = \delta_{k_2,0 }\frac{1}{ (1+|k_1|)^{{\lambda_1}}},
\end{align*}
which corresponds to the choice $\lambda_2 = \infty$ in \eqref{S1} and
yields a noise term which only depends on the space-like $x_1$ variable. The parabolic equations with constant diffusion 
coefficients driven by such a noise term has  recently been studied as \emph{parabolic Anderson model} in two and three spatial dimensions
\cite{GubinelliImkellerPerkowski,HairerLabbe1,HairerLabbe2,Bailleul}. Our theory applies without renormalization for all $\lambda_1 >-1$, which 
covers in particular the case of  one-dimensional spatial white noise, $\lambda_1 =0$. If we admit 
renormalization we can go all the way to $\lambda_1 >- \frac{5}{3} $.
% by choosing $\lambda_2 < 2$ and $\alpha>\frac23$ 
%as close to $2$ and $\frac{2}{3}$, respectively, as we please. 
This covers  the case $\lambda_1=-1$ for which the noise $f$ has the same scaling 
behaviour as spatial white noise in two dimensions (both are distributions of regularity $C^{-1-}$) but it does not cover the case $\lambda_1 = -2$ for  which the noise scales like spatial white 
noise in three dimensions.

%%%%%%%%%%%%%%%%%%%%%%%%%%%%%%%%%%%%%%%%%%%%%%%%%%%%%%%%%%%%%%%%%%%%%%%%%%%%%%%%%%%%%%%%%%%%%%%%%%%%%%%
%%%%%%%%%%%%%%%%%%%%%%%%%%%%%%%%%%%%%%%%%%%%%%%%%%%%%%%%%%%%%%%%%%%%%%%%%%%%%%%%%%%%%%%%%%%%%%%%%%%%%%%
%%%%%%%%%%%%%%%%%%%%%%%%%%%%%%%%%%%%%%%%%%%%%%%%%%%%%%%%%%%%%%%%%%%%%%%%%%%%%%%%%%%%%%%%%%%%%%%%%%%%%%%
%%%%%%%%%%%%%%%%%%%%%%%%%%%%%%%%%%%%%%%%%%%%%%%%%%%%%%%%%%%%%%%%%%%%%%%%%%%%%%%%%%%%%%%%%%%%%%%%%%%%%%%

\section{Proofs for the deterministic analysis}
\label{sec:DetProofs}
\subsection{Proof of Theorem \ref{Theo}}
%\addcontentsline{toc}{section}{\protect\numberline{}Proof of Theorem \ref{Theo}}%
%\addcontentsline{toc}{section}{\protect\numberline{}Theorem \ref{Theo}}%

\newcounter{TS} % proofstep = 0
\refstepcounter{TS} % increases value by 1

We write for abbreviation $[\cdot]=[\cdot]_\alpha$.
We consider the map defined through  
\begin{align}\label{wi34}
(\bar u,\bar a,\bar\sigma)\mapsto(\sigma:=\sigma(\bar u),a:=a(\bar u),\sigma\diam f,a\diam\partial_1^2v)
\mapsto (u,a,\sigma),
\end{align}
where $u$ is the solution provided by Proposition \ref{P}. This is the map
of which we seek to characterize the fixed point. Note that the right hand side depends on $\bar a$  and $\bar\sigma$
via the definition of the products $\sigma\diam f$ and $a\diam\partial_1^2$.

\medskip
{\sc Step} \arabic{TS}.\label{nonlinear}\refstepcounter{TS} Pointwise nonlinear transformation, application of Lemma \ref{L0}.
We work under the assumptions of part ii) of the theorem on the distributions $f_j$ and the off-line
products $v_i\diam f_j$, $v_i\diam\partial_1^2v_j$.
Suppose we are given two triplets $(\bar u_i,\bar a_i,\bar\sigma_i)$, $i=0,1$, of functions satisfying the constraints 
\begin{align}\label{wi28}
\bar\sigma_i\in[-1,1],\;\bar a_i\in[\lambda,\frac{1}{\lambda}],\quad[\bar\sigma_i],[\bar a_i]\le 1.
\end{align}
We measure the size of $\{(\bar u_i,\bar a_i,\bar\sigma_i)\}_{i}$ and their distance through
\begin{align}
\bar M&:=\max_{i}(M_{\bar u_i}+[\bar u_i])+N_0,\label{wi26}\\
\del\bar M&:=M_{\bar u_1-\bar u_0}+[\bar u_1-\bar u_0]+\|\bar u_1-\bar u_0\|\nonumber\\
&+N_0([\bar\sigma_1-\bar\sigma_0]+\|\bar\sigma_1-\bar\sigma_0\|
+[\bar a_1-\bar a_0]+\|\bar a_1-\bar a_0\|)+\del N_0,\label{wi27}
\end{align}
where $M_{\bar u_i}$ denotes the constant in the modelledness of $\bar u_i$ after $v_i$ according to $\bar a_i$ and $\bar\sigma_i$,
and where $M_{\bar u_1-\bar u_0}$ denotes the constant in the modelledness of
$\bar u_1-\bar u_0$ after $(v_1,v_0)$ according to $(\bar a_1,\bar a_0)$ and $(\bar\sigma_1,-\bar\sigma_0)$.

\medskip

We now consider $\sigma_i:=\sigma(\bar u_i)$ and $a_i:=a(\bar u_i)$. We claim
\begin{align}
\sigma_i&\in[-1,1],\;a_i\in[\lambda,\frac{1}{\lambda}],\;[\sigma_i] ,[a_i]\ll 1
\quad\mbox{provided}\;\max_i[\bar u_i]\ll 1,\label{wj24}\\
\tilde M&\lesssim \bar M\quad\mbox{provided}\;\max_i[\bar u_i]\le 1,\label{wi35}\\
\del\tilde M&\lesssim \del\bar M\quad\mbox{provided}\;\bar M\le 1,\label{wj16}
\end{align}
where we define in analogy with (\ref{wi26}) and (\ref{wi27}): 
\begin{align}
\tilde M&:=\max_{i}(M_{\sigma_i}+[\sigma_i]+M_{a_i}+[a_i])+N_0,\label{wj25}\\
\del\tilde M&:=M_{\sigma_1-\sigma_0}+[\sigma_1-\sigma_0]+\|\sigma_1-\sigma_0\|\nonumber\\
&+N_0\big([\omega_1-\omega_0]+\|\omega_1-\omega_0\|+[\bar a_1-\bar a_0]+\|\bar a_1-\bar a_0\|\big)\nonumber\\
&+M_{a_1-a_0}+[a_1-a_0]+\|a_1-a_0\|\nonumber\\
&+N_0\big([\mu_1-\mu_0]+\|\mu_1-\mu_0\|+[\bar a_1-\bar a_0]+\|\bar a_1-\bar a_0\|\big)+\del N_0,\label{wj21}
\end{align}
with the understanding that $\sigma_i$ is modelled after
$v_i$ according to $\bar a_i$ and $\omega_i:=\sigma'(\bar u_i)\bar\sigma_i$ and constant $M_{\sigma_i}$,
that $a_i$ is modelled after
$v_i$ according to $\bar a_i$ and $\mu_i:=a'(\bar u_i)\bar\sigma_i$ and constant $M_{a_i}$,
that $\sigma_1-\sigma_0$ is modelled after $(v_1,v_0)$ according to $(\bar a_1,\bar a_0)$ and $(\omega_1,-\omega_0)$ 
and a constant we name $M_{\sigma_1-\sigma_0}$, and 
that $a_1-a_0$ is modelled after $(v_1,v_0)$ according to $(\bar a_1,\bar a_0)$ and $(\mu_1,-\mu_0)$        
and a constant we name $M_{a_1-a_0}$. 

\medskip

It is obvious from (\ref{2.15})
that we have (\ref{wj24}) under the assumption $\max_i[\bar u_i]\ll 1$.
Estimate (\ref{wi35}) follows from part i) of Lemma \ref{L0} with $u$ replaced by $\bar u_i$ and
the generic nonlinearity $b$ replaced by $\sigma$ and by $a$, respectively, 
(using our assumptions (\ref{2.15})). More precisely, (\ref{wi35})
follows from (\ref{wj14}) by $[\bar u_i]\le 1$.
We now turn to (\ref{wj16}), which by definitions (\ref{wi27}) of $\del\bar M$ and (\ref{wj21}) of $\del\tilde M$
and because of $N_0\le 1$ we may split into the four statements
%x
\begin{align*}
M_{\sigma_1-\sigma_0}+[\sigma_1-\sigma_0]+\|\sigma_1-\sigma_0\|
&\lesssim M_{\bar u_1-\bar u_0}+[\bar u_1-\bar u_0]+\|\bar u_1-\bar u_0\|,\nonumber\\
[\omega_1-\omega_0]+\|\omega_1-\omega_0\|&\lesssim [\bar\sigma_1-\bar\sigma_0]+\|\bar\sigma_1-\bar\sigma_0\|\nonumber\\
&+[\bar u_1-\bar u_0]+\|\bar u_1-\bar u_0\|,\nonumber\\
M_{a_1-a_0}+[a_1-a_0]+\|a_1-a_0\|
&\lesssim M_{\bar u_1-\bar u_0}+[\bar u_1-\bar u_0]+\|\bar u_1-\bar u_0\|,\nonumber\\
[\mu_1-\mu_0]+\|\mu_1-\mu_0\|&\lesssim [\bar\sigma_1-\bar\sigma_0]+\|\bar\sigma_1-\bar\sigma_0\|\nonumber\\
&+[\bar u_1-\bar u_0]+\|\bar u_1-\bar u_0\|,\nonumber\\
\mbox{all provided}\;\max_i(M_{\bar u_i}+[\bar u_i])\le 1,
\end{align*}
where we also used the definition (\ref{wi26}) of $\bar M$. This is a consequence of part ii) of Lemma \ref{L0}
with $(\bar u_i,\bar\sigma_i, \bar{a}_i)$ playing the role of $(u_i,\sigma_i,a_i)$. 
The first two estimates follow from replacing the generic nonlinearity $b$ by $\sigma$,
the last two estimates from replacing it by $a$. The first and the third estimate are a consequence of (\ref{wj11}),
the second and fourth one of (\ref{wj12}), in which we use (\ref{wj24}). It is on all four 
we use our full assumptions (\ref{2.15}) on the nonlinearities $\sigma$ and $a$.

\medskip

%%%%%%%%%%%%%%%%%%%%%%%%%%%%%%%%%%%%%%%%%%%%%%%%%%%%%%%%%%%%%%%%%%%%%%%%%%%%%%%%%%%%%%%%%%%%%%%%%%%%%%%%%%%%%

{\sc Step} \arabic{TS}.\label{productsI}\refstepcounter{TS} 
Using the off-line products, application of Corollary \ref{C2}.
We claim that under the hypothesis of part ii) of the theorem
on the distributions $f_j$ and the off-line products $v_i\diam f_j$ \&
$v_i\diam\partial_1^2v_j$ we have the commutator estimates
\begin{align}
%\sup_{T\le 1}&(T^\frac{1}{4})^{2-2\alpha}
\|[\sigma_i,(\cdot)]\diam f_j\|_{2\alpha-2}
&\lesssim N_0\tilde M,\label{wi36}\\
%
%\sup_{T\le 1}&(T^\frac{1}{4})^{2-2\alpha}\|
 \| [\sigma_i,(\cdot)]\diam f_1-[\sigma_i,(\cdot)]\diam f_0\|_{2\alpha-2} &\lesssim \del N_0\tilde M,\label{wj44}\\
%\sup_{T\le 1}&(T^\frac{1}{4})^{2-2\alpha}
\|[\sigma_1,(\cdot)]\diam f_j-[\sigma_0,(\cdot)]\diam f_j\|_{2\alpha-2} &\lesssim N_0 \del\tilde M,\label{wj45}\\
%
%\sup_{T\le 1}&(T^\frac{1}{4})^{2-2\alpha}
\|
[a_i,(\cdot)]\diam\partial_1^2v_j\|_{2\alpha-2,2} &\lesssim N_0\tilde M,\label{wi37}\\
%\end{align}
%and
%\begin{align}
%\sup_{T\le 1}&(T^\frac{1}{4})^{2-2\alpha}
\|
[a_i,(\cdot)]\diam\partial_1^2v_1
-[a_i,(\cdot)]\diam\partial_1^2v_0  \|_{2\alpha-2,1}
&\lesssim \del N_0\tilde M,\label{wj46}\\
%\sup_{T\le 1}&(T^\frac{1}{4})^{2-2\alpha}
\|
[a_1,(\cdot)]\diam\partial_1^2v_j 
-[a_0,(\cdot)]\diam\partial_1^2v_j\|_{2\alpha-2,1}
&\lesssim N_0 \del\tilde M.\label{wj47}
\end{align}
This is an application of Corollary \ref{C2} with $(N_1,\del N_1)=(N_0,\del N_0)$.
Estimate (\ref{wi36}) is an application of Corollary \ref{C2} i) with $u$ replaced by $\sigma_i$; the hypotheses (\ref{wj36}) and (\ref{wj41})
are contained in the theorem's assumptions (\ref{x92}) and (\ref{x93}) (note that $f$ does not depend on an extra parameter $a_0'$).
%Corollary \ref{C2} i) is applied with $u$ replaced by $\sigma_i$; 
The output (\ref{wj40}) turns into (\ref{wi36})
since by definition (\ref{wj25}), $M_{\sigma_i}+N_0\le\tilde M$.
Estimate (\ref{wj44}) is an application of Corollary \ref{C2} ii)  still applied with $u$ replaced by $\sigma_i$; the hypotheses (\ref{wj52}) and (\ref{wj48})
are contained in the theorem's assumptions (\ref{y56}) and (\ref{y57}).
The output (\ref{wk35}) turns into (\ref{wj44})
as in the previous application.
Estimate (\ref{wj45}) is an application of Corollary \ref{C2} iii) now applied with $u_i$ replaced by $\sigma_i$ (and thus $(\sigma_i,a_i)$ replaced by $(\omega_i,\bar a_i)$); the hypotheses (\ref{wj50}) and (\ref{wj49})
are contained in the theorem's assumptions (\ref{v14}) and (\ref{y58}).
The output (\ref{wk36}) turns into (\ref{wj45}),
since by definition \eqref{wj21} we have
\begin{align*}
M_{\sigma_1-\sigma_0}
+N_0([\omega_1-\omega_0]+\|\omega_1-\omega_0\|+[\bar a_1-\bar a_0]+\|\bar a_1-\bar a_0\|)+\del N_0
\le \del \tilde M.
\end{align*}

\medskip

The arguments for (\ref{wi37}), (\ref{wj46}), and (\ref{wj47}) follow the same lines of those
for  (\ref{wi36}), (\ref{wj44}), and (\ref{wj45}), respectively. The only difference is that in all instances,
the distribution $f_j$ is replaced by the family of distributions $\partial_1^2v_j(\cdot,a_0)$
(and $a_i$ plays the role of $u$ in Corollary \ref{C2}). Hence the  hypotheses 
(\ref{wj36}) and (\ref{wj52}) in Corollary \ref{C2} turn into
\begin{align*}
\|
\partial_1^2v_{j}\|_{\alpha-2,2} \lesssim N_0,\qquad \qquad
\|
\partial_1^2(v_{1}-v_{0})\|_{\alpha-2,1}\lesssim \del N_0.
\end{align*}
This follows from Step \ref{C1S0} in the proof of Corollary \ref{C1} via (\ref{1.13}).

\medskip

%%%%%%%%%%%%%%%%%%%%%%%%%%%%%%%%%%%%%%%%%%%%%%%%%%%%%%%%%%%%%%%%%%%%%%%%%%%%%%%%%%%%%%%%%%%%%%%%%%%%%%%%%%%%%

{\sc Step} \arabic{TS}\label{TS3}\refstepcounter{TS}. Application of Proposition \ref{P}.
We claim that under the hypothesis of part ii) of the theorem
regarding the distributions $f_j$ and the off-line products $v_i\diam f_j$ and
$v_i\diam\partial_1^2v_j$
\begin{align}
M&\lesssim N_0(\tilde M+1)&\quad\mbox{provided}\;\max_i[\bar{u}_i]\ll 1,\label{wk37}\\
\max_iM_{u_i}&\lesssim N_0\tilde M&\quad\mbox{provided}\;\max_i[\bar{u}_i]\ll 1,\label{wk46}\\
\del M&\lesssim N_0\del\tilde M+\del N_0&\quad\mbox{provided in addition}\;\tilde M\lesssim 1,\label{wk38}\\
M_{u_1-u_0}&\lesssim N_0\del\tilde M+\del N_0\tilde M
&\quad\mbox{provided in addition}\;\tilde M\lesssim 1,\label{wk47}
\end{align}
where we define consistently with (\ref{wi26}) and (\ref{wi27})
\begin{align}
M:&=\max_{i}(M_{u_i}+[u_i])+N_0,\label{wk40}\\
\del M:&=M_{u_1-u_0}+[u_1-u_0]+\|u_1-u_0\| \nonumber\\
&+N_0\Big([\sigma_1-\sigma_0]
+\|\sigma_1-\sigma_0\|+[a_1-a_0]+\|a_1-a_0\|\Big)+\del N_0.\label{wk39}
\end{align}
Indeed, (\ref{wk37}) and (\ref{wk46}) are an application of part i) of Proposition \ref{P}: The hypothesis
(\ref{wk59}) of the proposition is built into the definition (\ref{wj25}) of $\tilde M$,
so that $\tilde M$ here plays the role of $N$ in the proposition. The hypothesis 
(\ref{wj58}) is identical to the theorem's assumption (\ref{v97}), hypothesis \eqref{wj59}   was 
established in (\ref{wj24}), 
hypotheses (\ref{wj56}) and (\ref{wj57}) are contained in (\ref{wi36}) 
and (\ref{wi37}) of Step \ref{productsI} which is consistent with $\tilde M$ playing the role of $N$ there. The combination
of (\ref{wk17}) and (\ref{wk18}) amounts to (\ref{wk37}) by definition (\ref{wk40}) of $M$.
Estimate (\ref{wk17}) by itself amounts to (\ref{wk46}).

\medskip

Estimate (\ref{wk38}) in turn is a consequence of part ii) of Proposition \ref{P}: Hypothesis
(\ref{wk50}) of the proposition is build into the definition (\ref{wj21}) of $\del\tilde M$,
so that $\del\tilde M$ here plays the role of $\del N$ in the proposition.
Hypotheses (\ref{wj97bis}) and (\ref{wj98bis}) are identical to (\ref{wi36})
and (\ref{wi37}) of Step \ref{productsI}.
Hypothesis (\ref{wk24}) is identical to our assumption (\ref{y56}), hypotheses
(\ref{wk13}), (\ref{wk25}), (\ref{wk14}), and (\ref{w26}) are identical to
(\ref{wj44}), (\ref{wj45}), (\ref{wj46}), and (\ref{wj47}) in Step \ref{productsI}.
The outcome (\ref{wk15}) of the proposition turns into (\ref{wk47}). The latter trivially 
for $\tilde M\lesssim 1$ implies
\begin{align*}
M_{u_1-u_0}\lesssim N_0\del\tilde M+\del N_0,
\end{align*}
whereas the outcome (\ref{wk16}) of the proposition assumes the form
\begin{align*}
[u_1-u_0]+\|u_1-u_0\|\lesssim N_0\del\tilde M+\del N_0.
\end{align*}
By definition (\ref{wj21}) of $\del\tilde M$  we have
\begin{align*}
[\sigma_1-\sigma_0]+\|\sigma_1-\sigma_0\|
+[a_1-a_0]+\|a_1-a_0\|\le \del\tilde M.
\end{align*}
The combination of the last three statement yields (\ref{wk38}) in view of definition (\ref{wk39}).

\medskip

%%%%%%%%%%%%%%%%%%%%%%%%%%%%%%%%%%%%%%%%%%%%%%%%%%%%%%%%%%%%%%%%%%%%%%%%%%%%%%%%%%%%%%%%%%%%%%%%%%%%%

{\sc Step} \arabic{TS}.\label{TS4}\refstepcounter{TS} 
Still under the assumptions of part ii) of the theorem
on the distributions $f_j$ and the off-line products $v_i\diam f_j$ and
$v_i\diam\partial_1^2v_j$, estimates \eqref{wi35} and \eqref{wj16} in Step \ref{nonlinear} and Step \ref{TS3} 
obviously combine to %\comment{check}
\begin{align}
M&\lesssim N_0(\bar M+1)&\quad\mbox{provided}\;\max_i[\bar u_i]\ll 1,\label{wk42}\\
\max_i M_{u_i}&\lesssim N_0\bar M&\quad\mbox{provided}\;\max_i[\bar u_i]\ll 1,\label{wk44}\\
\del M&\lesssim N_0\del\bar M+\del N_0&\quad\mbox{provided in addition}\;\bar M\le 1,\label{wk43}\\
M_{u_1-u_0}&\lesssim N_0\del\bar M+\del N_0\bar M
&\quad\mbox{provided in addition}\;\bar M\le 1.\label{wk45}
\end{align}

\medskip

%%%%%%%%%%%%%%%%%%%%%%%%%%%%%%%%%%%%%%%%%%%%%%%%%%%%%%%%%%%%%%%%%%%%%%%%%%%%%%%%%%%%%%%%%%%%%%%%%%%%%

{\sc Step} \arabic{TS}.\label{fixedpoint}\refstepcounter{TS} Contraction mapping argument.
We work under the assumptions of part ii) of the theorem
on the distributions $f_j$ and the off-line products $v_i\diam f_j,$ $v_i\diam\partial_1^2v_j$.
In this step, we specify to the case of a single model $f_1=f_0=:f$ with the corresponding constant-coefficient solution $v$;
this means that we may set $\del N_0=0$.

\medskip

We consider the space of all triplets $(\bar u,\bar a,\bar\sigma)$, where $\bar u$ is modelled after $v$ 
according to $\bar a$ and $\bar\sigma$, which fulfill the constraints (\ref{wi28}), and which satisfy
\begin{align}\label{wi30}
\bar M
\le N,
\end{align}
cf. \eqref{wi26}, 
for some constant $N$ to be fixed presently. We apply Step \ref{TS4} to $(f_i,\bar a_i,\bar\sigma_i)$
$=(f,\bar a,\bar\sigma)$. From (\ref{wi30}) and the definition (\ref{wi26}) of $\bar M$
we learn that the proviso of (\ref{wk42}) is fulfilled provided the constant $N$ is sufficiently small, 
which we now fix accordingly.
We thus learn from (\ref{wk42}), which by (\ref{wi30}) assumes the form of $M\lesssim N_0$, that the map 
defined through (\ref{wi34})
sends the set defined through (\ref{wi30}) into itself, provided $N_0\ll1$.

\medskip

For two triplets $(u_i,a_i,\sigma_i)$ as above  we first note that
\begin{align}\label{wi32}
d&\big((u_1,a_1,\sigma_1),(u_0,a_0,\sigma_0)\big):=M_{u_1-u_0}+[u_1-u_0]+\|u_1-u_0\|\nonumber\\
&+N_0([\sigma_1-\sigma_0]+\|\sigma_1-\sigma_0\|+[a_1-a_0]+\|a_1-a_0\|)
\end{align}
defines a distance function. 
Indeed, that also the modelledness constant $M_{u_1-u_0}$ satisfies a triange inequality 
in $(u_i, a_i, \sigma_i)$ can be seen by rewriting the definition~\eqref{1.3} as
\begin{align*}
&\sup_{x,R} \frac{1}{R^{2\alpha}} \inf_{\ell} \sup_{y : d(x,y) \leq R} |u_1(y) 
- \sigma_1(x) v(y,a_1(x)) \\
& \qquad \qquad \qquad \qquad - (u_0(y) - \sigma_0(x) v(y,a_0(x))) - \ell(y) |
\end{align*}
where $\ell$ runs over all linear functionals of the form $ay_1 +b$.
We now apply Step \ref{TS4} to the case of $(f_i,\bar a_i,\bar\sigma_i)=(f,\bar a_i,\bar\sigma_i)$.
From (\ref{wi30}) we learn that the proviso of (\ref{wk43}) is fulfilled; because of $\del N_0=0$,
(\ref{wk43}) assumes the form $\del M\lesssim N_0\del\bar M$. By definitions 
(\ref{wi27}) and (\ref{wk39}) of $\del\bar M$ and $\del M$, combined with $\del N_0=0$, this turns into
\begin{align*}
d\big((u_1,a_1,\sigma_1),(u_0,a_0,\sigma_0)\big)\lesssim N_0d\big((\bar u_1,\bar a_1,\bar\sigma_1),(\bar u_0,\bar a_0,\bar\sigma_0)\big).
\end{align*}
Hence the map (\ref{wi34}) is a contraction for $N_0\ll 1$.
We further note that the space of above triplets $(u,a,\sigma)$ endowed with the distance function (\ref{wi32}) is {\it complete};
and that the subset defined through the constraints (\ref{wi28}) and (\ref{wi30}) is {\it closed}. Hence by the contraction
mapping principle the map (\ref{wi34}) admits a unique fixed point on the set defined through (\ref{wi28}) and (\ref{wi30}).

\medskip

%%%%%%%%%%%%%%%%%%%%%%%%%%%%%%%%%%%%%%%%%%%%%%%%%%%%%%%%%%%%%%%%%%%%%%%%%%%%%%%%%%%%%%%%%%%%%%%%%%%%%%%%%%%%%

{\sc Step} \arabic{TS}.\label{postprocessI}\refstepcounter{TS}
Conclusion on part i) of the theorem.
Let $u$ now be as in part i) of the theorem. We note that the assumptions of part i) on the distribution
$f$ and the off-line products $v\diam f,$ $v\diam\partial_1^2v$ turn into the assumptions of part ii) with $\del N_0=0$.
We claim that $(u,a(u),\sigma(u))=:(u,a,\sigma)$ is a fixed point of the map 
(\ref{wi34}), which is obvious, that lies in the set defined through the constraints (\ref{wi28}) and (\ref{wi30}),
and therefore is unique. Indeed,
in view of $[a]\le\|a'\|[u]\le 1$, $[\sigma]\le\|\sigma'\|[u]\le 1$ by (\ref{2.15}) and (\ref{v96bis}),
the constraints (\ref{wi28}) are satisfied. The constraint (\ref{wi30}) will be an immediate consequence of the stronger
statement (\ref{v76}) (provided $N_0$ is sufficiently small). We thus turn to this a priori estimate (\ref{v76}) and apply Step \ref{TS4} to $(f_i,\bar a_i,\bar\sigma_i)=(f,a(u),\sigma(u))$. Since
we are dealing with fixed points, we have $\bar M=M$. By the theorem's assumption $[u]\ll 1$, the
provisos of (\ref{wk42}) and (\ref{wk44}) are satisfied so that because of $N_0\ll 1$, their application
yields
\begin{align}\label{wk56}
M\lesssim N_0\quad\mbox{and thus}\quad M_u\lesssim N_0^2.
\end{align}
By definition (\ref{wk40}) and the vanishing mean of $u$, this turns into (\ref{v76}).

\medskip

%%%%%%%%%%%%%%%%%%%%%%%%%%%%%%%%%%%%%%%%%%%%%%%%%%%%%%%%%%%%%%%%%%%%%%%%%%%%%%%%%%%%%%%%%%%%%%%%%%%%%%%

{\sc Step} \arabic{TS}.\label{postprocessII}\refstepcounter{TS} Conclusion on part ii) of the theorem.
Let $u_i$, $i=0,1$, now be as in part ii) of theorem. By Step \ref{postprocessI}, 
the two triplets $(u_i,a(u_i),\sigma(u_i))$ $=:(u_i,a_i,\sigma_i)$ satisfy the constraints 
(\ref{wi28}) and (\ref{wi30}) and each triplet is a fixed point of ``its own'' 
map (\ref{wi34}) (which depends on $i$ through the model $f_i$).
We apply Step \ref{TS4} to $(f_i,\bar a_i,\bar\sigma_i)=(f_i,a(u_i),\sigma(u_i))$. Since
we are dealing with fixed points, we have $\bar M=M$ and $\del\bar M=\del M$. By the
a priori estimate (\ref{v76}) and $N_0\ll1$, the two provisos of Step \ref{TS4} are satisfied. 
Because of $N_0\ll 1$, (\ref{wk43}) and (\ref{wk45}) turn into
\begin{align*}
\del M\lesssim\del N_0\quad\mbox{and then}\quad M_{u_1-u_0}\lesssim N_0\del N_0,
\end{align*}
where we used (\ref{wk56}). By definition (\ref{wk39}) of $\del M$, this turns into (\ref{wk54}) and (\ref{wk55}).

\bigskip

%%%%%%%%%%%%%%%%%%%%%%%%%%%%%%%%%%%%%%%%%%%%%%%%%%%%%%%%%%%%%%%%%%%%%%%%%%%%%%%%%%%%%%%%%%%%%%%
%%%%%%%%%%%%%%%%%%%%%%%%%%%%%%%%%%%%%%%%%%%%%%%%%%%%%%%%%%%%%%%%%%%%%%%%%%%%%%%%%%%%%%%%%%%%%%%%
%%%%%%%%%%%%%%%%%%%%%%%%%%%%%%%%%%%%%%%%%%%%%%%%%%%%%%%%%%%%%%%%%%%%%%%%%%%%%%%%%%%%%%%%%%%%%%%%
%%%%%%%%%%%%%%%%%%%%%%%%%%%%%%%%%%%%%%%%%%%%%%%%%%%%%%%%%%%%%%%%%%%%%%%%%%%%%%%%%%%%%%%%%%%%%%%%

\subsection{Proof of Proposition \ref{P}}
%\addcontentsline{toc}{section}{\protect\numberline{}Proposition \ref{P}}%
%{\sc Proof of Proposition \ref{P}}.

\newcounter{PS} % proofstep = 0
\refstepcounter{PS} % increases value by 1

We continue to abbreviate $[\cdot]=[\cdot]_\alpha$.
When a function $v$ depends on $a_0$ and $x$, we continue to write $\|v\|$ when we mean $\sup_{a_0}\|v(\cdot,a_0)\|$
and $[v]$ for $\sup_{a_0}[v(\cdot,a_0)]$. When we speak of a function $u$, we automatically mean
that it is H\"older continuous with exponent $\alpha$, that is, $[u]<\infty$;
when we speak of a distribution $f$, we imply that
it is of order $\alpha-2$ in the sense of $\|f\|_{\alpha-2}<\infty$.
When a distribution depends on the additional parameter $a_0$, we imply that the above bound
is uniform in $a_0$.

\medskip
%%%%%%%%%%%%%%%%%%%%%%%%%%%%%%%%%%%%%%%%%%%%%%%%%%%%%%%%%%%%%%%%%%%%%%%%%%%%%%%%%%%%%%%%%%%%%%%%%%%%%%%%%%%%%

{\sc Step} \arabic{PS}.\label{Unique}\refstepcounter{PS}
Uniqueness. Under the assumptions of part i) of the proposition
we claim that there is at most one mean-free $u$
modelled after $v$ according to $a$ and $\sigma$ satisfying the equation (\ref{wk07}).
Indeed, let $u'$ be another function with these properties; we trivially have by Definition \ref{D1} that
$u-u'$ is modelled after $v$ according to $a$ and to $0$ playing the role of $\sigma$.
We now apply Lemma \ref{L1} with $b$ replaced by $a$. We apply it
three times, namely
to $u$, to $u'$, and to $u-u'$. We obtain from these three versions of  (\ref{1.9})
and the triangle inequality that
\begin{align*}
\lim_{T\downarrow0}\|[a,(\cdot)_T]\diam\partial_1^2u
-[a,(\cdot)_T]\diam\partial_1^2u'-[a,(\cdot)_T]\diam\partial_1^2(u-u')\|=0
\end{align*}
and thus
$\lim_{T\downarrow0}\|(a\diam\partial_1^2u
-a\diam\partial_1^2u'-a\diam\partial_1^2(u-u'))_T\|=0$
so that $a\diam\partial_1^2u$ $-a\diam\partial_1^2u'$ $=a\diam\partial_1^2(u-u')$.
Hence  we obtain from taking the difference of the equations:
\begin{align}\label{x85}
\partial_2(u-u')-P a\diam\partial_1^2(u-u')=0.
\end{align}
We may also say that $u-u'$ is modelled after $0$ playing the role of $v$
and $0$ playing the role of $\sigma$; we call $\del M$ the corresponding modelling
constant. Hence we may apply Corollary \ref{C1} i)
with $f=0$ and thus $N_0=0$. We apply it with $u$ replaced by $u-u'$,
which we may thanks to (\ref{x85}).
In this context, the output (\ref{wk02}) of Corollary \ref{C1} assumes the form $[u-u']=0$.
%Since $u-u'$ is periodic, we first infer $\del\nu=0$ and then $u-u'=const$.
Since $u-u'$ has vanishing average, we obtain as desired $u-u'=0$.

\medskip
%%%%%%%%%%%%%%%%%%%%%%%%%%%%%%%%%%%%%%%%%%%%%%%%%%%%%%%%%%%%%%%%%%%%%%%%%%%%%%%%%%%%%%%%%%%%%%%%%%%%%%%%%%%%%%%%

{\sc Step} \arabic{PS}.\label{Regul}\refstepcounter{PS}
A special regularization. Under the assumptions of Lemma \ref{L1} and for $\tau>0$ and $i=1,\cdots,I$
we consider the convolution $v_{i\tau}$ of $v_i$ and {\it define}
\begin{align}\label{y15}
a\diam\partial_1^2 v_{i\tau}:=(a\diam \partial_1^2v_i)_\tau.
\end{align}
Then, we claim that for any function $u$ of class $C^{\alpha+2}$, which is
modelled after $(v_{1\tau},\cdots,v_{I\tau})$ according to $a$ and $(\sigma_1,\cdots,\sigma_I)$, we have
\begin{align}\label{x87}
a\diam\partial_1^2u=a\partial_1^2u-\sigma_i E[a,(\cdot)_\tau]\diam\partial_1^2v_i.
\end{align}

\medskip

Indeed, by Lemma \ref{L1} (with $b$ replaced by $a$) we understand the distribution $a\diam\partial_1^2u$ as defined by
\begin{align}\label{y16}
\lim_{T\downarrow0}\|[a,(\cdot)_T]\diam\partial_1^2u-\sigma_i E[a,(\cdot)_T]\diam\partial_1^2v_{i\tau}\|=0.
\end{align}
We note that (\ref{y15}) implies by the semi-group property
\begin{align}\label{hh1}
[a,(\cdot)_T]\diam\partial_1^2v_{i\tau}
=[a,(\cdot)_{T+\tau}]\diam\partial_1^2v_i,
\end{align}
which ensures that $[a,(\cdot)_T]\diam\partial_1^2v_{i\tau}\rightarrow[a,(\cdot)_\tau]\diam\partial_1^2v_i$
as $T\downarrow0$ uniformly in $x$ for fixed $a_0$. Thanks to the bound on the
$\frac{\partial}{\partial a_0}$-derivative in (\ref{1.5}),
this convergence is even uniform in $(x,a_0)$, so that (\ref{y16}) turns into
\begin{align*}
\lim_{T\downarrow0}\|[a,(\cdot)_T]\diam\partial_1^2u-\sigma_i E[a,(\cdot)_\tau]\diam\partial_1^2v_i\|=0.
\end{align*}
Since $u$ is of class $C^{\alpha+2}$, this further simplifies to
\begin{align*}
\lim_{T\downarrow0}\|a\partial_1^2u-(a\diam\partial_1^2u)_T-\sigma_i E[a,(\cdot)_\tau]\diam\partial_1^2v_i\|=0,
\end{align*}
from which we learn that the distribution $a\diam\partial_1^2u$ is actually the function given by
(\ref{x87}).

\medskip
%%%%%%%%%%%%%%%%%%%%%%%%%%%%%%%%%%%%%%%%%%%%%%%%%%%%%%%%%%%%%%%%%%%%%%%%%%%%%%%%%%%%%%%%%%%%%%%%%%%%%%%%%%%%%%

{\sc Step} \arabic{PS}. \label{Exist} \refstepcounter{PS} Existence in the regularized case.
Under the assumptions of part i) of this proposition and in line with Step \ref{Regul},
for $\tau>0$ we consider the mollification $f_\tau$ of $f$,
so that $v_\tau$ satisfies $(\partial_2-a_0\partial_1^2)v_\tau=Pf_\tau$, and complement definition (\ref{y15})
(without the index $i$) by
\begin{align}\label{wk08}
\sigma\diam f_{\tau}:=(\sigma\diam f)_\tau.
\end{align}
Then we claim that there exists a mean-free $u^\tau$ of class $C^{\alpha+2}$ modelled
after $v_\tau$ according to $a$ and $\sigma$ such that
\begin{align}\label{y14}
\partial_2u^\tau-P(a\diam\partial_1^2 u^\tau+\sigma\diam f_\tau)=0\quad\mbox{distributionally},
\end{align}
and at the same time
\begin{align}\label{y13}
\partial_2u^\tau-P(a\partial_1^2u^\tau-\sigma E[a,(\cdot)_\tau]\diam\partial_1^2v
+(\sigma\diam f)_\tau)=0\quad\mbox{classically}.
\end{align}

\medskip

We first turn to the existence of (\ref{y13}) and start by noting that the right hand side
$-\sigma E[a,(\cdot)_\tau]\diam\partial_1^2v+(\sigma\diam f)_\tau$ in (\ref{y13}) is of class $C^\alpha$.
Leveraging upon $[a]\ll 1$ we rewrite the equation as
$\partial_2u^\tau-a_0\partial_1^2u^\tau$ $=P((a-a_0)\partial_1^2u$ $-\sigma E[a,(\cdot)_\tau]\diam\partial_1^2v$
$+(\sigma\diam f)_\tau)$
for $a_0=a(0)$.
Using the invertibility of the constant-coefficient operator $\partial_2-a_0\partial_1^2$ on periodic mean-free
functions, and equipped with the corresponding Schauder estimates, see for instance \cite[Theorem 8.6.1]{Krylov}
lifted to the torus, we see that a solution of class $C^{\alpha+2}$ exists,
using a contraction mapping argument based on $\|a-a_0\|\ll 1$. Since both $u^\tau$ and $v_\tau(\cdot,a_0)$ are in particular of
class $C^{\alpha+1}$, $u$ is modelled after $v_\tau$ according to --- in fact any --- $a$ and $\sigma$.
By Step \ref{Regul} and definition (\ref{wk08}) we see that (\ref{y13}) may be rewritten as (\ref{y14}).

\medskip
%%%%%%%%%%%%%%%%%%%%%%%%%%%%%%%%%%%%%%%%%%%%%%%%%%%%%%%%%%%%%%%%%%%%%%%%%%%%%%%%%%%%%%%%%%%%%%%%%%%%%%

{\sc Step} \arabic{PS}.\label{Approx}\refstepcounter{PS} Basic construction.
We now work under the assumptions of part ii) of the proposition.
We interpolate the functions $\sigma_i$, $a_i$, and $v_i$ as well as the distribution $f_i$ linearly:
\begin{align}\label{y17}
\sigma_s:=s\sigma_1+(1-s)\sigma_0\quad\mbox{and the same for $a$, $f$, and $v$}.
\end{align}
We note that this preserves  \eqref{wj59}. We interpolate the products bi-linearly  
\begin{align}\label{y28}
\sigma_{s}\diam f_{s}&:=s^2\sigma_1\diam f_1+s(1-s)\sigma_1\diam f_0\nonumber\\
&+(1-s)s\sigma_0\diam f_1+(1-s)^2\sigma_0\diam f_0,\nonumber\\
\partial_{s}\sigma  \diam f_{s}&:=s\sigma_1\diam f_1+(1-s)\sigma_1\diam f_0
-s\sigma_0\diam f_1-(1-s)\sigma_0\diam f_0,\nonumber\\
\sigma_{s}\diam\partial_{s}f&:=s\sigma_1\diam f_1-s\sigma_1\diam f_0
+(1-s)\sigma_0\diam f_1-(1-s)\sigma_0\diam f_0,\nonumber\\
&\mbox{and the same for}\;a_{s}\diam\partial_1^2v_s,\;
\partial_{s}a\diam\partial_1^2v_s\;\mbox{and}\;a_{s}\diam\partial_1^2\partial_sv,
\end{align}
where here and below we use the convention that $\partial_s$ only acts on the object
directly following it (with argument suppressed), i.e. for example $\partial_{s}\sigma  \diam f_{s}= (\partial_s \sigma_s)\diam f_s$.

Thanks to the estimate (\ref{wj98bis}), which is preserved under bilinear interpolation, the
family of distributions $\{a_s\diam\partial_1^2v_s(\cdot,a_0)\}_{a_0}$ is continuously
differentiable in $a_0$ so that we may define
\begin{align}\label{v82}
a_s\diam\partial_1^2\frac{\partial v_s}{\partial a_0}(\cdot,a_0):=
\frac{\partial}{\partial a_0}a_s\diam\partial_1^2v_s(\cdot,a_0).
\end{align}
For given $0<\tau\le 1$, we define the singular products with the regularized distributions as in Step \ref{Regul},
namely
\begin{align}\label{y21}
\sigma_s\diam f_{s\tau}&:=(\sigma_s\diam f_s)_\tau\quad\mbox{and the same for}\;\partial_s\sigma\diam f_{s\tau},\;
\sigma_s\diam\partial_sf_\tau,\nonumber\\
&a_s\diam\partial_1^2v_{s\tau},\;
\partial_sa\diam\partial_1^2v_{s\tau},\;a_s\diam\partial_1^2\partial_sv_\tau,\;a_s\diam\partial_1^2\frac{\partial v_{s\tau}}{\partial a_0}.
\end{align}
We claim that there exists a curve $u_s^\tau$ of mean-free functions continuously differentiable in $s$
with respect to the class $C^{\alpha+2}$ such that
\begin{align}
u_s^\tau\;\mbox{is modelled after}\;v_{s\tau}\;\mbox{according to}\;a_s\;\mbox{and}\;\sigma_s,\label{y26}\\
\partial_2u_s^\tau-P(a_s\diam\partial_1^2u_s^\tau+\sigma_s\diam f_{s\tau})=0\quad
\mbox{distributionally}.\label{y25}
\end{align}
Furthermore, we claim that
\begin{align}
\partial_su^\tau\;\mbox{is modelled after}
\;(v_{s\tau},\frac{\partial v_{s\tau}}{\partial a_0},\partial_sv_\tau)
\;\mbox{according to}\;a_s\;\mbox{and}\;(\partial_s\sigma,\sigma_s\partial_s a,\sigma_s),\label{y29}\\
\partial_2\partial_su^\tau-P(a_s\diam\partial_1^2\partial_su^\tau+\partial_sa\diam\partial_1^2u_s^\tau+
\sigma_s\diam \partial_sf_{\tau}+\partial_s\sigma\diam f_{s\tau})=0\label{y27}
\end{align}
distributionally.
Note that (\ref{y27}) is what we get from formally applying $\partial_s$ to (\ref{y25}).

\medskip

Here comes the argument: By Steps \ref{Exist} and \ref{Unique} and our definitions of $\sigma_s\diam f_{s\tau}$
and $a_s\diam\partial_1^2v_{s\tau}$ by convolution, cf. (\ref{y21}),
there exists a unique mean-free $u_s^\tau$ of class $C^{\alpha+2}$
such that (\ref{y26}) and (\ref{y25}) hold.
Furthermore by Step \ref{Regul} $u_s^\tau$ is characterized as the classical solution of
\begin{align}\label{y24}
\partial_2u_s^\tau-P(a_s\partial_1^2u_s^\tau-\sigma_sE_s[a_s,(\cdot)_\tau]\diam\partial_1^2v_s
+(\sigma_s\diam f_{s})_\tau)=0.
\end{align}
In preparation of taking the $s$-derivative of (\ref{y24}) we note that
the definition (\ref{y28}) of $\sigma_s\diam f_s$ and $a_s\diam\partial_1^2v_s$ by (bi-)linear interpolation
ensures that Leibniz's rule holds:
\begin{align}
\partial_s(\sigma_s\diam f_s)&=\partial_s\sigma \diam f_s+\sigma_s\diam\partial_sf,\label{v92}\\
\partial_s(a_s\diam\partial_1^2v_s)&=\partial_s a  \diam
\partial_1^2v_s+a_s\diam\partial_1^2\partial_sv.\label{v81}
\end{align}
We recall that $E_s$ denotes the evaluation operator that evaluates a function
of $(x,a_0)$ at $(x,a_s(x))$; with the obvious commutation rule
$[\partial_s,E_s]=\partial_sa E_s\frac{\partial}{\partial a_0}$ we obtain from (\ref{v81})
and (\ref{v82}) 
\begin{align*}
\lefteqn{\partial_s(E_s a_s\diam\partial_1^2v_s)}\nonumber\\
&=E_s \partial_s a\diam\partial_1^2v_s+
\partial_s a E_s a_s\diam\partial_1^2\frac{\partial v_s}{\partial a_0}
+E_s a_s\diam\partial_1^2\partial_sv,
\end{align*}
which in conjunction with the classical differentiation rules extends to the commutator:
\begin{align}\label{v83}
\lefteqn{\partial_s(E_s[a_s,(\cdot)_\tau]\diam\partial_1^2v_s)
=E_s[\partial_s a,(\cdot)_\tau]\diam\partial_1^2v_s}\nonumber\\
&+\partial_s a E_s [a_s,(\cdot)_\tau]\diam\partial_1^2\frac{\partial v_s}{\partial a_0}+E_s [a_s,(\cdot)_\tau]\diam\partial_1^2\partial_sv.
\end{align}
Equipped with (\ref{v92}), (\ref{v81}), and (\ref{v83}) we learn from (\ref{y24})
by the argument of Step \ref{Exist} that $u_s^\tau$ is differentiable in $s$
with values in the class $C^{\alpha+2}$ and
\begin{align}\label{x53}
\lefteqn{\partial_2\partial_su^\tau
-P\big(a_s\partial_1^2\partial_su^\tau+\partial_sa\partial_1^2u^\tau_s-\sigma_sE_s[\partial_s a,(\cdot)_\tau]\diam\partial_1^2v_s}\nonumber\\
&-\partial_s\sigma E_s[a_s,(\cdot)_\tau]\diam\partial_1^2v_s
-\sigma_s\partial_sa E_s[a_s,(\cdot)_\tau]\diam\partial_1^2\frac{\partial v_s}{\partial a_0}\nonumber\\
&-\sigma_s E_s[a_s,(\cdot)_\tau]\diam\partial_1^2\partial_sv
+(\partial_s\sigma\diam f_{s})_\tau+(\sigma_s\diam \partial_s f)_{\tau}\big)=0.
\end{align}
Moreover, like in Step \ref{Exist}, (\ref{y29}) holds automatically because of the regularity
of $\partial_su^\tau$ and of $(v_{s\tau},\frac{\partial v_{s\tau}}{\partial a_0},\partial_sv_\tau)$.
In view of the definition (\ref{y21}) of $\partial_sa\diam\partial_1^2v_{s\tau}$
we have by Step \ref{Regul} applied to $u_s^\tau$ modelled according to
(\ref{y26})
\begin{align*}
\partial_sa\diam\partial_1^2u_s^\tau&
=\partial_sa\partial_1^2u_s^\tau-\sigma_sE_s[\partial_s a,(\cdot)_\tau]\diam\partial_1^2v_s.
\end{align*}
In view of the similar definition of $a_s\diam\partial_1^2\partial_sv$,
$a_s\diam\partial_1^2\frac{\partial v_{s\tau}}{\partial a_0}$, and $a_s\diam\partial_1^2\partial_sv_\tau$
we have by Step \ref{Regul} applied to $\partial_su^\tau$ modelled according to
(\ref{y29})
\begin{align*}
a_s\diam\partial_1^2\partial_su^\tau&=a_s\partial_1^2\partial_su^\tau
-\partial_s\sigma E_s[a_s,(\cdot)_\tau]\diam\partial_1^2v_s&\nonumber\\
&-\sigma_s\partial_sa E_s[a_s,(\cdot)_\tau]\diam\partial_1^2\frac{\partial v_s}{\partial a_0}
-\sigma_s E_s[a_s,(\cdot)_\tau]\diam\partial_1^2\partial_sv.
\end{align*}
Plugging these two formulas and the definition (\ref{y21}) of $\partial_s\sigma\diam f_\tau$
and $\sigma_s\diam\partial_sf_\tau$ into (\ref{x53}), we obtain (\ref{y27}).

\medskip
%%%%%%%%%%%%%%%%%%%%%%%%%%%%%%%%%%%%%%%%%%%%%%%%%%%%%%%%%%%%%%%%%%%%%%%%%%%%%%%%%%%%%%%%%%%%%%%%%%%%

{\sc Step} \arabic{PS}.\label{Interpol}\refstepcounter{PS}
We still work under the assumptions of part ii) of the proposition. We claim 
\begin{align}
[\sigma_s]+[a_s]&\le N,\label{wk52}\\
\|f_{s \tau}\|_{\alpha-2} &\lesssim N_0,\label{wk10}\\
\|[\sigma_s,(\cdot)_T]\diam f_{s\tau}\|_{2\alpha-2}&\lesssim N N_0,\label{wk11}\\
\|
[a_s,(\cdot)]\diam\partial_1^2v_{s\tau}\|_{2\alpha-2,2}&\lesssim N N_0\label{wk12}
\end{align}
and on the corresponding estimates on the infinitesimal level
\begin{align}
[\partial_s\sigma_s]+\|\partial_s\sigma_s\|+[\partial_sa]+\|\partial_sa\|&\le \del N,\label{wk53}\\
\|\partial_s f_\tau\|_{\alpha-2}&\le \del N_0,\label{wk19}\\
\|[\sigma_s,(\cdot)]\diam \partial_sf_\tau\|_{2\alpha-2}&\lesssim N \del N_0,\label{wk20}\\
%\sup_{T\le 1}(T^\frac{1}{4})^{2-2\alpha}
\|[\partial_s\sigma,(\cdot)]\diam f_{s\tau}\|_{2\alpha-2}&\lesssim \del N N_0,\label{wk21}\\
%
%\sup_{T\le 1}(T^\frac{1}{4})^{2-2\alpha}\sup_{a_0}
\|
%\{1,\frac{\partial}{\partial a_0}\}
[a_s,(\cdot)]\diam\partial_1^2 \partial_sv_\tau\|_{\alpha-2,1}
&\lesssim N \del N_0,\label{wk22}\\
%\sup_{T\le 1}(T^\frac{1}{4})^{2-2\alpha}\sup_{a_0}
\|
[\partial_s a_s,(\cdot)]\diam\partial_1^2v_{s\tau}\|_{2\alpha-2, 1}
&\lesssim \del N N_0.\label{wk23}
\end{align}
Indeed, (\ref{wk52}) and (\ref{wk53}) are immediate from our assumptions (\ref{wk59}) (with $i$) and (\ref{wk50}), respectively,
by the linear interpolation (\ref{y17}). For $\tau=0$ the remaining estimates, even with $\lesssim$ replaced by $\le$, follow from
the linear and bilinear interpolations (\ref{y17}) and (\ref{y28}) from the
assumptions of this proposition: inequality (\ref{wk10}) from (\ref{wj58}) (with $i$), (\ref{wk11}) from (\ref{wj97bis}),
(\ref{wk12}) from (\ref{wj98bis}).
Still for $\tau=0$, the five estimates  (\ref{wk19}), (\ref{wk20}), (\ref{wk21}), (\ref{wk22}), and (\ref{wk23}),
are direct consequences of  (\ref{wk24}), (\ref{wk13}), (\ref{wk25}), (\ref{wk14}), and (\ref{w26}), respectively.

\medskip

It remains to pass from $\tau=0$ to $0<\tau\le 1$ in the eight estimates of this step, based on our
definition (\ref{y21}) of singular products. This is done with help of the next step.

\medskip

%%%%%%%%%%%%%%%%%%%%%%%%%%%%%%%%%%%%%%%%%%%%%%%%%%%%%%%%%%%%%%%%%%%%%%%%%%%%%%%%%%%%%%%%%%%%%%%%%%%%%%%
{\sc Step} \arabic{PS}.\refstepcounter{PS}
Let the (generic) function $u$ and the (generic) distributions $f$ and $u\diam f$ be such that
\begin{align}
[u]\le N_0,\quad \|f\|_{\alpha-2}\le N_1\quad\mbox{and}\quad
%\sup_{T\le 1}(T^\frac{1}{4})^{2-2\alpha}
\|[u,(\cdot)]\diam f\|_{2\alpha-2}\le N_0N_1\label{y32}
\end{align}
for some constants $N_0$ and $N_1$. Then we claim that for $\tau\le 1$ the distributions $f_\tau$ and
$u\diam f_\tau:=(u\diam f)_\tau$ satisfy the same estimates:
\begin{align}\label{y31}
\| f_\tau\|_{\alpha-2}\lesssim N_1\quad\mbox{and}\quad\|[u,(\cdot)]\diam f_\tau\|_{2\alpha-2}\lesssim N_0N_1.
\end{align}
Indeed, by definition of $u\diam f_\tau$ we have like for \eqref{hh1} 
\begin{align*}
[u,(\cdot)_T]\diam f_\tau=[u,(\cdot)_{T+\tau}]\diam f,
\end{align*}
so that (\ref{y31}) follows automatically provided we can show that (\ref{y32})
extend from the range of $T\le 1$ to the range $T\le 2$ in form of
\begin{align}\label{y33}
\sup_{T\le 1}(T^\frac{1}{4})^{2-\alpha}\|f_{2T}\|\lesssim N_1,
\quad\sup_{T\le 1}(T^\frac{1}{4})^{2-2\alpha}\|[u,(\cdot)_{2T}]\diam f\|\lesssim N_0N_1.
\end{align}
For this, we appeal to the semi-group property giving us
\begin{align*}
f_{2T}=(f_T)_T\quad\mbox{and}\quad
[u,(\cdot)_{2T}]\diam f=([u,(\cdot)_T]\diam f)_T+[u,(\cdot)_T]f_T,
\end{align*}
so that by the boundedness of $(\cdot)_T$ in $\|\cdot\|$ indeed the last item in (\ref{y32})
entails (\ref{y33}), appealing to (\ref{wi2}) and using in addition that by the first items in (\ref{y32})
\begin{align*}
\|[u,(\cdot)_T]f_T\|\lesssim N_0(T^\frac{1}{4})^\alpha\|f_T\|\lesssim N_0N_1(T^\frac{1}{4})^{2\alpha-2}.
\end{align*}

\medskip
%%%%%%%%%%%%%%%%%%%%%%%%%%%%%%%%%%%%%%%%%%%%%%%%%%%%%%%%%%%%%%%%%%%%%%%%%%%%%%%%%%%%%%%%%%%%%%%%%%%%%%%%%%%%

{\sc Step} \arabic{PS}.\refstepcounter{PS} Application of Corollary \ref{C1}. We claim
for the modelling and H\"older constants of $u^\tau_s$ and $\partial_su^\tau$:
\begin{align}
M_s^\tau&\lesssim N_0N,\label{wk27}\\
[u_s^\tau]&\lesssim N_0(N+1),\label{wk28}\\\
\del M_s^\tau&\lesssim N_0\del N+\del N_0N\quad\mbox{provided}\;N\le 1,\label{wk29}\\
[\partial_su^\tau]&\lesssim N_0\del N+\del N_0\quad\mbox{provided}\;N\le 1.\label{wk30}
\end{align}
Indeed, for estimates (\ref{wk27}) and (\ref{wk28}) we apply Corollary \ref{C1} i) with
$(f,v,$ $\sigma,a,$ $\sigma\diam f,a\diam\partial_1^2v,$ $u)$ replaced by
$(f_{s\tau},v_{s\tau},$ $\sigma_s,a_s,$ $\sigma_s\diam f_{s\tau},a_s\diam\partial_1^2v_{s\tau},$ $u_s^\tau)$
(where it is clear that linear interpolation and convolution preserves the relation between $f_{s\tau}$ and
$v_{s\tau}$ through the constant coefficient equation).
As already remarked in Step \ref{Approx} 
the linear interpolation (\ref{y17}) preserves \eqref{wj59}. The hypotheses (\ref{wk59}), (\ref{wj58}), (\ref{wj56}), and 
(\ref{wj57})
were established in Step \ref{Interpol}, cf. (\ref{wk52}), (\ref{wk10}), (\ref{wk11}), and (\ref{wk12}), respectively. 
Hypothesis (\ref{wj55}) and the modelledness are clear by construction, cf. (\ref{y25}) and (\ref{y26}) in Step \ref{Approx}.
The  outputs (\ref{wj63}) and (\ref{wk02})  assume the form (\ref{wk27}) and
(\ref{wk28}).

\medskip

For the remaining estimates (\ref{wk29}) and (\ref{wk30}), we apply Corollary \ref{C1} ii)
with $(\del f,\del v,$ $\del \sigma,\del a,$ $\sigma\diam \del f,\del\sigma\diam f,$ $a\diam\partial_1^2\del v,\del a\diam\partial_1^2v,$
$\del u)$ replaced by
$(\partial_sf_{\tau},\partial_sv_{\tau},$ $\partial\sigma,\partial_sa,$ $\sigma_s\diam \partial_sf_{\tau},
\partial_s\sigma\diam f_{s\tau},$ $a_s\diam\partial_1^2\partial_sv_{\tau},\partial_sa\diam\partial_1^2v_{s\tau},$ $\partial_s u^\tau)$.
The six hypotheses (\ref{wk60})--% 
(\ref{wj92})
were established in Step \ref{Interpol}, cf. (\ref{wk53})--
 (\ref{wk23}).
Hypothesis (\ref{wj66}) and the corresponding modelledness are clear by construction, cf. (\ref{y27}) and (\ref{y29})
in Step \ref{Approx}. The  outputs (\ref{wj80}) and (\ref{wk03}) assume the form of (\ref{wk29}) and (\ref{wk30}).

\medskip
%%%%%%%%%%%%%%%%%%%%%%%%%%%%%%%%%%%%%%%%%%%%%%%%%%%%%%%%%%%%%%%%%%%%%%%%%%%%%%%%%%%%%%%%%%%%%%%%%%%%%%%%%%%%

{\sc Step} \arabic{PS}.\refstepcounter{PS} Integration. We claim that $u_1^\tau-u_0^\tau$ is modelled after
$(v_1^\tau,v_0^\tau)$ according to $(a_1,a_0)$ and $(\sigma_1,-\sigma_0)$ with the modelling constant and H\"older constant 
estimated as follows
\begin{align}
\del M^\tau&\lesssim N_0\del N+\del N_0N\quad\mbox{provided}\;N\le 1,\label{wk31}\\
[u_1^\tau-u_0^\tau]&\lesssim N_0\del N+\del N_0\quad\mbox{provided}\;N\le 1.\label{wk32}
\end{align}
Indeed, the H\"older estimate (\ref{wk32}) is obvious from (\ref{wk30}) by integration in $s\in[0,1]$.
The estimate on the modelling constant relies on the differentiation rule
\begin{align*}
\frac{\partial}{\partial s}&\big(u_s^\tau(y)-\sigma_s(x)v_{s\tau}(y,a_s(x))\big)
=\partial_su^\tau(y)-(\partial_s\sigma)(x) v_{s\tau}(y,a_s(x))\nonumber\\
&-(\sigma_s\partial_sa)(x)\frac{\partial v_{s\tau}}{\partial a_0}(y,a_s(x))
-\sigma_s(x)\partial_sv_\tau(y,a_s(x)),
\end{align*}
and on defining $\nu:=\int_0^1\nu_sds$, where $\nu$ belongs to $u_1^\tau-u_0^\tau$ and $\nu_s$ to $\partial_su^\tau$
in the sense of Definition \ref{D1}. This provides the link between (\ref{wk29}) and (\ref{wk31}) by
integration.

\medskip
%%%%%%%%%%%%%%%%%%%%%%%%%%%%%%%%%%%%%%%%%%%%%%%%%%%%%%%%%%%%%%%%%%%%%%%%%%%%%%%%%%%%%%%%%%%%%%%%%%%%%%%%%%%%

{\sc Step} \arabic{PS}.\refstepcounter{PS} Passage to limit.
We claim that we may pass to the limit $\tau\downarrow0$ in (\ref{wk27}) and (\ref{wk28}) with $s=0,1$,
recovering (\ref{wk17}) and (\ref{wk18}) in part i) of this proposition,
and in (\ref{wk31}) and (\ref{wk32}), recovering (\ref{wk15}) and (\ref{wk16}) in part ii) of the proposition.
Clearly, from the uniform-in-$\tau$ estimate (\ref{wk28}) (in conjunction with the vanishing mean of $u_i^\tau$
which provides the same bound on the supremum norm) we learn by Arzel\`a-Ascoli that there exists a subsequence
$\tau\downarrow 0$ (unchanged notation) and a continuous mean-free function $u_i$ to which $u_i^\tau$ converges
uniformly. Hence we may pass to the limit in the H\"older estimates (\ref{wk28}) and
(\ref{wk32}). Since also the convolution $v_{i\tau}$ converges to $v_i$ uniformly, we may pass to the limit
in the estimates (\ref{wk27}) and (\ref{wk31}) of the modelling constants. 
By uniqueness, cf. Step \ref{Unique}, it thus remains to
argue that $u_i$ solves (\ref{wk07}) (with $(f,\sigma,a)$ replaced by $(f_i,\sigma_i,a_i)$).
In order to pass from (\ref{y25}) to (\ref{wk07}) it remains to establish the distributional convergences
\begin{align}
\sigma_i\diam f_{i\tau}&\rightharpoonup \sigma_i\diam f_i,\label{y20}\\
a_i\diam\partial_1^2u_i^\tau&\rightharpoonup 
a_i\diam\partial_1^2u_i.\label{y22}
\end{align}
The convergence (\ref{y20}) is built-in by the definition (\ref{y21}) through convolution.
One of the ingredients for the convergence (\ref{y22}) is the analogue of (\ref{y20})
\begin{align*}
a_i\diam\partial_1^2v_{i\tau}(\cdot,a_0)
\rightharpoonup a_i\diam\partial_1^2v_{i}(\cdot,a_0),
\end{align*}
which in conjunction with the uniform convergence of $v_{i\tau}$ extends to the commutator
\begin{align*}
[a_i,(\cdot)_T]\diam\partial_1^2v_{i\tau}(\cdot,a_0)
\rightarrow[a_i,(\cdot)_T]\diam\partial_1^2v_{i}(\cdot,a_0).
\end{align*}
Since $\sup_{a_0}\|\frac{\partial}{\partial a_0}[a_i,(\cdot)_T]\diam\partial_1^2v_{i\tau}(\cdot,a_0)\|$
is uniformly bounded, cf. (\ref{wj98bis}) and (\ref{y21}) in conjunction with a formula of type \eqref{hh1}, we even have
\begin{align*}
[a_i,(\cdot)_T]\diam\partial_1^2v_{i\tau}(\cdot,a_0)
\rightarrow [a_i,(\cdot)_T]\diam\partial_1^2v_i(\cdot,a_0)\quad\mbox{uniformly in}\;a_0,
\end{align*}
so that
\begin{align*}
\sigma_iE_i[a_i,(\cdot)_T]\diam\partial_1^2v_{i\tau}
\rightarrow \sigma_i E_i[a_i,(\cdot)_T]\diam\partial_1^2v_i.
\end{align*}
In order to relate this to (\ref{y22}) we appeal to the modelledness of $u_i$
with respect to  $v_i$ according to $a_i$ and $\sigma_i$ which by (\ref{1.9}) in Lemma \ref{L1} yields
\begin{align*}
\lim_{T\downarrow0}\|[a_i,(\cdot)_T]\diam\partial_1^2u_i-\sigma_i E_i[a_i,(\cdot)_T]\diam\partial_1^2v_i\|=0.
\end{align*}
Likewise, the uniform modelledness of $u_i^{\tau}$, cf. (\ref{wk27}), in conjunction with
the uniform commutator bounds (\ref{wj57}) and the uniform bounds on $v_{i\tau}$,
we have, again by (\ref{1.9}) in Lemma \ref{L1}, the uniform convergence
\begin{align*}
\lim_{T\downarrow0}\sup_{\tau}
\|[a_i,(\cdot)_T]\diam\partial_1^2u_i^\tau-\sigma_{i}E_i[a_i,(\cdot)_T]\diam\partial_1^2v_{i\tau}\|=0.
\end{align*}
The combination of the three last statements implies
\begin{align*}
\lim_{T\downarrow0}\limsup_{\tau\downarrow 0}
\|[a_i,(\cdot)_T]\diam\partial_1^2u_i^\tau-[a_i,(\cdot)_T]\diam\partial_1^2u_i\|=0,
\end{align*}
which by the convergence of $u_i^\tau$ yields
\begin{align}\label{y07}
\lim_{T\downarrow 0}\limsup_{\tau\downarrow 0}\|(a_i\diam\partial_1^2u_i^\tau- a_i\diam\partial_1^2u_i)_T\|=0.
\end{align}
Now the next step shows that this implies (\ref{y22}).

\medskip
%%%%%%%%%%%%%%%%%%%%%%%%%%%%%%%%%%%%%%%%%%%%%%%%%%%%%%%%%%%%%%%%%%%%%%%%%%%%%%%%%%%%%%%%%%%%%%%%%%%%%%%%%%%%

{\sc Step} \arabic{PS}.\label{Aux}\refstepcounter{PS}
Let a sequence of distributions $\{f_n\}_{n\uparrow\infty}$ be bounded wrt $\|\cdot\|_{\alpha-2}$; then we
claim
\begin{align*}
\lim_{T\downarrow 0}\limsup_{n \uparrow \infty}\|f_{nT}\|=0\quad\Longrightarrow\quad
f_{n}\rightharpoonup 0.
\end{align*}
Indeed, we have for fixed $T>0$
and any $\tau\le T$ that $\|f_{n T}\|\lesssim\|f_{n\tau}\|$ and therefore
$\limsup_{n \uparrow \infty}\|f_{n T}\|\lesssim\limsup_{n \uparrow\infty}\|f_{n\tau}\|$
and  $\limsup_{n\uparrow\infty}\|f_{n T}\|
\lesssim\lim_{\tau\downarrow0}\limsup_{n \uparrow\infty}\|f_{n\tau}\|$. The latter is
equal to zero by assumption. Hence we have $f_{nT}\rightarrow 0$
for every $T>0$, which yields the claim by the boundedness of $f_n$ wrt $\|\cdot\|_{\alpha-2}$,
and then also in the more classical $C^{\alpha -2}$-norm, cf. \eqref{w30} in Step~\ref{LA1S1}
of Lemma~\ref{LA1}. 

\bigskip

%%%%%%%%%%%%%%%%%%%%%%%%%%%%%%%%%%%%%%%%%%%%%%%%%%%%%%%%%%%%%%%%%%%%%%%%%%%%%%%%%%%%%%%%%%%%%%%%%%%%%%%
%%%%%%%%%%%%%%%%%%%%%%%%%%%%%%%%%%%%%%%%%%%%%%%%%%%%%%%%%%%%%%%%%%%%%%%%%%%%%%%%%%%%%%%%%%%%%%%%%%%%%%%
%%%%%%%%%%%%%%%%%%%%%%%%%%%%%%%%%%%%%%%%%%%%%%%%%%%%%%%%%%%%%%%%%%%%%%%%%%%%%%%%%%%%%%%%%%%%%%%%%%%%%%%
%%%%%%%%%%%%%%%%%%%%%%%%%%%%%%%%%%%%%%%%%%%%%%%%%%%%%%%%%%%%%%%%%%%%%%%%%%%%%%%%%%%%%%%%%%%%%%%%%%%%%%%
%%%%%%%%%%%%%%%%%%%%%%%%%%%%%%%%%%%%%%%%%%%%%%%%%%%%%%%%%%%%%%%%%%%%%%%%%%%%%%%%%%%%%%%%%%%%%%%%%%%%%%%
%%%%%%%%%%%%%%%%%%%%%%%%%%%%%%%%%%%%%%%%%%%%%%%%%%%%%%%%%%%%%%%%%%%%%%%%%%%%%%%%%%%%%%%%%%%%%%%%%%%%%%%

\subsection{Proof of Corollary \ref{C1}}
%\addcontentsline{toc}{section}{\protect\numberline{}Corollary \ref{C1}}

We write $[\cdot]$ for $[\cdot]_\alpha$.

\newcounter{C1S} % proofstep = 0
\refstepcounter{C1S} % increases value by 1

{\sc Step} \arabic{C1S}.\label{C1S0}\refstepcounter{C1S} Application of Lemma \ref{LA1}.
We claim
\begin{align}
[v]_2\lesssim N_0,\label{wj75}\\
[\del v]_1\lesssim \del N_0,\label{wj77}
\end{align}
where we recall the notational convention (\ref{conv1}) for the $a_0$-derivatives.
The estimate (\ref{wj75})
is based on the identities following from differentiating (\ref{e:v-eq})
twice with respect to $a_0$
\begin{align}\label{x90}
(\partial_2-a_0\partial_1^2)\big\{v,\frac{\partial v}{\partial a_0},\frac{\partial^2 v}{\partial a_0^2}\big\}
=\big\{Pf,\partial_1^2v,2\partial_1^2\frac{\partial v}{\partial a_0}\big\}.
\end{align}
We now see that (\ref{wj75})  follows by an iterated application of Lemma \ref{LA1}:
From (\ref{wj58}) we first obtain the bound on $[v]$ by Lemma \ref{LA1}, then the bound on $\|\partial_1^2v\|_{\alpha-2}$
by (\ref{1.13}), then via (\ref{x90}) the bound on $[\frac{\partial v}{\partial a_0}]$ by Lemma \ref{LA1},
then the bound on $\|\partial_1^2\frac{\partial v}{\partial a_0}\|_{\alpha-2}$ by (\ref{1.13}), then via (\ref{x90})
finally the bound on $[\frac{\partial^2 v}{\partial a_0^2}]$ by Lemma \ref{LA1}.
The argument for (\ref{wj77}) is identical, just with $(f,v)$ replaced by $(\del f,\del v)$, cf. (\ref{wj65}),
and starting from (\ref{wj78}) instead of (\ref{wj58}) and thus with $N_0$ replaced by $\del N_0$.

\medskip

{\sc Step} \arabic{C1S}.\label{C1SL1}\refstepcounter{C1S} Application of Lemma \ref{L1}. We claim that
\begin{align}
\|[a,(\cdot)]\diam\partial_1^2u\|_{2\alpha-2} &\lesssim[a]M+NN_0,\label{wj85}\\
%\sup_{T\le 1}(T^\frac{1}{4})^{2-2\alpha}
\|[\del a,(\cdot)]\diam\partial_1^2u\|_{2\alpha-2}&\lesssim[\del a]M+\del NN_0,\label{wj86}\\
\|[a,(\cdot)]\diam\partial_1^2\del u\|_{2\alpha-2}&\lesssim[a]\del M
+N(N_0\del N+\del N_0).\label{wj87}
\end{align}
Here comes the argument: Estimate (\ref{wj85}) follows from Lemma \ref{L1}
with $b$ replaced by $a$, $I=1$ and $v_{i=1}=v$, so that the hypothesis (\ref{wi5})
is satisfied by (\ref{wj75}) in Step \ref{C1S0} with $N_0$ playing the role of $N_{i=1}$.
Hypothesis (\ref{1.5}) is satisfied by our assumption (\ref{wj57}) with $N$ playing the role of $N_0$.
In view of (\ref{wj59}), the outcome (\ref{x44}) of Lemma \ref{L1} turns into (\ref{wj85}).

\medskip

Estimate (\ref{wj86}) follows from applying Lemma \ref{L1}
with $b$ replaced by $\del a$, still $I=1$, $v_{i=1}=v$, and $N_0$ playing the role of $N_{i=1}$.
Hypothesis (\ref{1.5}) is satisfied by our assumption (\ref{wj92}) with $\del N$ playing the role of $N_0$.
In view of (\ref{wj59}), the outcome (\ref{x44}) of Lemma \ref{L1} turns into (\ref{wj86}).

\medskip
Finally, estimate (\ref{wj87}) follows from applying Lemma \ref{L1}
with $b$ again replaced by $a$, but this time $I=3$ and $(v_1,v_2,v_3)=(v,\frac{\partial v}{\partial a_0},\del v)$.
We learn from Step \ref{C1S0} that hypothesis (\ref{wi5}) is satisfied with $(N_1,N_2,N_3)=(N_0,N_0,\del N_0)$.
We now turn to the hypothesis (\ref{1.5}): For $i=1$ it is contained in our assumption (\ref{wj57}) with $N$
playing the role of $N_0$.
In preparation of checking hypothesis (\ref{1.5}) for $i=2$ we note that our assumption (\ref{wj57})
implies in particular that the family of distributions
$\{a\diam\partial_1^2 v(\cdot,a_0)\}_{a_0}$ is continuously differentiable in $a_0$.
This allows us to {\it define} the family of distributions
$\{a\diam\partial_1^2\frac{\partial v}{\partial a_0}(\cdot,a_0)\}_{a_0}$ via
\begin{align*}
a\diam\partial_1^2\frac{\partial v}{\partial a_0}:=
\frac{\partial}{\partial a_0}a\diam\partial_1^2v,
\end{align*}
which extends to the commutator:
\begin{align}\label{wj72}
[a,(\cdot)_T]\diam\partial_1^2\frac{\partial v}{\partial a_0}=
\frac{\partial}{\partial a_0}[a,(\cdot)_T]\diam\partial_1^2v.
\end{align}
Hence the hypothesis (\ref{1.5}) for $i=2$ is also satisfied by (\ref{wj57}) (here we use it up to
$\frac{\partial^2}{\partial a_0^2}$). Hypothesis (\ref{1.5}) for $i=3$ is identical to our assumption (\ref{wj69}).
We apply Lemma \ref{L1} with $\del u$ playing the role of $u$;
 the triple $(\del\sigma,\sigma\del a,\sigma)$ then plays the role of $(\sigma_1,\sigma_2,\sigma_3)$
and $\del M$ that of $M$. The outcome (\ref{x44}) of Lemma \ref{L1} assumes the form
\begin{align}\label{wj93}
&\|[a,(\cdot)]\diam\partial_1^2\del u\|_{2\alpha-2}\nonumber\\
&\qquad \lesssim[a]\del M
+N\big(N_0([\del\sigma]+\|\del\sigma\|+[\sigma\del a]+\|\sigma\del a\|)
+\del N_0([\sigma]+\|\sigma\|)\big).
\end{align}
We note that by (\ref{wj59}) and (\ref{wk60}) we have
\begin{align*}
&N_0([\del\sigma]+\|\del\sigma\|
+[\sigma\del a]+\|\sigma\del a\|)+\del N_0([\sigma]+\|\sigma\|)\nonumber\\
&\qquad \qquad \lesssim
N_0([\del\sigma]+\|\del\sigma\|+[\del a]+\|\del a\|)+\del N_0\lesssim
N_0\del N+\del N_0,\nonumber
\end{align*}
so that (\ref{wj93}) yields (\ref{wj87}).

\medskip

{\sc Step} \arabic{C1S}.\label{C1SComm}\refstepcounter{C1S} Commutator estimates. We claim
\begin{align}
\sup_{T\le 1}&(T^\frac{1}{4})^{2-2\alpha}\|\partial_2u_T-P(a\partial_1^2u_T+\sigma f_T)\|\lesssim[a]M+NN_0,\
\label{wj90}\\
\sup_{T\le 1}&(T^\frac{1}{4})^{2-2\alpha}\|\partial_2\del u_T
-P(a\partial_1^2\del u_T+\sigma\del aE\partial_1^2v_T+\sigma\del f_T+\del\sigma f_T)\|\nonumber\\
&\lesssim[a]\del M+([\del a]+\|\del a\|)M+N(N_0\del N+\del N_0)+\del N N_0.\label{wj91}
\end{align}
Indeed, we apply $(\cdot)_T$ to (\ref{wj55}) and rearrange terms:
\begin{align}\label{wj86bis}
\partial_2u_T-P(a\partial_1^2u_T+\sigma f_T)
&=-P([a,(\cdot)_T]\diam\partial_1^2u+[\sigma,(\cdot)_T]\diam f).
\end{align}
Similarly, we apply $(\cdot)_T$ to (\ref{wj66}) and rearrange terms:
\begin{align}
\partial_2\del u_T&-P(a\partial_1^2\del u_T+\sigma\del aE\partial_1^2v_T+\sigma\del f_T+\del\sigma f_T)\nonumber\\
&=-P\Big(-\del a(\partial_1^2u_T-\sigma E\partial_1^2v_T)
+[a,(\cdot)_T]\diam\partial_1^2\del u+[\del a,(\cdot)_T]\diam\partial_1^2u\nonumber\\
&\qquad +[\sigma,(\cdot)_T]\diam\del f+[\del\sigma,(\cdot)_T]\diam f\Big).\label{wj88}
\end{align}
By assumption (\ref{wj56}) and by (\ref{wj85}) in Step \ref{C1SL1} we obtain estimate (\ref{wj90}) from identity (\ref{wj86bis}).
By assumptions (\ref{wj67}) and (\ref{wj68}) and by (\ref{wj86}) and (\ref{wj87}) from Step \ref{C1SL1}
and from writing
\begin{align*}
\lefteqn{(\partial_1^2u_T-\sigma E\partial_1^2v_T)(x)
=\int dy\partial_1^2\psi_T(x-y)}\nonumber\\
&\times\big((u(y)-u(x))-\sigma(x)(v(y, a(x))-v(x,a(x)))-\nu(x)(y-x)_1\big),
\end{align*}
which entails with help of (\ref{1.13}) and (\ref{1.3})
\begin{align*}
\sup_{T\le 1}(T^\frac{1}{4})^{2-2\alpha}\|\del a(\partial_1^2u_T-\sigma E\partial_1^2v_T)\|
\lesssim\|\del a\|M,
\end{align*}
we obtain (\ref{wj91}) from formula (\ref{wj88}).

\medskip

{\sc Step} \arabic{C1S}.\refstepcounter{C1S} Application of Lemma \ref{L3} and conclusion.
We first apply Lemma \ref{L3} with $I=1$ and $f$ playing the role of $f_{i=1}$ (which does not depend on $a_0$).
The hypothesis (\ref{3.31}) is ensured by our assumption (\ref{wj58}) with $N_0$
playing the role of $N_{i=1}$. The hypothesis (\ref{y39}) is settled through (\ref{wj90}) in Step \ref{C1SComm}
with $N^2$ given by $[a]M+NN_0$.
Hence the two outputs (\ref{3.15}) and (\ref{wk01}) of Lemma \ref{L3} take the form of
\begin{align}
M&\lesssim [a]M+NN_0+N_0([\sigma]+\|\sigma\|[a]),\label{wj61}\\
[u]&\lesssim M+N_0\|\sigma\|.\label{wk04}
\end{align}
The smallness of $[a]$ and the boundedness of $\|\sigma\|$, cf. (\ref{wj59}), imply 
that (\ref{wj61}) simplifies to $M$ $\lesssim NN_0$ $+N_0([\sigma]+[a])$, which by (\ref{wk59})
means (\ref{wj63}).
Inserting (\ref{wj63}) into (\ref{wk04}) and using once more $\|\sigma\|\le 1$ yields (\ref{wk02}).

\medskip

We now apply Lemma \ref{L3} with $I=3$ and $(f,\partial_1^2v,\del f)$ playing the role of $(f_1,f_2,f_3)$;
by assumptions (\ref{wj58}), (\ref{wj78}) and by (\ref{wj75}), this triplet satisfies (\ref{3.31})
with $(N_1,N_2,N_3)=(N_0,N_0,\del N_0)$.
In view of (\ref{x90}) in Step \ref{C1S0}, and of
assumption (\ref{wj65}), the triplet $(v,\frac{\partial v}{\partial a_0},\del v)$
plays the role of $(v_1,v_2,v_3)$ in the sense of (\ref{v53}).
%; by (\ref{wj75})  \& (\ref{wj77}) in Step \ref{C1S0}, it also
%satisfies the estimates (\ref{3.31}) with $(N_1,N_2,N_3)=(N_0,N_0,\del N_0)$.
We apply Lemma \ref{L3} to $\del u$ playing the role of $u$,
$(\del\sigma,\sigma\del a,\sigma)$ playing the role of $(\sigma_1,\sigma_2,\sigma_3)$,
and $\del M$ playing the role of $M$. The hypothesis (\ref{y39}) is settled through Step \ref{C1SComm}
with $N^2$ estimated by the right hand side of (\ref{wj91}).
Hence the two outputs (\ref{3.15}) and (\ref{wk01}) of Lemma \ref{L3} take the form
\begin{align*}
\del M&\lesssim\mbox{expression on right hand side of (\ref{wj91})}+[a]\del M\nonumber\\
&+N_0([\del\sigma]+\|\del\sigma\|[a]+[\sigma\del a]+\|\sigma\del a\|[a])
+\del N_0([\sigma]+\|\sigma\|[a]),\\
[\del u]&\lesssim \del M+N_0(\|\del\sigma\|+\|\sigma\del a\|)
+\del N_0\|\sigma\|.
\end{align*}
Making use of the constraints (\ref{wj59}) on $\sigma$ and $a$, 
in particular to absorb $[a]\del M$ into the lhs, this simplifies to
\begin{align*}
\del M&\lesssim([\del a]+\|\del a\|)M+N(N_0\del N+\del N_0)\nonumber + \del N N_0\\
&+N_0([\del\sigma]+\|\del\sigma\|+[\del a]+\|\del a\|)
+\del N_0([\sigma]+[a]),\\
[\del u]&\lesssim \del M+N_0(\|\del\sigma\|+\|\del a\|)
+\del N_0.
\end{align*}
Inserting (\ref{wk59}) and  (\ref{wk60}), this reduces to
\begin{align}
\del M&\lesssim M\del N+N(N_0\del N+\del N_0)+N_0\del N, \label{hh2}\\
[\del u]&\lesssim \del M+N_0\del N \label{hh3}
+\del N_0.
\end{align}
Making use of the estimate (\ref{wj63}) on $M$ we just established, \eqref{hh2} implies
\begin{align*}
\del M\lesssim N(N_0\del N+\del N_0)+N_0\del N.%\quad\mbox{and}\quad
%[\del u]\lesssim \del M+N_0\del N+\del N_0.
\end{align*}
Clearly, this estimate implies the desired (\ref{wj80}). Plugging (\ref{wj80})
into \eqref{hh3} yields the desired (\ref{wk03}).

\bigskip

%%%%%%%%%%%%%%%%%%%%%%%%%%%%%%%%%%%%%%%%%%%%%%%%%%%%%%%%%%%%%%%%%%%%%%%%%%%%%%%%%%%%%%%%%%%%%%%%
%%%%%%%%%%%%%%%%%%%%%%%%%%%%%%%%%%%%%%%%%%%%%%%%%%%%%%%%%%%%%%%%%%%%%%%%%%%%%%%%%%%%%%%%%%%%%%%%
%%%%%%%%%%%%%%%%%%%%%%%%%%%%%%%%%%%%%%%%%%%%%%%%%%%%%%%%%%%%%%%%%%%%%%%%%%%%%%%%%%%%%%%%%%%%%%%%
%%%%%%%%%%%%%%%%%%%%%%%%%%%%%%%%%%%%%%%%%%%%%%%%%%%%%%%%%%%%%%%%%%%%%%%%%%%%%%%%%%%%%%%%%%%%%%%%
%%%%%%%%%%%%%%%%%%%%%%%%%%%%%%%%%%%%%%%%%%%%%%%%%%%%%%%%%%%%%%%%%%%%%%%%%%%%%%%%%%%%%%%%%%%%%%%%
%%%%%%%%%%%%%%%%%%%%%%%%%%%%%%%%%%%%%%%%%%%%%%%%%%%%%%%%%%%%%%%%%%%%%%%%%%%%%%%%%%%%%%%%%%%%%%%%

\subsection{Proof of Lemma \ref{L3}}
%\addcontentsline{toc}{section}{\protect\numberline{}Lemma \ref{L3}}

All functions are periodic if not stated otherwise.

\newcounter{proofstep} % proofstep = 0
\refstepcounter{proofstep} % increases value by 1

{\sc Step} \arabic{proofstep}.\label{L3step0}
\refstepcounter{proofstep} Estimate of $v_i$ and $\frac{\partial v_i}{\partial a_0}$. We claim
\begin{align}\label{v03}
[v_i]_{\alpha,1}\lesssim N_i,
\end{align}
where we recall the abbreviation (\ref{conv1}).
This follows immediately from assumption (\ref{3.31}) on $f_i$ and the definition (\ref{v53}) of $v_i$
via Lemma \ref{LA1} and the argument of Step~\ref{C1S0} of Corollary~\ref{C1}.

\medskip

{\sc Step} \arabic{proofstep}.\label{L3step2}
\refstepcounter{proofstep} Freezing-in the coefficients.
We claim that we have for all points $x_0$
\begin{align}\label{w92}
(\partial_2-a(x_0)\partial_1^2)\big(u_T
&-\sigma_i(x_0)v_{iT}(\cdot,a(x_0))\big) =Pg^T_{x_0},
\end{align}
where the function $g^T_{x_0}$ is estimated as follows
\begin{align}\label{w98}
|g^T_{x_0}(x)|\lesssim \tilde N^2\big((T^\frac{1}{4})^{2\alpha-2}
+(T^\frac{1}{4})^{\alpha-2}d^\alpha(x,x_0)\big)\quad\mbox{for}\;T\le 1
\end{align}
with the abbreviation
\begin{align}\label{v02}
\tilde N^2:=N^2+[a]_\alpha[u]_\alpha+N_i([\sigma_i]_\alpha+\|\sigma_i\|[a]_\alpha).
\end{align}
Indeed, making use of $P^2 = P$ we write
\begin{align}\label{w93}
(\partial_2-a(x_0)\partial_1^2)u_T=P(\sigma_i(x_0)f_{iT}(\cdot,a(x_0))+g^T_{x_0})
\end{align}
with $g^T_{x_0}$ defined through
\begin{align}\label{w96}
\lefteqn{g^T_{x_0}:=\partial_2u_T-P(a\partial_1^2u_T+\sigma_i Ef_{iT})+(a-a(x_0))\partial_1^2u_T}\nonumber\\
&+(\sigma_i-\sigma_i(x_0))Ef_{iT}+\sigma_i(x_0)(Ef_{iT}-f_{iT}(\cdot,a(x_0))).
\end{align}
By definition \eqref{v53} of $v_i(\cdot,a_0)$, to which we apply $(\cdot)_T$, which we evaluate for $a_0=a(x_0)$, and which
we contract with $\sigma_i(x_0)$ we obtain
\begin{align}\label{w94}
(\partial_2-a(x_0)\partial_1^2)\sigma_i(x_0)v_{iT}(\cdot,a(x_0))=P \sigma_i(x_0)f_{iT}(\cdot,a(x_0)).
\end{align}
From the combination of (\ref{w93}) and (\ref{w94}) we obtain (\ref{w92}),
so that it remains to estimate $g^T_{x_0}$. Making use of the assumption (\ref{y39}) we obtain from (\ref{w96})
\begin{align*}
|g^T_{x_0}(x)|&\le N^2(T^\frac{1}{4})^{2\alpha-2}+d^\alpha(x,x_0)\big([a]_\alpha\|\partial_1^2u_T\|\nonumber\\
&+[\sigma_i]_\alpha\sup_{a_0}\|f_{iT}\|
+\|\sigma_i\|[a]_\alpha\sup_{a_0}\|(\frac{\partial f_{i}}{\partial a_0})_T\|\big),
\end{align*}
so that by (\ref{1.13}) and by assumption (\ref{3.31})
\begin{align*}
|g^T_{x_0}(x)| &\lesssim N^2(T^\frac{1}{4})^{2\alpha-2}+(T^\frac{1}{4})^{\alpha-2}d^\alpha(x,x_0)
\big([a]_\alpha[u]_\alpha+N_i([\sigma_i]_\alpha+\|\sigma_i\|[a]_\alpha)\big)
\end{align*}
which can be consolidated into the estimate (\ref{w98}).

\medskip

{\sc Step} \arabic{proofstep}.\label{L3step3}
\refstepcounter{proofstep} PDE estimate.
Under the outcome of Step \ref{L3step2}, 
we have for all points $x_0$ and radii $R\ll L$
\begin{align}\label{w93bis}
\lefteqn{\frac{1}{R^{2\alpha}}\inf_\ell\|u_T-\sigma_i(x_0)v_{iT}(\cdot,a(x_0))-\ell\|_{B_R(x_0)}}\nonumber\\
&\lesssim (\frac{R}{L})^{2(1-\alpha)}\frac{1}{L^{2\alpha}}\inf_\ell\|u_T-\sigma_i(x_0)v_{iT}(\cdot,a(x_0))-\ell\|_{B_L(x_0)}\nonumber\\
&+\tilde N^2\big(\frac{L^2}{R^{2\alpha}(T^\frac{1}{4})^{2-2\alpha}}
+\frac{L^{2+\alpha}}{R^{2\alpha}(T^\frac{1}{4})^{2-\alpha}}\big),%\nonumber\\
\end{align}
where $\ell$ runs over all functions spanned by $1$ and $x_1$
and  $\|\cdot\|_{B_R(x_0)}$ denotes the supremum norm restricted to the ball $B_R(x_0)$
in the intrinsic metric (\ref{1.11}) with center $x_0$ and radius $R$. This step mimics the
heart of the kernel-free approach of  Safonov to the classical Schauder theory,
see \cite[Theorem 8.6.1]{Krylov}.
Here comes the argument: Wlog we restrict to $x_0=0$ and write $B_R=B_R(0)$ and $\|\cdot\|_R:=\|\cdot\|_{B_R}$.
Let $w_>$ be the (non-periodic) solution of
\begin{align*}
(\partial_2-a(0)\partial_1^2)w_>=I(B_L)g_{0}^T,
\end{align*}
where $I(B_L)$ denotes the indicator function of the set $B_L$. 
Hence in view of (\ref{w92}), where we write $Pg_0^T = g_0^T +c $ with $c = -\int_{[0,1)^2} g_0^T$, the function
\begin{align}\label{v06}
w_<:=u_T-\sigma_i(0)v_{iT}(\cdot,a(0))-w_>
\end{align}
satisfies
\begin{align}\label{v01}
(\partial_2-a(0)\partial_1^2)w_<=c\quad\mbox{in}\;B_L.
\end{align}
By standard estimates for the heat equation we have
\begin{align}
\|w_>\|&\lesssim L^2\|g_0^T\|_L,\label{3.2}\\
\|\{ \partial_1^2,\partial_2\}w_<\|_{\frac{L}{2}}&\lesssim L^{-2}\|w_<-\ell_L\|_L\label{3.1}
\end{align}
for any function $\ell_L\in{\rm span}\{1,x_1\}$. The interior estimate (\ref{3.1}) is slightly
non-standard because of the non-vanishing right hand side $c$ but can be easily reduced
to the case of $c=0$: First of all, replacing $w$ by $w-\ell_L$ in (\ref{v01}) and (\ref{3.1})
we may reduce to the case of $\ell_L=0$. Testing (\ref{v01}) with a cut-off function for $B_L$
that is smooth on scale $L$ we learn that $|c|\lesssim L^{-2}\|w_<\|_L$.
We then may replace $w$ by $w+cx_2$ which reduces the further estimate to the
standard case of $c=0$. We refer to \cite[Theorem 8.4.4]{Krylov} for an elementary argument for (\ref{3.1})
in case of $c=0$ only relying on the maximum principle via Bernstein's argument. We refer
to \cite[Exercise 8.4.8]{Krylov} for the statement (\ref{3.2}) via the representation through
the heat kernel.
Since by construction, cf. (\ref{v06}), we have $u_T$ $-\sigma_i(0)v_{iT}(\cdot,a(0))$ 
$= w_<+w_>$
we obtain by the triangle inequality for a suitably chosen $\ell_R\in{\rm span}\{1,x_1\}$
\begin{align*}
\|u_T&-\sigma_i(0)v_{iT}(\cdot,a(0))-\ell_R\|_R
\nonumber\\
&\le\|w_<-\ell_R\|_R+\|w_>\|_R
\lesssim R^2\|\{\partial_1^2,\partial_2\}w_<\|_{R}+\|w_>\|_R.
\end{align*}
Inserting (\ref{3.1}) for $R\ll L$, and by another application of the triangle inequality this yields
\begin{align*}
\|u_T&-\sigma_i(0)v_{iT}(\cdot,a(0))-\ell_R\|_R\nonumber\\
&\lesssim L^{-2}R^2\|w_<-\ell_L\|_{L}+\|w_>\|_R\nonumber\\
&\le L^{-2}R^2\|u_T-\sigma_i(0)v_{iT}(\cdot,a(0))-\ell_L\|_{L}+2\|w_>\|.
\end{align*}
Inserting (\ref{3.2}) \& (\ref{w98}) this gives
\begin{align}\label{3.10}
\lefteqn{\inf_{\ell}\|u_T-\sigma_i(0)v_{iT}(\cdot,a(0))-\ell\|_R}\nonumber\\
&\lesssim L^{-2}R^2\inf_{\ell}\|u_T-\sigma_i(0)v_{iT}(\cdot,a(0))-\ell\|_{L}
+\tilde N^2 L^2((T^\frac{1}{4})^{2\alpha-2}+L^\alpha(T^\frac{1}{4})^{\alpha-2}),
\end{align}
where we recall that $\ell$ runs over ${\rm span}\{1,x_1\}$. Dividing by $R^{2\alpha}$ gives (\ref{w93bis}).

\medskip

{\sc Step} \arabic{proofstep}.\label{L3step4}
\refstepcounter{proofstep} Equivalence of norms. We claim that the modelling constant
$M$ of $u$ is estimated by the expression appearing in Step \ref{L3step3}:
\begin{align}\label{3.8}
M\lesssim M',
\end{align}
where we have set for abbreviation
\begin{align}\label{3.6}
M':=\sup_{x_0}&\sup_{R\le 1}R^{-2\alpha}\inf_{\ell}\|u-\sigma_i(x_0)v_i(\cdot,a(x_0))-\ell\|_{B_R(x_0)}
\end{align}
and where the maximal radius $1$ is chosen such that a ball of that covers a periodic cell. 
In fact, also the reverse estimate
holds, highlighting once more that the modulation function $\nu$ in the definition of modelledness (Definition~\ref{D1}) plays
a small role compared to $\sigma_i$. The equivalence of (\ref{3.8}) and (\ref{3.6})
on the level of standard H\"older spaces is the starting point for the approach to Schauder theory
by  Safonov, see \cite[Theorem 8.5.2]{Krylov}.
We first argue that the $\ell$ in (\ref{3.6}) may be chosen to be independent of $R$, that is,
\begin{align}\label{3.7}
\sup_{x_0}\inf_{\ell}\sup_{R\le 1}R^{-2\alpha}\|u-\sigma_i(x_0)v_i(\cdot,a(x_0))-\ell\|_{B_R(x_0)}\lesssim M'.
\end{align}
Indeed, fix $x_0$, say $x_0=0$, and let $\ell_R=\nu_Rx_1+c_R$ be (near) optimal in (\ref{3.6}), then we have by
definition of $M'$ and by the triangle inequality $R^{-2\alpha}\|\ell_{2R}-\ell_{R}\|_R\lesssim M'$.
This implies $R^{1-2\alpha}|\nu_{2R}-\nu_{R}|+R^{-2\alpha}|c_{2R}-c_R|\lesssim M'$,
which thanks to $\alpha>\frac{1}{2}$ yields by telescoping $R^{1-2\alpha}|\nu_{R}-\nu_{R'}|+R^{-2\alpha}|c_{R}-c_{R'}|\lesssim M'$
for all $R'\le R$ and thus the existence of $\nu,c\in\mathbb{R}$ such that $R^{1-2\alpha}|\nu_{R}-\nu|+R^{-2\alpha}|c_{R}-c|\lesssim M'$, so that $\ell:=\nu x_1+c$ satisfies
\begin{align}\label{1.21}
R^{-2\alpha}\|\ell_R-\ell\|_R\lesssim M'.
\end{align}
Hence we may pass from (\ref{3.6})
to (\ref{3.7}) by the triangle inequality.

\medskip

It is clear from (\ref{3.7}) that necessarily for any $x_0$, say $x_0=0$,
the optimal $\ell$ must be of the form $\ell(x)$
$=u(0)$ $-\sigma_i(0)v_i(0,a(0))$ $-\nu(0)x_1$.
This establishes the
main part of (\ref{3.8}), namely the modelledness (\ref{1.3}) for any ``base'' point $x$
and any $y$ of distance at most $1$.
Since $B_1(x)$ covers a periodic cell, by periodicity of $y\mapsto (u(y)-u(x))$ $-\sigma_i(x)(v_i(y,a(x))-v_i(x,a(x)))$
we extract $|\nu(x)|\lesssim M'$. Since $\alpha\ge\frac{1}{2}$,
this implies that $|\nu(x)(x-y)_1|\lesssim M'd^{2\alpha}(x,y)$ for all $y\not\in B_1(x)$. Hence
once again by periodicity of
$y\mapsto (u(y)-u(x))$ $-\sigma_i(x)(v_i(y,a(x))-v_i(x,a(x)))$,
(\ref{1.3}) holds also for $y\not\in B_1(x)$.

\medskip

{\sc Step} \arabic{proofstep}.\label{L3step5}
\refstepcounter{proofstep} Modelledness implies approximation property.
We claim that for any mollification parameter $0<T\le 1$, radius $L$, and point $x_0$ we have
\begin{align}\label{3.11}
\frac{1}{(T^\frac{1}{4})^{2\alpha}}&\|(u_T-u)-\sigma_i(x_0)(v_{iT}-v_i)(\cdot,a(x_0))\|_{B_L(x_0)}
\lesssim M+\tilde N^2 (\frac{L}{T^\frac{1}{4}})^\alpha.
\end{align}
%
%where we recall the definition (\ref{v02}) of $\tilde N$.
Wlog we consider $x_0=0$ %, write $v_i(y,x)=v_i(y,a(x))$, 
and recall that the first moment of $\psi_T$ vanishes, so that
\begin{align*}
\lefteqn{(u_T-u)(x)-\sigma_i(0)(v_{iT}-v_i)(x,a(0))}\nonumber\\
&=\int dy\psi_T(x-y)\big((u(y)-u(x))-\sigma_i(0)(v_i(y,a(0))-v_i(x,a(0)))\nonumber\\
&-\nu(x)(y-x)_1\big).
\end{align*}
We split the right hand side into three terms:
\begin{align*}
\lefteqn{(u_T-u)(x)-\sigma_i(0)(v_{iT}-v_i)(x,a(0))}\nonumber\\
&=\int dy\psi_T(x-y)\big((u(x)-u(y))-\sigma_i(x)(v_i(y,a(x))-v_i(x,a(x)))\nonumber\\
&\qquad-\nu(x)(x-y)_1\big)\nonumber\\
&+\int dy\psi_T(x-y)(\sigma_i(x)-\sigma_i(0))(v_i(y,a(0))-v_i(x,a(0)))\nonumber\\
&+\int dy\psi_T(x-y)\sigma_i(x)\big((v_i(y,a(x))-v_i(y,a(0)))-(v_i(x,a(x))-v_i(x,a(0))\big).\nonumber
\end{align*}
For the first right-hand-side term we appeal to the modelledness assumption (\ref{1.3}), which implies that the
integrand is estimated by $|\psi_T(x-y)|$ $M$ $d^{2\alpha}(x,y)$. Hence by (\ref{1.13}) the integral is estimated
by $M$ $(T^\frac{1}{4})^{2\alpha}$.
The integrand of the second rhs term is estimated by $|\psi_T(x-y)|$
$[\sigma_i]_{\alpha}$ $d^\alpha(x,0)$
$[v_i(\cdot,0)]_\alpha$ $d^\alpha(x,y)$ so that by
(\ref{1.13}) and (\ref{v03}) the integral is controlled by
$\lesssim [\sigma_i]_{\alpha}$ $d^\alpha(x,0)$
$N_i$ $(T^\frac{1}{4})^\alpha$; since $x\in B_L(0)$ it
is controlled by $\lesssim [\sigma_i]_{\alpha}$ $L^\alpha$ $N_i$
$(T^\frac{1}{4})^\alpha$. Using the identity (and dropping the index $i$)
\begin{align*}
\lefteqn{(v(y,a(x))-v(y,a(0)))-(v(x,a(x))-v(x,a(0)))=(a(x)-a(0))}\nonumber\\
&\times\int_0^1ds\big(\frac{\partial v}{\partial a_0}(y,sa(x)+(1-s)a(0))
-\frac{\partial v}{\partial a_0}(x,sa(x)+(1-s)a(0)))\big),
\end{align*}
we see that the integrand of the third right-hand-side term is estimated by $|\psi_T(x-y)|\|\sigma_i\|$
$d^\alpha (x,y)$ $[a]_{\alpha}$
$\sup_{a_0}[\frac{\partial v_i}{\partial a_0}(\cdot,a_0)]_\alpha$ $d^\alpha(x,0)$;
hence in view of (\ref{v03}) the third term itself is estimated by $\|\sigma_i\|$ $N_i$ $(T^\frac{1}{4})^\alpha$
$[a]_{\alpha}$ $L^\alpha$.
Collecting these estimates we obtain for $x\in B_L(0)$
\begin{align*}
|(u_T-u)(x)-\sigma_i(0)(v_{iT}-v_i)(x,0)|
\lesssim
M(T^\frac{1}{4})^{2\alpha}
+N_i([\sigma_i]_{\alpha}+\|\sigma_i\|[a]_{\alpha}) 
L^\alpha (T^\frac{1}{4})^\alpha.
\end{align*}
In view of the definition (\ref{v02}) of $\tilde N^2$, this yields (\ref{3.11}).

\medskip

{\sc Step} \arabic{proofstep}.\label{L3step6}
\refstepcounter{proofstep} Estimate of $M$. We claim that
\begin{align}\label{v54}
M\lesssim \tilde N^2.
\end{align}
Indeed, we can now close the argument and to this purpose rewrite (\ref{w93bis})
from Step \ref{L3step3} with help of the triangle inequality as
\begin{align}\nonumber
\lefteqn{\frac{1}{R^{2\alpha}}\inf_{\ell}\|u-\sigma_i(x_0)v_i(\cdot,a(x_0))-\ell\|_{B_R(x_0)}}\nonumber\\
&\lesssim \Big(\frac{R}{L} \Big)^{2-2\alpha}\frac{1}{L^{2\alpha}}\inf_{\ell}\|u-\sigma_i(x_0)v_i(\cdot,a(x_0))-\ell\|_{B_L(x_0)}\nonumber\\
&+\tilde N^2 \Big(\frac{L^2}{R^{2\alpha}(T^\frac{1}{4})^{2-2\alpha}}
+\frac{L^{2+\alpha}}{R^{2\alpha}(T^\frac{1}{4})^{2-\alpha}}\Big)\nonumber\\
&+\Big(\frac{T^\frac{1}{4}}{R} \Big)^{2\alpha}\frac{1}{(T^\frac{1}{4})^{2\alpha}}
\|(u_T-u)-\sigma_i(x_0)(v_{iT}-v_i)(\cdot,a(x_0))\|_{B_L(x_0)}.\nonumber
\end{align}
We now insert (\ref{3.11}) from Step \ref{L3step5} to obtain
\begin{align}\nonumber
\frac{1}{R^{2\alpha}}\inf_{\ell}&\|u-\sigma_i(x_0)v_i(\cdot,a(x_0))-\ell\|_{B_R(x_0)}\nonumber\\
&\lesssim \Big(\frac{R}{L} \Big)^{2-2\alpha}M+\tilde N^2 \Big(\frac{L^2}{R^{2\alpha}(T^\frac{1}{4})^{2-2\alpha}}
+\frac{L^{2+\alpha}}{R^{2\alpha}(T^\frac{1}{4})^{2-\alpha}}\Big)\nonumber\\
&+\Big(\frac{T^\frac{1}{4}}{R}\Big)^{2\alpha}M+\tilde N^2\frac{L^\alpha (T^\frac{1}{4})^\alpha}{R^{2\alpha}}.\label{v04}
\end{align}
Here we have used that
\begin{align*}
\sup_{x_0}\sup_{L}\frac{1}{L^{2\alpha}}\inf_{\ell}\|u-\sigma_i(x_0)v_i(\cdot,a(x_0))-\ell\|_{B_L(x_0)}\lesssim M
\end{align*}
by the definition
of the modelling constant $M$ with $\ell_{x_0}(x)$ $=u(x_0)$ $-\sigma_i(x_0)v_i(x_0,a(x_0))$ $-\nu(x_0)(x-x_0)_1$.
Relating the length scales $T^\frac{1}{4}$ and $L$ to the given $R\le 1$ in (\ref{v04}) via
$T^\frac{1}{4}=\epsilon R$ (so that in particular as required $T\le 1$ since we think of $\epsilon\ll 1$)
and $L=\epsilon^{-1}R$, taking the supremum over $R\le 1$ and $x_0$ yields by definition (\ref{3.6}) of $M'$
\begin{align}\nonumber
M'
\lesssim (\epsilon^{2-2\alpha}+\epsilon^{2\alpha}) M
+\big(\epsilon^{2\alpha-4}+\epsilon^{-4}+1\big)\tilde N^2.\nonumber
\end{align}
By (\ref{3.8}) in Step \ref{L3step4}, this implies
\begin{align}\nonumber
M
\lesssim (\epsilon^{2-2\alpha}+\epsilon^{2\alpha}) M+\epsilon^{-4}\tilde N^2.\nonumber
\end{align}
Since $0<\alpha<1$, we may choose $\epsilon$ sufficiently small such that the first right-hand-side term may be absorbed into
the lhs yielding the desired estimate $M\lesssim\tilde N^2$ (note that $M<\infty$ is part of our assumption).

\medskip

{\sc Step} \arabic{proofstep}.
\refstepcounter{proofstep} Conclusion. Clearly, (\ref{3.15}) and (\ref{wk01}) immediately follow from the
combination of
\begin{align*}
M\lesssim N^2+[a]_\alpha[u]_\alpha+N_i([\sigma_i]_\alpha+\|\sigma_i\|[a]_\alpha),\quad [u]_\alpha\lesssim M+N_i\|\sigma_i\|.
\end{align*}
The first estimate is identical to (\ref{v54}) in Step \ref{L3step6} into which we plug the definition (\ref{v02})
of $\tilde N$. The second estimate is an application of Step \ref{L2aS-1} in the proof of Lemma \ref{L2a} with
$v(y,x):=\sigma_i(x)v_i(y,a_i(x))$, so that the hypothesis (\ref{wi75}) holds with
$N$ replaced by $\|\sigma_i\|N_i$, cf. (\ref{v03}) in Step \ref{L3step0}.

\bigskip

%%%%%%%%%%%%%%%%%%%%%%%%%%%%%%%%%%%%%%%%%%%%%%%%%%%%%%%%%%%%%%%%%%%%%%%%%%%%%%%%%%%%%%%%%%%%%%%%
%%%%%%%%%%%%%%%%%%%%%%%%%%%%%%%%%%%%%%%%%%%%%%%%%%%%%%%%%%%%%%%%%%%%%%%%%%%%%%%%%%%%%%%%%%%%%%%%
%%%%%%%%%%%%%%%%%%%%%%%%%%%%%%%%%%%%%%%%%%%%%%%%%%%%%%%%%%%%%%%%%%%%%%%%%%%%%%%%%%%%%%%%%%%%%%%%
%%%%%%%%%%%%%%%%%%%%%%%%%%%%%%%%%%%%%%%%%%%%%%%%%%%%%%%%%%%%%%%%%%%%%%%%%%%%%%%%%%%%%%%%%%%%%%%%
\subsection{Proof of Lemma \ref{L2a}}
%\addcontentsline{toc}{section}{\protect\numberline{}Lemma \ref{L2a}}
%{\sc Proof of Lemma \ref{L2a}}.

We write for abbreviation $[\cdot]:=[\cdot]_\alpha$ and $E:=\Ed$.
\newcounter{L2aS} % proofstep = 0
\refstepcounter{L2aS} % increases value by 1

{\sc Step} \arabic{L2aS}\label{L2aS0}\refstepcounter{L2aS}.
We claim
\begin{align}\label{wi74}
[\nu]_{2\alpha-1}\lesssim M+N.
\end{align}
Indeed, introducing $\ell_{x}(y):=\nu(x)y_1$ we see that (\ref{wi81}) can be rewritten as
\begin{align*}
|(u-v(\cdot,x)-\ell_{x})(y)-(u-v(\cdot,x)-\ell_{x})(x)|&\le Md^{2\alpha}(y,x),
\end{align*}
so that we obtain by the triangle inequality
\begin{align}\label{wi73}
|(u-v(\cdot,x)-\ell_{x})(y)-(u-v(\cdot,x)-\ell_{x})(y')|
&\le M(d^{2\alpha}(y,x)+d^{2\alpha}(y',x)).
\end{align}
In combination with (\ref{wi80}) this yields by the triangle inequality
\begin{align*}
|&(u-v(\cdot,x')-\ell_{x})(y)-(u-v(\cdot,x')-\ell_{x})(y')|\nonumber\\
&\qquad \le M(d^{2\alpha}(y,x)+d^{2\alpha}(y',x))+Nd^\alpha(x,x')d^\alpha(y,y').
\end{align*}
We now take the difference of this with (\ref{wi73}) with $x$ replaced by $x'$ to obtain, once more by the triangle inequality,
\begin{align*}
|&(\ell_{x}-\ell_{x'})(y)-(\ell_{x}-\ell_{x'})(y')|\nonumber\\
&\le M\big(d^{2\alpha}(y,x)+d^{2\alpha}(y',x)
+d^{2\alpha}(y,x')+d^{2\alpha}(y',x')\big)
+Nd^\alpha(x,x')d^\alpha(y,y').
\end{align*}
By definition of $\ell$ and with the choice of $y=x$ and $y'=x+(R,0)$, this assumes the form
\begin{align*}
|\nu(x)-\nu(x')|R\le M(R^{2\alpha}+d^{2\alpha}(x,x')+(R+d(x,x'))^{2\alpha})
+Nd^\alpha(x,x')R^\alpha.
\end{align*}
With the choice of $R=d(x,x')$ this turns into
\begin{align*}
|\nu(x)-\nu(x')|d(x,x')\lesssim (M+N)d^{2\alpha}(x,x'),
\end{align*}
which amounts to the desired (\ref{wi74}).

\medskip

{\sc Step} \arabic{L2aS}\label{L2aS-1}\refstepcounter{L2aS}.
Under our additional assumption (\ref{wi75}) we claim
\begin{align}\label{wi76}
[u]+\|\nu\|\lesssim M+N.
\end{align}
By the triangle inequality on (\ref{wi81}) we obtain for all pairs of points
$|\nu(x)(x-y)_1|$ $\le|u(x)-u(y)|$ $+[v(\cdot,x)] d^\alpha(y,x)$ $+Md^{2\alpha}(x,y)$. Choosing
$y=x+(1,0)$, appealing to the space-time periodicity of $u$, taking the supremum over $x$, and appealing
to (\ref{wi75}), this turns into the $\nu$-part of (\ref{wi76}):
\begin{align}\label{C0.2}
\|\nu\|\lesssim M+N.
\end{align}
We now consider pairs of points $(x,y)$ with $d(x,y)\le 1$. By the triangle inequality
from (\ref{wi81}) we get
\begin{align*}
\frac{1}{d^\alpha(x,y)}|u(x)-u(y)|\lesssim M+N+\|\nu\|.
\end{align*}
By space-time periodicity, this extends to all pairs so that
\begin{align*}
[u]\lesssim M+N+\|\nu\|.
\end{align*}
Inserting (\ref{C0.2}) into this yields the $u$-part of (\ref{wi76}).

\medskip

{\sc Step} \arabic{L2aS}\label{L2aS1}\refstepcounter{L2aS}.
Dyadic decomposition. For $\tau<T$ (with $T$ a dyadic multiple of $\tau$) we claim that
\begin{align}\label{w1}
%\lefteqn{
\big(uf_T-E [v,(\cdot)_T]\diam f &
-\nu[x_1,(\cdot)_T]f\big)
%}\nonumber\\
-\big(uf_\tau-E[v,(\cdot)_\tau]\diam f-\nu[x_1,(\cdot)_\tau]f\big)_{T-\tau}\nonumber\\
&=\sum_{\tau\le t<T}\Big(\big([u,(\cdot)_t]-E[v,(\cdot)_t]
-\nu[x_1,(\cdot)_t]\big)f_t
\nonumber\\
& \qquad\qquad -[\nu,(\cdot)_t][x_1,(\cdot)_t]f-[E,(\cdot)_t][v,(\cdot)_t]\diam f
\Big)_{T-2t},
\end{align}
where the sum runs over the dyadic ``times'' $t=\frac{T}{2},\frac{T}{4},\cdots,\tau$.
By telescoping based on the semi-group property (\ref{1.10}) this reduces to
\begin{align*}
 &\big(uf_{2t}-E[v,(\cdot)_{2t}]\diam f
-\nu[x_1,(\cdot)_{2t}]f\big)
-\big(uf_t-E[v,(\cdot)_t]\diam f
-\nu[x_1,(\cdot)_t]f\big)_{t}\nonumber\\
& \qquad =\big([u,(\cdot)_t]-E[v,(\cdot)_t]
-\nu[x_1,(\cdot)_t]\big)f_t
-[\nu,(\cdot)_t][x_1,(\cdot)_t]f-[E,(\cdot)_t][v,(\cdot)_t]\diam f,
\end{align*}
and splits into the three statements
\begin{align}\label{w2}
uf_{2t}-(uf_t)_t&=[u,(\cdot)_t]f_t,\\
\nu[x_1,(\cdot)_{2t}]f-(\nu[x_1,(\cdot)_{t}]f)_t&=\nu[x_1,(\cdot)_t]f_t+[\nu,(\cdot)_t][x_1,(\cdot)_t]f,
\nonumber\\
E[v,(\cdot)_{2t}]\diam f-(E[v,(\cdot)_t]\diam f)_t&=E[v,(\cdot)_t]f_t+[E,(\cdot)_t][v,(\cdot)_t]\diam f.\nonumber
\end{align}
Plugging in the definition of the commutator $[\nu,(\cdot)_t]$, the middle statement
reduces to
\begin{align}\label{w3}
[x_1,(\cdot)_{2t}]f-([x_1,(\cdot)_{t}]f)_t&=[x_1,(\cdot)_t]f_t.
%E[v,(\cdot)_{2t}]\diam f-(E[v,(\cdot)_t]\diam f)_t&=E[v,(\cdot)_t]f_t
%+[E,(\cdot)_t][v,(\cdot)_t]\diam f.\nonumber
\end{align}
By the definition of the commutator $[E,(\cdot)_t]$, the last statement reduces to
\begin{align}\label{wi16}
[v,(\cdot)_{2t}]\diam f-([v,(\cdot)_t]\diam f)_t&=[v,(\cdot)_t]f_t,
\end{align}
which by definition of $[v,(\cdot)_{T}]\diam f$ splits into
\begin{align}
vf_{2t}-(vf_t)_t=[v,(\cdot)_t]f_t&\quad\mbox{and}\quad
(v\diam f)_{2t}-((v\diam f)_t)_t=0.\label{w4}
\end{align}
Now identities (\ref{w2}), (\ref{w3}), and (\ref{w4}) follow immediately from the semi-group property.

\medskip

{\sc Step} \arabic{L2aS}\label{L2aS2}\refstepcounter{L2aS}.
For $\tau<T\le 1$ (with $T$ still a dyadic multiple of $\tau$) we claim the estimate
\begin{align}\label{w1bis}
\lefteqn{\|\big(uf_T-E[v,(\cdot)_T]\diam f
-\nu[x_1,(\cdot)_T]f\big)}\nonumber\\
&-\big(uf_\tau-E[v,(\cdot)_\tau]\diam f
-\nu[x_1,(\cdot)_\tau]f\big)_{T-\tau}\|\nonumber\\
&\lesssim (M+N)N_1(T^\frac{1}{4})^{3\alpha-2}.
\end{align}
Indeed, by the dyadic representation (\ref{w1}), the triangle inequality in $\|\cdot\|$
and the fact that $(\cdot)_{T-2t}$ is bounded in that norm, cf. (\ref{1.13}), it is enough
to show that the right-hand-side term of (\ref{w1}) under the parenthesis is estimated by
$(M+N)N_1$
$(t^\frac{1}{4})^{3\alpha-2}$ for all $t\le 1$; here we crucially use that by assumption $3\alpha-2>0$ for the 
convergence of the geometric series.
Using Step \ref{L2aS0} to control $[\nu]_{2\alpha-1}$ in (\ref{w5}) by $M+N$,
this estimate splits into
\begin{align}
\|\big([u,(\cdot)_t]-E[v,(\cdot)_t]
-\nu[x_1,(\cdot)_t]\big)f_t\|&\lesssim MN_1 (t^\frac{1}{4})^{3\alpha -2},\nonumber\\
\|[\nu,(\cdot)_t][x_1,(\cdot)_t]f\|&\lesssim [\nu]_{2\alpha-1}N(t^\frac{1}{4})^{3\alpha -2}, \label{w5}\\
\|[E,(\cdot)_t][v,(\cdot)_t]\diam f\|
&\lesssim NN_1(t^\frac{1}{4})^{3\alpha -2}.\label{w8}
\end{align}
Appealing to our assumptions (\ref{wi78}) \& (\ref{wi79}) and to Lemma \ref{LA2}, these three estimates reduce to
\begin{align}
\|\big([u,(\cdot)_t]&-E[v,(\cdot)_t]
-\nu[x_1,(\cdot)_t]\big)\tilde f\|\lesssim M \|\tilde f\|(t^\frac{1}{4})^{2\alpha},\nonumber\\
\|[\nu,(\cdot)_t]\tilde f\|&\lesssim [\nu]_{\beta}\|\tilde f\|(t^\frac{1}{4})^{\beta},\label{wi2}\\
\|[E,(\cdot)_t]\tilde v\|
&\lesssim \sup_{x,x'}\frac{1}{d^\alpha(x,x')}\|\tilde v(\cdot,x)-\tilde v(\cdot,x')\|(t^\frac{1}{4})^{\alpha},\label{w6}
\end{align}
where $\tilde f=\tilde f(y)$ plays the role of $f_t$ or $[x_1,(\cdot)_t]f$,
and $\tilde v=\tilde v(x,y)$ plays the role of $([v(\cdot,x),(\cdot)_t]\diam f)(y)$,
but now can be, like $\nu$, generic functions; similarly, $\beta$ plays the role of $2\alpha-1$ but could
be any exponent in $[0,1]$.
Using the definition of $E$, we may rewrite these estimates more explicitly as
\begin{align}
\Big|\int dy\psi_t(x-y)\Big(\big(u(x)-u(y)\big)-\big(v(x,x)&-v(y,x)\big)\nonumber\\
-\nu (x)(x-y)_1\Big)\tilde f(y)\Big| &\lesssim M \|\tilde f\|(t^\frac{1}{4})^{2\alpha},\nonumber\\
\Big|\int dy\psi_t(x-y)\big(\nu(x)-\nu(y)\big)\tilde f(y) \Big|
&\lesssim [\nu]_{\beta} \|\tilde f\|(t^\frac{1}{4})^{\beta},\nonumber\\
\Big|\int dy\psi_t(x-y)\big(\tilde v(y,x)-\tilde v(y,y)\big)\Big|
%\nonumber\\&
&\lesssim \sup_{x,x'}\frac{1}{d^\alpha(x,x')}\|\tilde v(\cdot,x)-\tilde v(\cdot,x')\|(t^\frac{1}{4})^{\alpha}.\nonumber
\end{align}
All three estimates rely on the moment bounds (\ref{1.13}),
the first estimate is then an immediate consequence of (\ref{wi81})
and the two last ones tautological.

\medskip

{\sc Step} \arabic{L2aS}\label{L2aS3}\refstepcounter{L2aS}. For
\begin{align*}
F^\tau:=uf_\tau-E[v,(\cdot)_\tau]\diam f-\nu[x_1,(\cdot)_\tau]f
\end{align*}
and under our addditional assumptions (\ref{wi75}) \& (\ref{wi85}) we claim the estimates
\begin{align}
\sup_{T\le 1}&(T^\frac{1}{4})^{2-2\alpha}
\|uf_T-F^\tau_{T-\tau}\|\lesssim
(M+N)N_1,\quad
\|F^\tau\|_{\alpha-2}%=\sup_{T\le 1}&(T^\frac{1}{4})^{2-\alpha}\|F^\tau_T\|
\lesssim(M+N+\|u\|)N_1.\label{v44}
\end{align}
Indeed, the first item in (\ref{v44}) follows from (\ref{w1bis}) in Step \ref{L2aS2} 
via the triangle inequality and
\begin{align*}
\|E[v,(\cdot)]\diam f\|_{2\alpha-2}\stackrel{(\ref{wi85})}{\le} N N_1,\quad
\|\nu[x_1,(\cdot)]f\|_{2\alpha-2}\lesssim(M+N)N_1,\nonumber
\end{align*}
the latter being a consequence of (\ref{wi76}) in Step \ref{L2aS-1}, (\ref{LA2.2})
in Lemma \ref{LA2}, and our assumption (\ref{wi78}); here, we make extensively use
of $T\le 1$. The second item in (\ref{v44}) in turn follows from (\ref{v44}) via
$\|F_T^\tau\|=\|(F_{T-\tau}^\tau)_\tau\|\lesssim\|F_{T-\tau}^\tau\|$ (cf. (\ref{1.10}) and (\ref{1.13}))
by the triangle inequality, (\ref{wi76}), and (\ref{wi78}), again making use of $T\le 1$.

\medskip

{\sc Step} \arabic{L2aS}\label{L2aS4}\refstepcounter{L2aS}. Conclusion:
By the second item in (\ref{v44}) in Step \ref{L2aS3}, the sequence $\{F^\tau\}_{\tau\downarrow0}$ 
is bounded wrt $\|\cdot\|_{\alpha-2}$. By standard weak compactness based on the equivalence of norms
from Step \ref{LA1S1} in the proof of Lemma \ref{LA1}, there exists a subsequence $\tau_n\downarrow0$
and a distribution we give the name of $u\diam f$ such that $F^{\tau_n}\rightharpoonup u\diam f$.
By standard lower semi-continuity, we may pass to the limit in (\ref{v44}) in Step \ref{L2aS3} to obtain
(\ref{wi83}). Likewise, we may pass to the limit
in (\ref{w1bis}) in Step \ref{L2aS2} to obtain (\ref{wi84}). 
Note that our additional assumptions (\ref{wi75}) \& (\ref{wi85}) were only qualitatively used
in deriving (\ref{wi84}) by ensuring the above boundedness of $\{F^\tau\}_{\tau\downarrow0}$.

\bigskip
%%%%%%%%%%%%%%%%%%%%%%%%%%%%%%%%%%%%%%%%%%%%%%%%%%%%%%%%%%%%%%%%%%%%%%%%%%%%%%%%%%%%%%%%%%%%%
%%%%%%%%%%%%%%%%%%%%%%%%%%%%%%%%%%%%%%%%%%%%%%%%%%%%%%%%%%%%%%%%%%%%%%%%%%%%%%%%%%%%%%%%%%%%%
%%%%%%%%%%%%%%%%%%%%%%%%%%%%%%%%%%%%%%%%%%%%%%%%%%%%%%%%%%%%%%%%%%%%%%%%%%%%%%%%%%%%%%%%%%%%%
%%%%%%%%%%%%%%%%%%%%%%%%%%%%%%%%%%%%%%%%%%%%%%%%%%%%%%%%%%%%%%%%%%%%%%%%%%%%%%%%%%%%%%%%%%%%%
%%%%%%%%%%%%%%%%%%%%%%%%%%%%%%%%%%%%%%%%%%%%%%%%%%%%%%%%%%%%%%%%%%%%%%%%%%%%%%%%%%%%%%%%%%%%%

\subsection{Proof of Lemma \ref{L1}}
%\addcontentsline{toc}{section}{\protect\numberline{}Lemma \ref{L1}}
%{\sc Proof of Lemma \ref{L1}}.

The proof follows the lines of Steps \ref{L2aS1} through \ref{L2aS4} of the proof of Lemma \ref{L2a}.

\newcounter{L1S} % proofstep = 0
\refstepcounter{L1S} % increases value by 1

{\sc Step} \arabic{L1S}\label{L1S1}\refstepcounter{L1S}.
For $\tau<T$ (with $T$ a dyadic multiple of $\tau$) we claim the formula
\begin{align}\label{wi8}
\lefteqn{\big(b\partial_1^2u_T-\sigma_i E[b,(\cdot)_T]\diam\partial_1^2v_i\big)
-\big(b\partial_1^2u_\tau-\sigma_i E[b,(\cdot)_\tau]\diam\partial_1^2v_i\big)_{T-\tau}}\nonumber\\
&=\sum_{\tau\le t<T}\Big(\big([b,(\cdot)_t]\partial_1^2u_t-\sigma_i E[b,(\cdot)_t]\partial_1^2v_{it}\big)\nonumber\\
&-[\sigma_i,(\cdot)_t]E[b,(\cdot)_t]\diam \partial_1^2v_i
-\sigma_i[E,(\cdot)_t][b,(\cdot)_t]\diam\partial_1^2v_i\Big)_{T-2t},
\end{align}
where the sum runs over $t = \frac{T}{2}, \frac{T}{4}, \ldots, \tau$.
By telescoping based on the semi-group property the formula reduces to
\begin{align*}
\lefteqn{\big(b\partial_1^2u_{2t}-\sigma_i E[b,(\cdot)_{2t}]\diam\partial_1^2v_i\big)
-\big(b\partial_1^2u_t-\sigma_i E[b,(\cdot)_t]\diam\partial_1^2v_i\big)_{t}}\nonumber\\
&=\big([b,(\cdot)_t]\partial_1^2u_t-\sigma_i E[b,(\cdot)_t]\partial_1^2v_{it}\big)\nonumber\\
&-[\sigma_i,(\cdot)_t]E[b,(\cdot)_t]\diam\partial_1^2v_i
-\sigma_i[E,(\cdot)_t][b,(\cdot)_t]\diam\partial_1^2v_i,
\end{align*}
and splits into the two statements
\begin{align}\label{wi6}
b\partial_1^2u_{2t}-(b\partial_1^2u_t)_t&=[b,(\cdot)_t]\partial_1^2u_t,\\
\sigma_i E[b,(\cdot)_{2t}]\diam\partial_1^2v_i-(\sigma_i E[b,(\cdot)_t]\diam\partial_1^2v_i)_t
&=\sigma_i E[b,(\cdot)_t]\partial_1^2v_{it}\nonumber\\
+[\sigma_i,(\cdot)_t]E[b,(\cdot)_t]\diam \partial_1^2v_i&+\sigma_i[E,(\cdot)_t][b,(\cdot)_t]\diam\partial_1^2v_i.\nonumber
\end{align}
By definition of the commutator $[\sigma_i,(\cdot)_t]$, the last statement
reduces to
\begin{align*}
E[b,(\cdot)_{2t}]\diam\partial_1^2v_i-(E[b,(\cdot)_t]\diam\partial_1^2v_i)_t
=E[b,(\cdot)_t]\partial_1^2v_{it}
+[E,(\cdot)_t][b,(\cdot)_t]\diam\partial_1^2v_i,
\end{align*}
and by the definition of $[E,(\cdot)_t]$ further to
\begin{align}\label{wi7}
[b,(\cdot)_{2t}]\diam\partial_1^2v_i-([b,(\cdot)_t]\diam\partial_1^2v_i)_t&=[b,(\cdot)_t]\partial_1^2v_{it}.
\end{align}
Now (\ref{wi6}) and (\ref{wi7}) are consequences of the semi-group property.

\medskip

{\sc Step} \arabic{L1S}\label{L1S2}\refstepcounter{L1S}.
We claim the estimate
\begin{align*}
\|\big(b\partial_1^2u_T&-\sigma_i E[b,(\cdot)_T]\diam\partial_1^2v_i\big)
-\big(b\partial_1^2u_\tau-\sigma_i E[b,(\cdot)_\tau]\diam\partial_1^2v_i\big)_{T-\tau}\|\nonumber\\
&\lesssim \big([b]_\alpha M+ N_0N_i([\sigma_i]_\alpha+\|\sigma_i\|[a]_\alpha)
\big)(T^\frac{1}{4})^{3\alpha-2}.
\end{align*}
In view of (\ref{wi8}) this estimate splits into
\begin{align}
\|[b,(\cdot)_t]\partial_1^2u_t-\sigma_i E[b,(\cdot)_t]\partial_1^2v_{it}\|&
\lesssim [b]_\alpha M (t^\frac{1}{4})^{3\alpha -2},\label{wi3}\\
\|[\sigma_i,(\cdot)_t]E[b,(\cdot)_t]\diam\partial_1^2v_i\|
&\lesssim N_0N_i[\sigma_i]_\alpha (t^\frac{1}{4})^{3\alpha -2},\label{wi1}\\
\|[E,(\cdot)_t][b,(\cdot)_t]\diam\partial_1^2v_i\|
&\lesssim N_0N_i[a]_\alpha(t^\frac{1}{4})^{3\alpha -2}.\label{w71}
\end{align}
Estimate (\ref{wi1}) follows from (\ref{wi2}) (with $\sigma_i$ playing the role of $\nu$, 
$E[b,(\cdot)_T]$ $ \diam \,\partial_1^2 v_i$ playing the role of $\tilde f$,
 and $\alpha$ playing the role of $\beta$)
and our assumption (\ref{1.5}) (without $\frac{\partial}{\partial a_0}$).
Estimate (\ref{w71}) from (\ref{w6}) (with $[b,(\cdot)_t] \diam \partial_1^2 v_i$ playing 
the role of $\tilde v$)
and our assumptions (\ref{wi5}) and (\ref{1.5}) (with $\frac{\partial}{\partial a_0}$):
\begin{align*}
\frac{1}{d^{\alpha}(x,x')}\lefteqn{\|([b,(\cdot)_t]\diam\partial_1^2v_i)(\cdot,a(x))
-([b,(\cdot)_t]\diam\partial_1^2v_i)(\cdot,a(x'))\|}\nonumber\\
&\le[a]_\alpha\sup_{a_0}\|\frac{\partial}{\partial a_0}[b,(\cdot)_t]\diam\partial_1^2v_i\|
\le[a]_\alpha N_0N_i(t^\frac{1}{4})^{2\alpha-2}.
\end{align*}
For (\ref{wi3}) we write
\begin{align}\label{wi4}
\lefteqn{\big([b,(\cdot)_t]\partial_1^2u_t-\sigma_i E[b,(\cdot)_t]\partial_1^2v_{it}\big)(x)}\nonumber\\
&=\int dy\psi_t(x-y)(b(x)-b(y))\big(\partial_1^2u_t(y)-\sigma_i(x)\partial_1^2v_{it}(y,a(x))\big)
\end{align}
and
\begin{align*}
\lefteqn{\partial_1^2u_t(y)-\sigma_i(x)\partial_1^2v_{it}(y,a(x))
=\int dz\partial_1^2\psi_t(y-z)\times}\nonumber\\
&\big(u(z)-u(x)-\sigma_i(x)(v_i(z,a(x))-v_i(x,a(x)))-\nu(x)(z-x)_1\big).
\end{align*}
Hence by the modelledness assumption of $u$, the triangle inequality $d(z,x)\le d(z,y)+d(y,x)$,
and (\ref{1.13}) we obtain
\begin{align*}
%\lefteqn{
|\partial_1^2u_t(y)-\sigma_i(x)\partial_1^2v_{it}(y,a(x))|%}\nonumber\\
&\lesssim M((t^\frac{1}{4})^{2\alpha-2}+(t^\frac{1}{4})^{-2}d^{2\alpha}(y,x)).
\end{align*}
Plugging this into (\ref{wi4}), we obtain using (\ref{1.13}) once more
\begin{align*}
\big|[b,(\cdot)_t]\partial_1^2u_t-\sigma_i E[b,(\cdot)_t]\partial_1^2v_{it}\big|(x)
\lesssim [b]_\alpha M (t^\frac{1}{4})^{3\alpha-2},
\end{align*}
as desired.

\medskip

The further two steps are as Steps \ref{L2aS3} and \ref{L2aS4} in Lemma \ref{L2a}.

\bigskip

%%%%%%%%%
%%%%%%%%%
%%%%%%%%%  new version
%%%%%%%%%

\subsection{Proof of Corollary \ref{C2}}
%\addcontentsline{toc}{section}{\protect\numberline{}Corollary \ref{C2}}

This is a corollary to Lemma \ref{L2a} in the sense that we specify the families $\{v(\cdot,x)\}_x$ and
$\{v(\cdot,x)\diamond f\}_x$ there to be given by
$\{\sigma(x)v(\cdot,a(x))\}_x$ and $\{\sigma(x)\,v(\cdot,a(x))\diamond f\}_x$, respectively.
Step \ref{L2S0} provides the necessary translations of the continuity and boundedness assumptions.
In addition, for part i) of this corollary, we need to deal with (up to second) derivatives in the parameter $a_0'$, which
on the level of Lemma \ref{L2a} is taken care of in Step \ref{S.three}. For part ii), next to the parameter derivatives,
we need to deal with differences in $f$, which is tackled in Step \ref{S.two}. Finally, for part iii),
again next to parameter derivatives, we are confronted with differences in $v$, which is taken care of
in Step \ref{S.one}. We write $[\cdot]$ for $[\cdot]_\alpha$.

\medskip

\newcounter{L2S} % proofstep = 0
\refstepcounter{L2S} % increases value by 1

{\sc Step } \arabic{L2S}.\label{S.one} \refstepcounter{L2S}
Differences in $v$ in Lemma \ref{L2a}. Suppose we are given {\it two} families of
functions $\{v_i(\cdot,x)\}_x$, $i=0,1$, and {\it two} families of distributions $\{v_i(\cdot,x)\diam f\}_x$
both satisfying (\ref{wi80}) \& (\ref{wi79}) \& (\ref{wi75}) \& (\ref{wi85}), and satisfying the analogue for
the difference, which with the abbreviations $\delta\hspace{-0.3ex}v:=v_1-v_0$,
$\delta\hspace{-0.3ex}v(\cdot,x)\diamond f:=v_1(\cdot,x)\diamond f-v_0(\cdot,x)\diamond f$ can be written as
\begin{align}
[\delta\hspace{-0.3ex}v(\cdot,x)]&\le \delta\hspace{-0.3ex}N,\label{lo10}\\
[\delta\hspace{-0.3ex}v(\cdot,x)-\delta\hspace{-0.3ex}v(\cdot,x')]&\le \delta\hspace{-0.3ex}Nd^\alpha(x,x'),\label{lo11}\\
\|[\delta\hspace{-0.3ex}v(\cdot,x),(\cdot)]\diamond f\|_{2\alpha-2}&\le \delta\hspace{-0.3ex}N N_1,\label{lo12}\\
\|[\delta\hspace{-0.3ex}v(\cdot,x ),(\cdot)]\diamond f
- [\delta\hspace{-0.3ex}v(\cdot,x'),(\cdot)]\diamond f\|_{2\alpha-2}&\le \delta\hspace{-0.3ex}N N_1 d^\alpha(x,x')\label{lo13}
\end{align}
for some constant $\delta\hspace{-0.3ex}N$. Suppose further we are given two functions $u_i$ both satisfying
(\ref{wi81}) and their difference $\delta\hspace{-0.3ex}u:=u_1-u_0$ satisfying the analogue statement
for some constant $\delta\hspace{-0.3ex}M$ and function $\delta\nu$:
\begin{align}\label{lo14}
|(\delta\hspace{-0.3ex}u(y)-\delta\hspace{-0.3ex}u(x)) %\nonumber\\
 -(\delta\hspace{-0.3ex}v(y,x)-\delta\hspace{-0.3ex}v(x,x))
-\delta\hspace{-0.3ex}\nu(x)(y-x)_1|\le\delta\hspace{-0.3ex}M d^{2\alpha}(y,x).
\end{align}
We claim that (\ref{wi83}) holds in form of
\begin{align}\label{lo17}
\|[u_1,(\cdot)]\diam f-[u_0,(\cdot)]\diam f\|_{2\alpha-2}\lesssim(\delta\hspace{-0.3ex}M+\delta\hspace{-0.3ex}N)N_1.
\end{align}
Indeed, we start by applying Lemma \ref{L2a} with $(u,v,M,N)$ replaced by
$(\delta\hspace{-0.3ex}u,\delta\hspace{-0.3ex}v,\delta\hspace{-0.3ex}M,\delta\hspace{-0.3ex}N)$:
There exists $\delta\hspace{-0.3ex}u$ such that (\ref{wi83}) takes the form
\begin{align}\label{lo16}
\|[\delta\hspace{-0.3ex}u,(\cdot)]\diamond f\|_{2\alpha-2}\lesssim(\delta\hspace{-0.3ex}M
+\delta\hspace{-0.3ex}N)N_1.
\end{align}
Note that (\ref{wi81}) holds for $(u,v,\nu)$ replaced by
$(\delta\hspace{-0.3ex}u,\delta\hspace{-0.3ex}v,\delta\hspace{-0.3ex}\nu)$,
$(u_1,v_1,\nu_1)$ and $(u_0,v_0,\nu_0)$.
Because of the definition $(\delta\hspace{-0.3ex}u,\delta\hspace{-0.3ex}v)=(u_1-u_0,v_1-v_0)$
we thus obtain from the triangle inequality that
$|\delta\nu-(\nu_1-\nu_0)|(x)|(y-x)_1|$ $\le(\delta\hspace{-0.3ex}M+2M)d^{2\alpha}(y,x)$,
which for $y\rightarrow x$ yields $\delta\nu=\nu_1-\nu_0$.
Note that (\ref{wi84}) holds with $(u,v,M,N)$ replaced by
$(\delta\hspace{-0.3ex}u,\delta\hspace{-0.3ex}v,\delta\hspace{-0.3ex}M,\delta\hspace{-0.3ex}N)$,
$(u_1,v_1,M,N)$ and $(u_0,v_0,M,N)$. Because of $(\delta\hspace{-0.3ex}u,\delta\hspace{-0.3ex}v,
\delta\hspace{-0.3ex}v\diamond f,\delta\hspace{-0.3ex}\nu)$
$=(u_1-u_0,v_1-v_0,v_1\diamond f-v_0\diamond f,\nu_1-\nu_0)$ we obtain from the triangle inequality
in $\|\cdot\|_{3\alpha-2}$
that $\lim_{T\downarrow0}\|(\delta\hspace{-0.3ex}u\diamond f-(u\diamond f_1-u\diamond f_0))_T\|=0$ and thus
$\delta\hspace{-0.3ex}u\diamond f=u_1\diamond f-u_0\diamond f$. Therefore (\ref{lo16}) turns into (\ref{lo17}).

\medskip

{\sc Step} \arabic{L2S}.\label{S.two} \refstepcounter{L2S}
Differences in $f$ in Lemma \ref{L2a}.
Suppose we are given {\it two} distributions $f_j$, $j=0,1$, and {\it two} families of distributions
$\{v(\cdot,x)\diamond f_j\}_x$ both satisfying (\ref{wi78}) \& (\ref{wi79}) \& (\ref{wi85}),
and satisfying the analogue for the difference,
which introducing the abbreviations $\delta\hspace{-0.3ex}f:=f_1-f_0$ and $v(\cdot,x)\diamond \delta\hspace{-0.3ex}f:=
v(\cdot,x)\diamond f_1-v(\cdot,x)\diamond f_0$, we may rewrite as
\begin{align}
\|\delta\hspace{-0.3ex}f\|_{\alpha-2}&\le \delta\hspace{-0.3ex}N_1,\label{lo01}\\
\|[v(\cdot,x),(\cdot)]\diamond \delta\hspace{-0.3ex}f\|_{2\alpha-2}&\le N\delta\hspace{-0.3ex}N_1,\label{lo02}\\
\|[v(\cdot,x ),(\cdot)]\diamond \delta\hspace{-0.3ex}f
 -[v(\cdot,x'),(\cdot)]\diamond \delta\hspace{-0.3ex}f\|_{2\alpha-2}&\le N\delta\hspace{-0.3ex}N_1 d^\alpha(x,x')\label{lo03}
\end{align}
for some constant $\delta\hspace{-0.3ex}N_1$. Then we claim the analogue of (\ref{wi83}), namely
\begin{align}\label{lo04}
\|[u,(\cdot)]\diamond f_1-[u,(\cdot)]\diamond f_0\|_{2\alpha-2}\lesssim(M+N)\delta\hspace{-0.3ex}N_1.
\end{align}
Indeed, from (\ref{lo01}) - (\ref{lo03}) together with the remaining assumptions of Lemma \ref{L2a}
we learn from the latter that there exists a distribution we call $u\diamond\delta\hspace{-0.3ex}f$
such that (\ref{wi84}) holds with $(f,N_1)$ replaced by $(\delta\hspace{-0.3ex}f,\delta\hspace{-0.3ex}N_1)$.
Since it also holds with $(f_j,N_1)$, we obtain from the triangle inequality and the above definition
of $v(\cdot,x)\diamond \delta\hspace{-0.3ex}f$ that $\lim_{T\downarrow 0}\|(u\diamond\delta\hspace{-0.3ex}f-
(u\diamond f_1-u\diamond f_0))_T\|=0$, which gives $u\diamond\delta\hspace{-0.3ex}f=u\diamond f_1-u\diamond f_0$
and thus (\ref{wi83}), still with $(f,N_1)$ replaced by $(\delta\hspace{-0.3ex}f,\delta\hspace{-0.3ex}N_1)$,
turns into (\ref{lo04}).

\medskip

{\sc Step} \arabic{L2S}.\label{S.three} \refstepcounter{L2S}
$C^m$-dependence of $(f,v\diamond f)$ on a parameter $a_0'\in[\lambda,\frac{1}{\lambda}]$ in Lemma \ref{L2a}.
Suppose $f\in C_{a_0'}^m(C^{\alpha-2})$ and that $[v(\cdot,x),(\cdot)]\diamond f$ is of class $C^m_{a_0'}(C^{2\alpha-2})$
uniformly in $x$, see below for the precise meaning.
We claim that this is preserved: $[u,(\cdot)]\diam f$ is of class $C^m_{a_0'}(C^{2\alpha-2})$.
Moreover, if (\ref{wi78}) \& (\ref{wi79}) \& (\ref{wi85}) are strengthened to
\begin{align}
\|f\|_{\alpha-2,m}&\le N_1,\label{lo05}\\
\|[v(\cdot,x),(\cdot)]\diamond f\|_{2\alpha-2,m}&\le N N_1,\label{lo06}\\
\|[v(\cdot,x ),(\cdot)]\diamond f
 -[v(\cdot,x'),(\cdot)]\diamond f\|_{2\alpha-2,m}&\le N N_1 d^\alpha(x,x'),\label{lo07}
\end{align}
cf. (\ref{conv1}) \& (\ref{conv2}), then (\ref{wi83}) improves likewise:
\begin{align}\label{li08}
\|[u,(\cdot)]\diamond f\|_{2\alpha-2,m}\lesssim(M+N) N_1.
\end{align}

\medskip

In virtue of Lemma \ref{L2a} and fixing $a_0'$ and $j$, we may associate to
$\big((\frac{\partial}{\partial a_0'})^jf,
\{(\frac{\partial}{\partial a_0'})^j(v(\cdot,x)\diamond f)\}_x\big)$ a
distribution we call $u\diamond(\frac{\partial}{\partial a_0'})^jf$. Under the assumptions of
Lemma \ref{L2a} enhanced by (\ref{lo05})-(\ref{lo07}), (\ref{wi83}) turns into
\begin{align}\label{lo09}
\|[u,(\cdot)]\diamond (\frac{\partial}{\partial a_0'})^jf\|_{2\alpha-2}\lesssim(M+N)N_1.
\end{align}
It is convenient to abbreviate by ${\mathcal R}f(\tilde a_0',a_0')$
$:=f(\tilde a_0')-\sum_{j=0}^m(\tilde a_0'-a_0')^j(\frac{\partial}{\partial a_0'})^jf(a_0')$ Taylor's remainder
for a generic (Banach space-valued) function $f$ of $a_0'$. Our $C^m$-assumption on the input includes
$\lim_{\tilde a_0'\rightarrow a_0'}$ $\|{\mathcal R}f(\tilde a_0',a_0')\|_{\alpha-2}=0$
and $\lim_{\tilde a_0'\rightarrow a_0'}\sup_{x}\|{\mathcal R}([v(x,\cdot),(\cdot)]\diamond f)(\tilde a_0',a_0')\|_{2\alpha-2}$ $=0$.
From the latter we learn that $u\diamond f\in C^m_{a_0'}(C^{\alpha-2})$
with $(\frac{\partial}{\partial a_0'})^j\big(u\diamond f\big)$ $=u\diam(\frac{\partial}{\partial a_0'})^jf$,
so that in particular (\ref{lo09}) turns into (\ref{li08}).
From the former we therefore learn that
$\lim_{\tilde a_0'\rightarrow a_0'}\|{\mathcal R}([u,(\cdot)]\diam f)(\tilde a_0',a_0')\|_{2\alpha-2}=0$, so that
the $C^m_{a_0'}(C^{2\alpha-2})$ property is transmitted.

\medskip

{\sc Step} \arabic{L2S}.\label{L2S0} \refstepcounter{L2S} Some algebra. 
Suppose that $\{v(\cdot,a_0)\}_{a_0}$ and $\{v_i(\cdot,a_0)\}_{a_0}$, $i=0,1$, are three families of functions and
$|\cdot|$ a semi-norm on functions of $x$ (like $[\cdot]$) such that
%%
%\begin{align}
%\sup_{a_0}|\{1,\frac{\partial}{\partial a_0}\}v(\cdot,a_0)|\le N_0,\label{wj1}\\
%\sup_{a_0}|\{1,\frac{\partial}{\partial a_0},\frac{\partial^2}{\partial a_0^2}\}v_{i}(\cdot,a_0)|
%\le N_0,\label{wj2}\\
%\sup_{a_0}|\{1,\frac{\partial}{\partial a_0}\}(v_{1}-v_{0})(\cdot,a_0)|
%\le \del N_0\label{wj3}
%\end{align}
%
\begin{align}
|v|_{1}\le N_0,\label{wj1}\\
|v_{i}|_2
\le N_0,\label{wj2}\\
|v_{1}-v_{0}|_{1}
\le \del N_0\label{wj3}
\end{align}
for some constants $N_0,\del N_0$ (here as in \eqref{conv1} the subscripts  in $|\cdot|_1$ and $|\cdot|_2$
refer to the number of parameter derivatives with respect to $a_0$). 
The reason for this more general framework
is useful because in Step \ref{S.five} we apply it with $v(\cdot,a_0)$ replaced by $[v(\cdot,a_0),(\cdot)_T]\diam f$
and with the supremum norm $N_1^{-1}(T^\frac{1}{4})^{2-2\alpha}\|\cdot\|$ playing the role of $|\cdot|$.
%and with $(N_0,\del N_0)$ replaced by $(N_1N_0(T^\frac{1}{4})^{2\alpha-2},$
%$N_1\del N_0(T^\frac{1}{4})^{2\alpha-2})$. 
We claim that this entails
\begin{align}
|&\sigma(x)v(\cdot,a(x))|\le N_0\|\sigma\|,\label{wi92}\\
|&\sigma(x)v(\cdot,a(x))-\sigma(x')v(\cdot,a(x'))|\le
N_0([\sigma]+\|\sigma\|[a])d^\alpha(x,x'),\label{wi41}\\
|&\sigma_1(x)v_1(\cdot,a_1(x))-\sigma_0(x)v_0(\cdot,a_0(x))|\nonumber\\
& \qquad \le N_0(\|\sigma_1-\sigma_0\|+\max_i\|\sigma_i\|\|a_1-a_0\|)+\del N_0\max_{i}\|\sigma_i\|,\label{wi93}\\
|&\big(\sigma_1(x)v_1(\cdot,a_1(x))-\sigma_0(x)v_0(\cdot,a_0(x))\big) \nonumber\\
& \qquad \qquad  -\big(\sigma_1(x')v_1(\cdot,a_1(x'))-\sigma_0(x')v_0(\cdot,a_0(x'))\big)|\nonumber\\
& \qquad \le\Big(N_0\max_{i,j}
\big([\sigma_1-\sigma_0]+\|\sigma_i\|[a_1-a_0]+[\sigma_i]\|a_1-a_0\|\nonumber\\
&\qquad +\|\sigma_1-\sigma_0\|[a_i]+\|\sigma_i\|[a_j]\|a_1-a_0\|\big)\nonumber\\
&\qquad +\del N_0 \max_{i}\big([\sigma_i]+\|\sigma_i\|[a_i]\big)\Big)d^\alpha(x,x').\label{wi42}
\end{align}
Estimate (\ref{wi92}) follows immediately from (\ref{wj1}).
%estimates  \eqref{wi41} and \eqref{wi93} follow from \eqref{wj1}--\eqref{wj3}
%by applying the discrete Leibniz rule.
We treat (\ref{wi41}), (\ref{wi93}), and (\ref{wi42}) along the same lines, which is a bit of an overkill
for (\ref{wi41}) and (\ref{wi93}).
%For \eqref{wi42} we start with the elementary, and purposefully symmetric, formula
We start with the two elementary, and purposefully symmetric, formulas
\begin{align}\label{wi46}
\sigma v-\sigma'v'=\frac{1}{2}(\sigma-\sigma')(v+v')+\frac{1}{2}(\sigma+\sigma')(v-v'),
\end{align}
%
%and
%
\begin{align}\label{wi59}
\lefteqn{(\sigma_1v_1-\sigma_0v_0)-(\sigma_1'v_1'-\sigma_0'v_0')}\nonumber\\
&=\frac{1}{4}((\sigma_1-\sigma_0)-(\sigma_1'-\sigma_0'))(v_1+v_1'+v_0+v_0')\nonumber\\
&+\frac{1}{4}((\sigma_1+\sigma_1'+\sigma_0+\sigma_0'))((v_1-v_0)-(v_1'-v_0'))\nonumber\\
&+\frac{1}{4}((\sigma_1-\sigma_1')+(\sigma_0-\sigma_0'))((v_1-v_0)+(v_1'-v_0'))\nonumber\\
&+\frac{1}{4}((\sigma_1-\sigma_0)+(\sigma_1'-\sigma_0'))((v_1-v_1')+(v_0-v_0')).
\end{align}

  We use the first formula twice. The first application is for
   $\sigma=\sigma(x)$ and $\sigma'=\sigma(x')$, $v=v(\cdot,a(x))$, and
   $v'=v(\cdot,a(x'))$ to obtain using the triangle inequality %
   \begin{align*}
    \lefteqn{|\sigma(x) v(\cdot,a(x))-\sigma(x')v'(\cdot,a(x'))|}\nonumber\\
    &\le[\sigma]d^\alpha(x,x')\sup_{a_0}|v(\cdot,a_0)|
    +\|\sigma\|\sup_{a_0}|\frac{\partial v}{\partial a_0}(\cdot,a_0)|[a]d^\alpha(x,x').
    \end{align*}
    In view of the assumption (\ref{wj1}) this yields (\ref{wi41}).
    The second application is for
    $\sigma=\sigma_1(x)$ and $\sigma'=\sigma_0(x)$, $v=v_1(\cdot,a_1(x))$, and
    $v'=v_0(\cdot,a_0(x))$. We obtain the inequality
    \begin{align}\label{wi94}
    \lefteqn{|\sigma_1(x)v_1(\cdot,a_1(x))-\sigma_0(x)v_0(\cdot,a_0(x))|}\nonumber\\
    &\le\|\sigma_1-\sigma_0\|\max_i\sup_{a_0}|v_i(\cdot,a_0)|+\max_i\|\sigma_i\||v_1(\cdot,a_1(x))-v_0(\cdot,a_0(x))|.
    \end{align}
    In view of the assumption (\ref{wj2}), the first right-hand-side term is estimated as desired. For the second rhs term
    we interpolate linearly in the sense of $v_s$ $:=sv_1$ $+(1-s)v_0$ and $a_s$ $:=sa_1$ $+(1-s)a_0$, to the effect of
    \begin{align}\label{wi45}
    \lefteqn{v_1(\cdot,a_1(x))-v_0(\cdot,a_0(x))}\nonumber\\
    &=\int_0^1ds\big((v_1-v_0)(\cdot,a_s(x))+\frac{\partial v_s}{\partial a_0}(\cdot,a_s(x))(a_1-a_0)(x),
    \end{align}
    from which we learn
    \begin{align}\label{wi43}
    |v_1(\cdot,a_1(x))-v_0(\cdot,a_0(x))|\le\sup_{a_0}|v_1-v_0|+\max_i\sup_{a_0}|\frac{\partial v_i}{\partial a_0}|\|a_1-a_0\|.
    \end{align}
    Inserting this into (\ref{wi94}) and in view of the assumption (\ref{wj2})\&(\ref{wj3}) we obtain
    the remaining part of (\ref{wi93}).
    
%\medskip

We use the second  formula \eqref{wi59} 
for $\sigma_i=\sigma_i(x)$, $\sigma_i'=\sigma_i(x')$, $v_i=v_i(\cdot,a_i(x))$, and
$v_i'=v_i(\cdot,a_i(x'))$ to obtain
\begin{align*}
|&\big(\sigma_1(x)v_1(\cdot,a_1(x))-\sigma_0(x)v_0(\cdot,a_0(x))\big)%\nonumber\\
-\big(\sigma_1(x')v_1(\cdot,a_1(x'))-\sigma_0(x')v_0(\cdot,a_0(x'))\big)|\nonumber\\
&\le[\sigma_1-\sigma_0]d^\alpha(x,x')\max_i\sup_{a_0}|v_i(\cdot,a_0)|\nonumber\\
&+\max_i\|\sigma_i\||(v_1(\cdot,a_1(x))-v_0(\cdot,a_0(x)))-(v_1(\cdot,a_1(x'))-v_0(\cdot,a_0(x')))|\nonumber\\
&+\max_i[\sigma_i]d^\alpha(x,x')\sup_y|v_1(\cdot,a_1(y))-v_0(\cdot,a_0(y))|\nonumber\\
&+\|\sigma_1-\sigma_0\|\max_i\sup_{a_0} \; \Big|\frac{\partial v_i}{\partial a_0}(\cdot,a_0) \Big| \;[a_i]d^\alpha(x,x').
\end{align*}
In order to deduce (\ref{wi42}) from this inequality, in view of (\ref{wi43})
and of our assumption (\ref{wj2}) \& (\ref{wj3}),
it remains to show for the second right-hand-side terms
\begin{align}
|&(v_1(\cdot,a_1(x))-v_0(\cdot,a_0(x)))-(v_1(\cdot,a_1(x'))-v_0(\cdot,a_0(x')))|\nonumber\\
&\le\sup_{a_0} \Big|\frac{\partial}{\partial a_0}(v_1-v_0)(\cdot,a_0) \Big|\max_i[a_i]d^\alpha(x,x')\nonumber\\
&+\max_i\sup_{a_0} \Big|\frac{\partial^2 v_i}{\partial a_0^2}(\cdot,a_0) \Big|\max_j[a_j]d^\alpha(x,x')\|a_1-a_0\|\nonumber\\
&+\max_i\sup_{a_0}\Big|\frac{\partial v_i}{\partial a_0}(\cdot,a_0)\Big|[a_1-a_0]d^\alpha(x,x').\label{wi44}
\end{align}
We appeal again to the outcome (\ref{wi45}) of the linear interpolation,
which immediately yields the first right-hand-side term (\ref{wi44}) from the first right-hand-side term
in (\ref{wi45}). For the second right-hand-side term in (\ref{wi44}), we appeal  once more to formula
(\ref{wi46}) (applied to $\sigma$ $=\frac{\partial v_s}{\partial a_0}(\cdot,a_s(x))$,
$\sigma'$ $=\frac{\partial v_s}{\partial a_0}(\cdot,a_s(x'))$, $v$ $=(a_1-a_0)(x)$, and
$v'$ $=(a_1-a_0)(x')$).

\medskip

%%%%%%%%%%%%%%%%%%%%%%%%%%END IGNORE%%%%%%%%%%%%%%%%%%%%%%%%%%%%%%%%%%%%%%%%%%%%%%%%%%%%%%%%%%%%%%%%%%%%%%%%%%%%%%

{\sc Step} \arabic{L2S}.\label{S.five}\refstepcounter{L2S} Conclusion
We start with part i) of this corollary; we apply Lemma \ref{L2a}, in form of Step \ref{S.three}
with $m=2$, to the families given by distributions $\{f(\cdot,a_0')\}_{a_0'}$,
the functions $\{\sigma(x)v(\cdot,a(x))\}_x$, and the products $\{\sigma(x)\,v(\cdot,a(x))\diamond f(\cdot,a_0')\}_{x,a_0'}$.
To this purpose we verify the hypotheses; hypothesis (\ref{lo05}) on the distribution $f$ is identical
to the corollary's hypothesis (\ref{wj36}).
We now turn to those on the function $v$, namely (\ref{wi75}) and (\ref{wi80}).
Using (\ref{wk58}), these follow, with $N_0$ playing the role of $N$, from (\ref{wi92}) and (\ref{wi41}) of Step \ref{L2S0}
provided the generic semi-norm $|\cdot|$ there is chosen to be $[\cdot]$.
The relevant hypothesis (\ref{wj1}) of Step \ref{L2S0} is identical to the
corollary's hypothesis (\ref{L0.1}).
We last turn to the hypothesis on the product $v\diamond f$,
that is, (\ref{lo06}) \& (\ref{lo07});
to this purpose, we fix a convolution parameter $T\in(0,1]$,
an order of differentiation $j=0,1,2$ and the parameter $a_0'$.
These hypotheses follow again from (\ref{wi92}) and (\ref{wi41}) of Step \ref{L2S0}, this time
with $\{(\frac{\partial}{\partial a_0'})^j[v(\cdot,a_0),(\cdot)_T]\diamond f(\cdot,a_0')\}_{a_0}$
playing the role of $\{v(\cdot,a_0)\}_{a_0}$ and the norm $N_1^{-1}(T^\frac{1}{4})^{2-2\alpha}\|\cdot\|$
replacing $|\cdot|$. The relevant hypothesis (\ref{wj1})
then holds by the corollary's hypothesis (\ref{wj41}); the outputs (\ref{wi92}) \& (\ref{wi41}) indeed
turn into (\ref{lo06}) \& (\ref{lo07}), still with $N_0$ playing the role of $N$.
Finally, the outcome (\ref{li08}) of Step \ref{S.three} turns into the desired (\ref{wj40}).

\medskip

We now turn to part ii) of this corollary. Again, we apply Lemma \ref{L2a}, this time in form of Step \ref{S.two},
upgraded by Step \ref{S.three} with $m=1$ in the sense that the expressions $(\|\cdot\|_{\alpha-2},\|\cdot\|_{2\alpha-2})$
are replaced by $(\|\cdot\|_{\alpha-2,1},\|\cdot\|_{2\alpha-2,1})$. The argument follows the lines
of the one for part i): When it comes to the product $(\frac{\partial}{\partial a_0'})^m(v\diamond f_1-v\diamond f_0)$,
for fixed $m=0,1$ and parameter $a_0'$,
the presence of an $a_0$-derivative in the corollary's hypothesis (\ref{wj48})
feeds into Step \ref{L2S0}'s hypothesis (\ref{wj1}) with the semi-norm
$|\cdot|$ $=(\delta\hspace{-0.3ex}N_1)^{-1}(T^\frac{1}{4})^{2-2\alpha}\|\cdot\|$.
Step \ref{L2S0}'s output (\ref{wi92}) \& (\ref{wi41}) provides Step \ref{S.two}'s input (\ref{lo02}) \& (\ref{lo03}).
Step \ref{S.two}'s output (\ref{lo04}) is identical to the corollary's claim (\ref{wk35}).

\medskip

We finally turn to part iii) of this corollary. A last time, we apply Lemma \ref{L2a}, now in form of Step \ref{S.one},
upgraded in terms of differentiability in the parameter $a_0'$ by Step \ref{S.three} with $m=1$.
We apply Step \ref{S.one} to
the families given by distributions $\{f(\cdot,a_0')\}_{a_0'}$,
the functions $\{\sigma_i(x)v_i(\cdot,a_i(x))\}_x$,
and the products $\{\sigma_i(x)\,v(\cdot,a_i(x))\diamond f(\cdot,a_0')\}_{x,a_0'}$.
We start with the hypotheses (\ref{lo10}) \& (\ref{lo11}) on the difference of the functions
and apply Step \ref{L2S0} to $|\cdot|=[\cdot]$: The relevant input (\ref{wj2}) \& (\ref{wj3}) of that step is provided
by the corollary's assumptions (\ref{wi20}) \& (\ref{wi21}). In view of (\ref{wk58}), the output
(\ref{wi93}) \& (\ref{wi42}) of Step \ref{L2S0} turns into the hypotheses (\ref{lo10}) \& (\ref{lo11})
with $\delta\hspace{-0.3ex}N:=N_0([\sigma_1-\sigma_0]+\|\sigma_1-\sigma_0\|+[a_1-a_0]+\|a_1-a_0\|)+\delta\hspace{-0.3ex}N_0$.
We now turn to the hypotheses (\ref{lo12}) \& (\ref{lo13}) on the difference of the products
and apply Step \ref{L2S0} to $\{(\frac{\partial}{\partial a_0'})^jv_i(\cdot,a_0)\diam f(\cdot,a_0')\}_{a_0}$, $i=0,1$,
playing the role of $\{v_i(\cdot,a_0)\}_{a_0}$ and with
$|\cdot|=N_1^{-1}(T^\frac{1}{4})^{2-2\alpha}\|\cdot\|$
for fixed $T$, $j=0,1$ and $a_0'$:
The relevant input (\ref{wj2}) \& (\ref{wj3}) of that step is provided
by the corollary's assumptions (\ref{wj50}) \& (\ref{wj49}). The output
(\ref{wi93}) \& (\ref{wi42}) of Step \ref{L2S0} turns into the hypotheses (\ref{lo12}) \& (\ref{lo13}) with
the above definition of $\delta\hspace{-0.3ex}N$.
Finally, we note that the modelledness assumption of our corollary assumes the form (\ref{lo14}).
The output (\ref{lo17}) of Step \ref{S.one} turns into the desired (\ref{wk36}).

\subsection{Proof of Corollary \ref{Coro}}
%\addcontentsline{toc}{section}{\protect\numberline{}Corollary \ref{Coro}}
%{\sc Proof of Corollary \ref{Coro}}. 

\newcounter{LCoroS} % proofstep = 0
\refstepcounter{LCoroS} % increases value by 1

{\sc Step} \arabic{LCoroS}\label{LCoroS1}\refstepcounter{LCoroS}. Proof of (i) $\Rightarrow$ (ii).
As $v$ is a $C^{\beta+2}$ function the assumption that $u$ is modelled after $v$ according to $a(u)$, $\sigma(u)$ implies that $u$ is of class $C^{2\alpha }$, in particular $\partial_1 u$ is a function of class $C^{2\alpha -1}$ (of course, as we will see below, $u$ is actually of class $C^{\beta+2}$ but 
we do not have this information to our disposal yet). Together with the regularity assumption on $f$ this implies that there is a classical interpretation of the products $\sigma(u) f$ and $a(u) \partial_1^2 u$ the latter as a distribution. In fact, this is obvious for 
$\sigma(u) f$ and for $a(u) \partial_1^2u$ we can set, for example,
\begin{align}\label{classical-Produkt11}
a(u) \partial_1^2 u := \partial_1 (a(u) \partial_1 u) - \partial_1 a(u) \partial_1 u.
\end{align}
The claim  then follows from standard parabolic regularity theory as soon as we have established that 
\begin{align}
\label{AMclaim1}
\sigma(u) \diam f &= \sigma(u) f + \sigma'(u) \sigma(u) g_1( \cdot, a(u))\\
\label{AMclaim2}
a(u) \diam \partial_1^2 u &= a(u) \partial_1^2 u + a'(u) \sigma^2(u) g_2( \cdot, a(u), a(u)).
\end{align}
\medskip

We first argue that \eqref{AMclaim1} holds. To see this, first by Lemma~\ref{L0} $\sigma(u)$
is modelled after $v$ according to $a(u)$ and $\sigma'(u) \sigma(u)$.
Then, Corollary~\ref{C2} characterizes $\sigma(u) \diam f$ as the unique distribution  
for which 
\begin{align}\label{AM99}
\lim_{T \downarrow 0}\|[\sigma(u), ( \cdot)_T] \diam f - \sigma'(u) \sigma(u)E [v, ( \cdot)_T] \diam f - \nu [x_1, (\cdot)] f\| =0.
\end{align}
By the $C^\beta$ regularity of $f$ as well as the $C^{2 \alpha}$ regularity of $\sigma(u)$ one sees immediately that each of the commutators in this expression 
goes to zero if  $\diam$ is replaced by the classical product
\begin{align*}
\|[\sigma(u), ( \cdot)_T]  f \| , \; \| \sigma'(u) \sigma(u)E [v, ( \cdot)_T]  f \| , \;  \| \nu [x, (\cdot)] f\|  \to 0
\end{align*}
for $T \to 0$.
Hence \eqref{AM99} turns into 
\begin{align*}
\lim_{T \downarrow 0} \| \sigma(u) f - (\sigma(u) \diam f)_T - \sigma'(u) \sigma(u) g_{1,T}( \cdot, a(u)) \| =0.
\end{align*}
Since, $g( \cdot, a_0) \in C^{\beta}$ by assumption, this yields \eqref{AMclaim1}.
 In the same way, one can see that for any $a_0'$ we have
\begin{align}\label{AM101}
a(u) \diam \partial_1^2 v( \cdot, a_0') = a(u) \partial_1^2 v(\cdot, a_0') + a'(u) \sigma(u) g_2 ( \cdot, a(u), a_0')
\end{align}
(the classical definition of $ a(u) \partial_1^2 v(\cdot, a_0')$ poses no problem because $v $ is of class $C^{\beta+2}$). 

\medskip

It remains to upgrade \eqref{AM101} to \eqref{AMclaim2}, i.e. the second factor $\partial_1^2 v$ in \eqref{AM101} should be replaced by $\partial_1^2 u$. To this end 
we make the ansatz 
\begin{align}\label{ansatz-bad-B}
a(u) \diam \partial_1^2 u &= a(u) \partial_1^2 u+ a'(u) \sigma^2(u) g_2( \cdot, a(u), a(u)) +B,
\end{align}
and aim to show that $B=0$. Recalling once more that $u$ is modelled after $v$ according to  $a(u)$, $\sigma(u)$ we 
 invoke Lemma~\ref{L1} and plug in our ansatz \eqref{ansatz-bad-B} to obtain 
\begin{align} \notag
\lim_{T \downarrow 0}  &\| [a(u), (\cdot)_T]  \partial_1^2 u - ( a'(u) \sigma^2(u) g_2( \cdot, a(u), a(u)))_T +(B)_T \\
\label{AM307}
& \qquad - \sigma(u)E [a(u), (\cdot)_T] \diam \partial_1^2 v \| =0.
\end{align}
%
%As before, the classical commutator $\sigma(u)E [a(u), (\cdot)_T] \partial_1^2 v$ goes to zero uniformly in $x$ as $T \downarrow 0$. 
Plugging \eqref{AM101} into \eqref{AM307} we obtain
\begin{align} \notag
\lim_{T \downarrow 0}  &\| [a(u), (\cdot)_T]  \partial_1^2 u - ( a'(u) \sigma^2(u) g_2( \cdot, a(u), a(u)))_T +(B)_T \\
\label{AM306}
&  - \sigma(u)E [a(u), (\cdot)_T]  \partial_1^2 v -   a'(u) \sigma^2(u) E (g_2 ( \cdot, a(u), a_0'))_T \| =0.
\end{align}
Now according to our regularity assumptions we have both
\begin{align*}
\|  ( a'(u) \sigma^2(u) g_2( \cdot, a(u), a(u)))_T - a'(u) \sigma^2(u) E (g_2 ( \cdot, a(u), a_0'))_T \| &\to 0\\
 \, \| \sigma(u)E [a(u), (\cdot)_T]  \partial_1^2 v \| &\to 0,
\end{align*}
for $T \to 0$, which reduces \eqref{AM306} to 
\begin{align*}
\lim_{T \downarrow 0}\| [a(u), (\cdot)_T]  \partial_1^2 u - B_T \| =0,
\end{align*}
where we recall that the classical commutator is defined based on \eqref{classical-Produkt11}. Now, according to its definition \eqref{classical-Produkt11} 
we have $ [a(u), (\cdot)_T] \rightharpoonup 0 $, which characterizes $B$ as $0$.

\medskip
{\sc Step} \arabic{LCoroS}\label{LCoroS2}\refstepcounter{LCoroS}. Proof of (ii) $\Rightarrow$ (i). If $u$ as well as all the $v( \cdot, a_0)$ are  of class $C^{\beta +2}$, then 
$u$ is automatically modelled after $v$ according to $a(u)$ and $\sigma(u)$. Thus we 
can conclude from Step~\ref{LCoroS1} that \eqref{AMclaim1} and \eqref{AMclaim2} hold 
which in turn implies that $u$ solves $\partial_2u-P(a(u)\diam \partial_1^2u+\sigma(u)\diam f)=0$ distributionally.

\bigskip

%%%%%%%%%%%%%%%%%%%%%%%%%%%%%%%%%%%%%%%%%%%%%%%%%%%%%%%%%%%%%%%%%%%%%%%%%%%%%%%%%%%%%%%%%%%%%%%%
%%%%%%%%%%%%%%%%%%%%%%%%%%%%%%%%%%%%%%%%%%%%%%%%%%%%%%%%%%%%%%%%%%%%%%%%%%%%%%%%%%%%%%%%%%%%%%%%
%%%%%%%%%%%%%%%%%%%%%%%%%%%%%%%%%%%%%%%%%%%%%%%%%%%%%%%%%%%%%%%%%%%%%%%%%%%%%%%%%%%%%%%%%%%%%%%%

\section{Proofs of the stochastic bounds}
\label{s:stoch-proofs}

\subsection{Proof of Lemma \ref{lem5}}
%\addcontentsline{toc}{section}{\protect\numberline{}Lemma \ref{lem5}}

%{\sc Proof of Lemma \ref{lem5}}.

\newcounter{Llem55S} % proofstep = 0
\refstepcounter{Llem55S} % increases value by 1

{\sc Step} \arabic{Llem55S}\label{Llem55S1}\refstepcounter{Llem55S}. Proof of \eqref{SB2}.
By stationarity of $f_T = f \ast \psi_T$ we have for $T \leq 1$
%
%Assumption \eqref{A1} and the stationarity and periodicity of $f$ imply that for $T\leq 1$ 
%
\begin{align*}
&\langle f_T^2(0) \rangle 
 = \Big\langle \int_{[0,1)^2} f_T^2  \, dx \Big\rangle 
\overset{\eqref{SB1}}{=} \sum_{k \in (2\pi \Z)^2 }  \hat{\psi}^2_T (k) \hat{C}(k)\\
& \overset{\eqref{SB2}}{\leq}   (T^{\frac14})^{-3+\lambda_1 + \lambda_2 }  \sum_{k \in (2\pi \Z)^2 \setminus\{0\} }  (T^{\frac14})^3  \frac{e^{- 2(T^{\frac14} k_1)^4 - 2(T^{\frac12} k_2)^2}}{ (T^{\frac14 }(1+|k_1|))^{\lambda_1} (T^{\frac12 }(1+|k_2|))^{\frac{\lambda_2}{2}}} \\
&\lesssim   (T^{\frac14})^{2 \alpha'-4}.
\end{align*}
In the last estimate we have used that  for $T \downarrow 0$  the sum in the third line is a Riemann sum approximation of the  integral 
$\int  e^{-2\hat{k}_1^4 - 2\hat{k}_2^2} |\hat{k}_1|^{-\lambda_1}|\hat{k}_2|^{-\frac{\lambda_2}{2}} d\hat{k}$ which converges due to $\lambda_1, \frac{\lambda_2 }{2}< 1$.

\medskip

The fact that $f_T$ is Gaussian and stationary implies that we have $\langle |f_T(x)|^p \rangle \lesssim \langle f^2_T(0) \rangle^{\frac{p}{2}}$, which  permits to write
\begin{align*}
\Big\langle \int_{[0,1)^2} |f_T|^p dx \Big\rangle^{\frac{1}{p}} \lesssim \langle f_T^2(0) \rangle^{\frac{1}{2}} \lesssim   (T^{\frac14})^{ \alpha'-2}.
\end{align*}
In order to upgrade this $L^p$ bound to an $L^\infty$ bound  under the expectation  we
observe that by the semi-group property \eqref{1.10} we have $f_T = (f_{T/2})_{T/2}$ 
such that H\"older's inequality implies
\begin{align*}
\| f_T \| \lesssim  \| f_{T/2} \|_{L^p}  \| \psi_{T/2, \text{per}} \|_{L^{p'}} 
\end{align*}
where as before $\| \cdot \|$ refers to the supremums norm over $\R^2$  (or equivalently $[0,1)^2$ by periodicity) and $\| \cdot \|_{L^p}$ refers to the $L^p$ norm over $[0,1)^2$, $p' := \frac{p}{p-1}$ is the dual exponent of $p$, and  $\psi_{T, \text{per}}(x) = \sum_{z \in \Z^2} \psi_T(x+z)$ is the periodization of $\psi_T$.  By observing that for small $T$ the difference $\big| \| \psi_{T, \text{per}} \|_{L^{p'}} - \big(\int_{\R^2} |\psi_{T}|^{p'}dx \big)^{\frac{1}{p'}}\big|$ stays bounded, and scaling we get 
 $ \| \psi_{T, \text{per}} \|_{L^{p'}} \lesssim (T^{\frac14})^{-\frac{3}{p}}$ such that finally
\begin{align*}
\big\langle \|f_T\|^p \big\rangle^{\frac{1}{p}} &\lesssim (T^{\frac14})^{-\frac{3}{p}}  \big\langle  \|f_T\|_{L^p}^p \big\rangle^{\frac{1}{p}}   
\lesssim  (T^{\frac14})^{ \alpha'-2-\frac{3}{p}}.
\end{align*}
To also accommodate for the supremum over the scales $T$ we first note that  $\| f_{T+t}\| \lesssim \| f_T\|$ implies 
\begin{equation*}
\| f\|_{\alpha-2}   = \sup_{T \leq 1}(T^{\frac14})^{2-\alpha} \|f_T \| \lesssim \sup_{T \leq 1, \mathrm{dyadic}}(T^{\frac14})^{2-\alpha} \|f_T \|,
\end{equation*}
where the subscript dyadic means that this supremum is only taken over dyadic $T$. Then we write
\begin{align*}
\Big\langle  \Big(\sup_{T \leq 1, \mathrm{dyadic}}(T^{\frac14})^{2-\alpha} \|f_T \| \Big)^p \Big\rangle & \leq \sum_{T \leq 1, \mathrm{dyadic}} (T^{\frac14})^{p(2-\alpha)}   \big\langle \|f_T \|^p \big\rangle \\
&\lesssim 
\sum_{T \leq 1, \mathrm{dyadic}} (T^{\frac14})^{p(2-\alpha)}  (T^{\frac14})^{p( \alpha'-2)-3},
\end{align*}
which converges as soon as $p>\frac{3}{\alpha' - \alpha}$ and thus establishes
$
 \big\langle   \| f \|_{\alpha-2}^p \big\rangle^{\frac{1}{p}} \lesssim 1
$
 for large $p$. The same bound for smaller $p$ can be derived from the bound for large $p$ and Jensen's inequality.
Finally, because of $(f_\eps)_T=  \varphi_\eps \ast f_T$ and because the operators $\varphi_\eps \ast$ are bounded 
with respect to $\| \cdot \|$ uniformly in $\eps$, the bound holds uniformly in the regularization leading to the 
desired estimate \eqref{SB2}.

\medskip

{\sc Step} \arabic{Llem55S}\label{Llem55S2}\refstepcounter{Llem55S}. Proof of \eqref{SB2A}.
The bound on the $\eps$-differences follows from \eqref{SB2} as soon as   we have established the deterministic bound 
\begin{align}\label{e-diff1}
 \| (f_\eps)_T - f_T \| \lesssim \min \Big\{ \Big( \frac{\eps}{T}\Big)^{\frac14}, 1\Big\}  \| f_{T/2}\|,
\end{align}
which by the semi-group property reduces to 
\[
 \| (f_\eps)_T - f_T \| \lesssim \min \Big\{ \Big( \frac{\eps}{T}\Big)^{\frac14}, 1\Big\}  \| f \|.
\]
Since   $( \cdot)_T$ and $\varphi_\eps \ast $ are bounded with respect to $\| \cdot \|$, it suffices 
to consider $\eps \leq T$. 
%
%
% Given the  bound  
%\begin{align}\label{e-diff2}
% \| (f_\eps)_T - f_T \| \leq  \| (f_\eps)_T \| +  \| f_T \|  \lesssim  \| f_{T/2}\|
%\end{align}
%%
%which follows from the triangle inequality, the semi-group property in the form $(f_\eps)_T = (f_{T/2})_{T/2} \ast \ph_{\eps}$ and $f_T = (f_{T/2})_{T/2}$
%as well as the fact that the  operators $(\cdot)_{T/2}$ and $\ph_\eps \ast$ are bounded with respect to $\| \cdot \|$ uniformly  in $T$ and $\eps$, it suffices to 
%establish \eqref{e-diff1} for $\eps \leq \frac14 T$. 
We then write
\begin{align*}
 \| (f_\eps)_T - f_T \| & = \| (\psi_{T} \ast \ph_{\eps} - \psi_{T} ) \ast f \| 
 \leq \int_{\R^2} |\psi_{T} \ast \ph_{\eps} - \psi_{T} | dx \; \| f\|,
\end{align*}
and have thereby reduced \eqref{e-diff1} (and hence \eqref{SB2A}) to establishing that 
\begin{align}\notag
\int_{\R^2} |\psi_{T} \ast \ph_{\eps} - \psi_{T} | dx  \lesssim \Big( \frac{\eps}{T}\Big)^{\frac14} \qquad \text{for } \eps \leq T.
\end{align}
By scaling (recalling that $\psi_T(x_1,x_2) = T^{-\frac34} \psi_1(T^{-\frac14}x_1, T^{-\frac12}x_2)$), it suffices to show this bound for $T=1$,
in which case it turns into 
\begin{align*}
\int_{\R^2} |\psi_1 \ast  \ph_\eps - \psi_1| \lesssim \eps^{\frac14} \qquad \text{for $\eps \leq 1$}
\end{align*}
 which is immediate for Schwartz kernels $\psi_1$, $\varphi$ 
and in view of the definition \eqref{i-def-fepsveps} of $\varphi_\eps$.
%(x_1,x_2):=\frac{1}{\eps^\frac{3}{4}}\varphi ( \frac{x_1}{\eps^\frac{1}{4}},\frac{x_2}{\eps^\frac{1}{2}})$.

\bigskip

\subsection{Proof of Lemma~\ref{cch1}}
%\addcontentsline{toc}{section}{\protect\numberline{}Lemma~\ref{cch1}}
%{\sc Proof of Lemma \ref{cch1}}.
%
% 
%
\newcounter{Lcch1S} % proofstep = 0
\refstepcounter{Lcch1S} % increases value by 1

For $a_0 \in [\lambda,\frac{1}{\lambda}]$ let $G( \cdot, a_0)$ be the (periodic) Green function of $( \partial_2 -a_0 \partial_1^2 )$, where the heat operator is endowed with periodic and zero average time-space boundary conditions. Its Fourier series 
is given by 
\begin{equation}\label{e:Greens}
\hat{G}( k, a_0) = 
\begin{cases}
\frac{1}{ a_0 k_1^2 -  i k_2} =  \frac{a_0  k_1^2 + i k_2}{ a_0^2 k_1^4 +  k_2^2} \qquad &\text{ for } k \in (2 \pi \Z)^2 \setminus \{0\}, \\
0 \qquad &\text{ for } k = 0.
\end{cases}
\end{equation}
With this notation in place, $v(\cdot, a_0)$ 
is characterized by its discrete Fourier transforms 
 $\hat{v} (k, a_0) = \hat{G}( k, a_0) \hat{f}(k)$. 
% $\hat{v}_\eps (k , a_0)= \hat{G}( k, a_0) \hat{f}_\eps(k)$ for $k \in (2 \pi \Z)^2$. 
Throughout the proof the parameter dependence on $a_0$ only appears in $\hat{G}(k, a_0)$
for which only the bound 
\begin{equation}\label{e:GreensBound}
|\hat{G}(k,a_0)| \lesssim \frac{1}{k_1^2 + |k_2|}
\end{equation}
is used. We thus suppress the $a_0$-dependence in all expressions.

\medskip

{\sc Step} \arabic{Lcch1S}\label{Lcch1S1}\refstepcounter{Lcch1S}. Bound on the expectation.
We claim that
\begin{align}\label{psn1}
 \langle [v' , (\cdot)_T]  \diam f' \rangle \lesssim  (T^{\frac14})^{2 \alpha' -2 - \kappa_1 - \kappa_2}.
 \end{align}
By stationarity $\langle (v'  \diam f' )_T \rangle  = \langle v'  \diam f'  \rangle =0$. Furthermore, by stationarity
and \eqref{SB1}, \eqref{A1} we have
\begin{align}
\notag
&\big| \langle v' f'_T  \rangle \big| = \Big| \sum_k \langle \hat{v'}(-k) \hat{\psi}_T(k) \hat{f'}(k) \rangle \Big| 
%\\
%\notag
%&
= \Big| \sum_k (\hat{M}_2\hat{G})(-k)  \hat{\psi}_{T}(k) \hat{M}_1(k)   \hat{C}(k) \Big|\\
\notag
&\overset{\eqref{mn00},\eqref{e:GreensBound}}{\leq}   (T^{\frac14})^{-3+2 + \lambda_1 + \lambda_2 - \kappa_1 - \kappa_2} \\
\notag
& \qquad  \times  \sum_k (T^{\frac14})^3\frac{\hat\psi_{T}(k)}{( T^{\frac14} k_1)^2 +|T^{\frac12}k_2|}   \frac{(( T^{\frac14} k_1)^4 +( T^{\frac12}k_2)^2)^{\frac{\kappa_1 + \kappa_2}{4}}}{(T^{\frac14}(1+|k_1|))^{\lambda_1} (T^{\frac12}(1+|k_2|))^{\frac{\lambda_2}{2}}}\\
\notag
&\lesssim   (T^{\frac14})^{ 2\alpha'-2-\kappa_1 - \kappa_2},
\end{align}
where the sum is taken over $(2 \pi \Z)^2 \setminus\{0\}$.
In the last step we have used the fact that the Riemann sum in the third line approximates the integral 
$\int \hat\psi_1(\hat{k}) \frac{(\hat{k}_1^4 + \hat{k}_2^2)^{\frac{\kappa_1 + \kappa_2}{4}}}{\hat{k}_1^2 + \hat{k}_2} \frac{1}{|\hat{k}_1|^{\lambda_1} |\hat{k}_2|^{\lambda_2/2}}d\hat{k}$. This integral 
converges because the singularities on the axes $\hat{k}_1 = 0$ and $\hat{k}_2=0$ are integrable because of $\lambda_1, \frac{\lambda_2}{2}<1$ and the singularity near the origin is integrable due to $2 + \lambda_1 + \lambda_2 = 1+2\alpha' <3$, where we appeal to the fact that the parabolic dimension is $3$ (alternatively,
one may split the integral into $|x_1| \leq \sqrt{|x_2|}$ and its complement). 
 This establishes \eqref{psn1}.
 
 \medskip
{\sc Step} \arabic{Lcch1S}\label{Lcch1S2}\refstepcounter{Lcch1S}. Preparation for bound on the variance.
For the variances we seek the bound
\begin{align*}
\Big| \langle ([v', (\cdot)_T] \diam f')^2 \rangle -  \langle v' f'_T  \rangle^2 \Big|^{\frac12} \lesssim (T^{\frac14})^{2\alpha'-2 -\kappa_1 - \kappa_2 } ,
\end{align*}
which by definition of $\diam$ can be expressed equivalently without the renormalization as
\begin{align}\label{SB11}
 \Big| \langle ([v', (\cdot)_T] f')^2 \rangle  - \langle [v', (\cdot)_T]  f' \rangle^2  \Big|^{\frac12}\lesssim (T^{\frac14})^{2\alpha'-2- \kappa_1 - \kappa_2}.
\end{align}
To  derive the estimate in the form \eqref{SB11} we write using  once more  stationarity
\begin{align*}
\langle ([v', (\cdot)_T] f')^2 \rangle &= \Big\langle  \int_{[0,1)^2} ([v', (\cdot)_T] f')^2 dx \Big\rangle 
= \sum_{k \in (2 \pi \Z)^2} \Big\langle  |\reallywidehat{[v', (\cdot)_T] f'}(k)|^2  \Big\rangle.
\end{align*}
The expression appearing in the last expectation can be evaluated according to its definition
\begin{align}
%\notag
 \reallywidehat{[v', (\cdot)_T] f'}(k) 
\notag
=& \sum_{\ell \in ( 2 \pi \Z)^2} ( \hat{\psi}_T(\ell) - \hat{\psi}_T(k) ) \hat{v'}(k-\ell) \hat{f'}(\ell)  \\
\label{SBB2}
=& \sum_{\ell \in ( 2 \pi \Z)^2} ( \hat{\psi}_T(\ell) - \hat{\psi}_T(k) ) (\hat{M_2}\hat{G})(k-\ell) \hat{f}(k-\ell)\hat{M}_1(\ell) \hat{f}(\ell),
\end{align}
which permits to write
\begin{align}
\notag
\langle ([v', (\cdot)_T] f')^2 \rangle
%\notag
 & = \sum_{k} \sum_{\ell} \sum_{\ell'} ( \hat{\psi}_T(\ell) - \hat{\psi}_T(k) )  ( \hat{\psi}_T(-\ell') - \hat{\psi}_T(-k) )\\
 \notag
&\quad \times (\hat{M_2} \hat{G})(k-\ell) (\hat{M_2}\hat{G})(-(k-\ell') ) \hat{M_1}(\ell) \hat{M_1}(-\ell')\\
\label{L1-1}
& \quad \times \big\langle \hat{f}(k-\ell) \hat{f}(\ell)  \hat{f}(-(k-\ell')) \hat{f}(-\ell') \big \rangle,
\end{align}
where all sums are taken over $(2\pi \Z)^2$.
We now use \eqref{SB1} and the Gaussian identity
\begin{align}
 \label{L1-2}
&\langle  \hat{f}(k-\ell) \hat{f}(\ell) \hat{f}(-(k-\ell')) \hat{f}(-\ell')\rangle  \\
&  =  \delta_{k,0} \hat{C}(\ell)  \hat{C}(\ell')   
\notag
 + \delta_{\ell, \ell'} \hat{C}(k-\ell)  \hat{C}(\ell)  + \delta_{k-\ell, \ell'} \hat{C}(k-\ell) \hat{C}(\ell) .
\end{align}
Plugging this identity into  \eqref{L1-1} results in three terms which we bound one by one.
%
%
%and bound the three terms resulting from plugging this identity into one by one. 
The first term coincides with the square of the expectation  
(which is subtracted on the left hand side of \eqref{SB11})
 \begin{align*}
&\sum_{\ell} \sum_{\ell'} ( \hat{\psi}_T(\ell) - \hat{\psi}_T(0) )  ( \hat{\psi}_T(-\ell') - \hat{\psi}_T(0) ) \\
& \qquad (\hat{M}_2 \hat{G})(-\ell) (\hat{M}_2\hat{G})(\ell' )(\hat{M}_1\hat{C})(\ell) (\hat{M}_1\hat{C})(-\ell') \\
& = \Big( \sum_{\ell} ( \hat{\psi}_T(\ell) - \hat{\psi}_T(0) ) (\hat{M}_2 \hat{G})(-\ell) (\hat{M}_1\hat{C})(\ell) \Big)^2 
=  \langle  [v', (\cdot)_T] f' \rangle^2  
\end{align*}
so that the required bound \eqref{SB11} follows as soon as we can bound the remaining two terms.
The term originating from the third 
contribution on the right hand side of \eqref{L1-2} can be absorbed into the second term 
$\sum_{k,\ell}  
	( \hat{\psi}_T(\ell) - \hat{\psi}_T(k) )^2  
	|(\hat{M}_2\hat{G})(k-\ell )|^2 
	|\hat{M}_1(\ell)|^2  \hat{C}(\ell)\hat{C}(k-\ell)$
using the Cauchy-Schwarz inequality. Indeed, we may write 
\begin{align}
\notag
&\sum_{k} \sum_{\ell}  ( \hat{\psi}_T(\ell) - \hat{\psi}_T(k) )  ( \hat{\psi}_T(-(k-\ell)) - \hat{\psi}_T(-k) )\\
\notag
&\qquad \qquad \times  (\hat{M}_2 \hat{G})(k-\ell) (\hat{M}_2 \hat{G})(-\ell ) (\hat{M}_1\hat{C})(\ell) (\hat{M}_1\hat{C})(-(k-\ell))\\
  \notag
 &  \leq \big( \sum_{k,\ell}  ( \hat{\psi}_T(\ell) - \hat{\psi}_T(k) )^2  
 |(\hat{M}_2\hat{G})(k-\ell )|^2 |\hat{M}_1(\ell)|^2  \hat{C}(-(k-\ell))\hat{C}(\ell) \Big)^{\frac12}\\
 \notag
 &  \qquad \times \Big( \sum_{k,\ell} ( \hat{\psi}_T(- (k-\ell )) - \hat{\psi}_T(-k) )^2
  |(\hat{M}_2\hat{G})(-\ell)|^2 \\
  \notag
  & \qquad \qquad\qquad  \qquad\times |\hat{M}_1(-(k-\ell))|^2 \hat{C}(-(k-\ell))\hat{C}(\ell) \Big)^\frac{1}{2},
\end{align}
and the second factor on the right hand side can be seen to coincide with the first one by performing the change of 
variables $k' = -k$ and $\ell' =  \ell-k$ and the symmetry $\hat C(k) = \hat{C}(-k)$.  Hence, it only remains to bound the term coming from 
the second contribution on the right hand side of \eqref{L1-2}. We use the assumptions~\eqref{A1} and \eqref{mn00}
to bound this term as follows
\begin{align}
\notag&
\sum_{k,\ell}  
	( \hat{\psi}_T(\ell) - \hat{\psi}_T(k) )^2  
	|(\hat{M}_2\hat{G})(k-\ell )|^2 
	|\hat{M}_1(\ell)|^2  \hat{C}(\ell)\hat{C}(k-\ell)\\
\notag &
\overset{\eqref{mn00}\eqref{e:GreensBound}}{\lesssim  }
\sum_{k \neq \ell} 
	\frac{
		( \hat{\psi}_T(\ell) - \hat{\psi}_T(k) )^2  
	}{ 
		(|(k-\ell)_1^2 + |(k-\ell)_2| )^{2}
	}  
	\frac{
		(\ell_1^4 + \ell_2^2)^{\frac{\kappa_1}{2}}
	}{ 
		(1+|\ell_1|)^{\lambda_1}  (1+|\ell_2|)^{\frac{\lambda_2}{2}}
	}\\
\notag& \qquad \qquad 
\times  
	\frac{
		( (k-\ell)_1^4 +  (k-\ell)_2^2  )^{\frac{\kappa_2}{2}}
	}{ 
		(1+|(k-\ell)_1)|)^{\lambda_1}  (1+|(k-\ell)_2|)^{\frac{\lambda_2}{2}}
	} \\
\notag & 
= 
(T^{\frac14})^{4\alpha' -4 - 2 \kappa_1 - 2 \kappa_2}   
\sum_{k \neq \ell} 
	\big(T^{\frac14} \big)^6 
	\Big( 
		\frac{   
			\hat{\psi}_T(\ell) - \hat{\psi}_T(k)  
		}{ 
			( T^{\frac14}(k-\ell)_1)^2 + |T^{\frac12}(k-\ell)_2|  
		} 
	\Big)^2 \\
 &\qquad  \qquad \times
 \notag
  \frac{
  		(
			(T^{\frac14}\ell_1)^4   +   
			( T^{\frac12}  \ell_2)^2
		)^{\frac{\kappa_1}{2}}
	}{
		(T^{\frac14}  (1+ |\ell_1|)  )^{\lambda_1}  (T^{\frac12}  (1+|\ell_2|)  )^{\frac{\lambda_2}{2}}
	}  \\
  \label{SBBB4}
& \qquad \qquad \times    
	\frac{
		((T^{\frac14}(k-\ell)_1)^4 
		+ (T^{\frac12}(k-\ell)_2)^2)^{\frac{\kappa_2}{2}}
	}{
		(T^{\frac14}(1+ |(k-\ell)_1|))^{\lambda_1}  
		(T^{\frac12}  (1+|(k-\ell)_2|)  )^{\frac{\lambda_2}{2}}
	}   .
\end{align}

\medskip

{\sc Step} \arabic{Lcch1S}\label{Lcch1S3}\refstepcounter{Lcch1S}. Bound on an integral.
In order to show that the expression \eqref{SBBB4} is bounded by $\lesssim   (T^{\frac14})^{4\alpha' -4 - 2 \kappa_1 - 2 \kappa_2}  $ which in turn establishes \eqref{SB11},
it remains to show the convergence of the integral which is  approximated
by the Riemann sum in the last lines:
\begin{align}\label{impint}
&\int \!\! \int  \Big(\frac{\hat{\psi}_1(\hat{\ell}) - \hat{\psi}_1(\hat{k})}{(\hat{k}-\hat{\ell})_1^2 + |\hat{k}_2-\hat{\ell}_2|} \Big)^2
\bar{C_1}(\hat{\ell}) \bar{C}_2(\hat{k}-\hat{\ell}) d\hat{\ell} d\hat{k},
%% 
%& \qquad \qquad  \qquad \qquad   \times \frac{((\hat{k}-\hat{\ell})_1^4 + (\hat{k}-\hat{\ell})_2^2)^{\frac{\kappa_2}{2}}}{|\hat{k}_1-\hat{\ell}_1|^{\lambda_1}|\hat{k}_2-\hat{\ell}_2|^{\frac{\lambda_2}{2}}} d\hat{\ell} d\hat{k},
%&=\int \!\!\int  \Big(\frac{\hat{\psi}_1(\hat{\ell}) - \hat{\psi}_1(\hat{\ell}+\hat{h})}{\hat{h}_1^2 + |\hat{h}_2|} \Big)^2 \frac{(\hat{\ell}_1^4 + \hat{\ell}_2^2)^{\frac{\kappa_1}{2}}}{|\hat{\ell}_1|^{\lambda_1}|\hat{\ell}_2|^{\frac{\lambda_2}{2}}} \frac{(\hat{h}_1^4 + \hat{h}_2^2)^{\frac{\kappa_2}{2}}}{|\hat{h}_1|^{\lambda_1}|\hat{h}_2|^{\frac{\lambda_2}{2}}}  d\hat{\ell} d\hat{h} .
\end{align}
where momentarily we use the short-hand
\[
\bar{C}_i(\hat{\ell}):= \frac{(\hat{\ell}_1^4 + \hat{\ell}_2^2)^{\frac{\kappa_i}{2}}}{|\hat{\ell}_1|^{\lambda_1}|\hat{\ell}_2|^{\frac{\lambda_2}{2}}}  \qquad i=1,2.
\]
As a first step we deal with the integral near the diagonal, where $|(\hat{\ell}-\hat{k})_1| + |(\hat{\ell}-\hat{k})_2|\leq1$. For these values the change of variables $\hat{h} = \hat{k}-\hat{\ell}$ is useful.
We furthermore make use of the bound $|\hat{\psi}_1(\hat{\ell}) - \hat\psi_1(\hat{\ell}+\hat{h})| \lesssim (|\hat{h}_1| + |\hat{h}_2| ) \int_0^1 | \nabla \hat\psi_1(\hat{\ell} + \theta \hat{h} )| d \theta$
and brutally bound $ (|\hat{h}_1| + |\hat{h}_2| )  \leq \sqrt{\hat{h}_1^2 +|\hat{h}_2|}$
 so that we need to address the convergence of  
 \[
 \int \!\! \int_{|\hat{h}_1|+|\hat{h}_2| \leq 1}  \frac{\max_{\theta \in [0,1]}  |\nabla \hat{\psi}_1|^2(\ell + \theta \hat{h})}{\hat{h}_1^2 + |\hat{h}_2|} 
\bar{C_1}(\hat{\ell}) \bar{C}_2(\hat{h})  d\hat{h} d\hat{\ell}.
 \]
 For the $d\hat{h}$ integral over a finite volume it suffices to assert that the singularities near the axes $\hat{h}_1=0$ and $\hat{h}_2=0$ are integrable  
 due to   $\lambda_1, \frac{\lambda_2}{2}<1$
 and that the singularity near the origin is integrable because by assumption~\eqref{A1} $2+\lambda_1+\lambda_2=1+2\alpha'$ which is less than the parabolic 
 dimension $3$. The singularities for the $d\hat{\ell}$ integral are only better behaved  and the convergence of the integral for $|\hat{\ell}| \to \infty$
 is guaranteed by the exponential decay of $\nabla \hat{\psi}_1$. 
 
We now discuss the convergence of \eqref{impint} for $|(\hat{\ell}-\hat{k})_1| + |(\hat{\ell}-\hat{k})_2|> 1$: For these values we write $(\hat{\psi}_1(\hat{\ell}) - \hat{\psi}_1(\hat{k}))^2 \lesssim \hat{\psi}_1^2(\hat{\ell}) + \hat{\psi}_1^2(\hat{k})$
and treat the resulting integrals separately. For the integral coming from $\hat{\psi}_1^2(\hat{\ell}) $ we use the same change of variables $\hat{h} = \hat{k}-\hat{\ell}$ which leads us to consider the integral
\begin{align*}
%&\int \!\! \int_{|(\ell-k)_1| + |(\ell-k)_2|> 1}  \frac{\hat{\psi}_1^2(\hat{\ell})}{((\hat{k}-\hat{\ell})_1^2 + |\hat{k}_2-\hat{\ell}_2|)^2} 
%\bar{C_1}(\hat{\ell}) \bar{C}_2(\hat{k}-\hat{\ell}) d\hat{\ell} d\hat{k}\\
\int \!\! \int_{|\hat{h}_1| + |\hat{h}_2|> 1}  \frac{\hat{\psi}_1^2(\hat{\ell})}{(\hat{h}_1^2 + |\hat{h}_2|)^2} 
\bar{C_1}(\hat{\ell}) \bar{C}_2(\hat{h}) d\hat{\ell} d\hat{h}.
\end{align*}
 As above, the $d \hat{\ell}$ integral converges because the singularities of $\bar{C}_1(\ell)$ near the axes $\ell_1=0$ and $\ell_2=0$ as 
 well as the singularity near the origin are integrable and because of the exponential decay of $\hat{\psi}_1$ at infinity. The singularities 
 of the $d\hat{h}$ integral near the axes are also integrable and its convergence for $|\hat{h}| \to \infty$ is guaranteed by 
 the fact that $4+ \lambda_1 + \lambda_2 = 3+2\alpha'$ which is larger than the parabolic dimension $3$ and by $\kappa_2\ll 1$.
  
  It remains to treat the integral coming from $\hat{\psi}_1^2(\hat{k}) $:
\begin{align*}
&\int \!\! \int_{|(\ell-k)_1| + |(\ell-k)_2|> 1}  \frac{\hat{\psi}_1^2(\hat{k})}{((\hat{k}-\hat{\ell})_1^2 + |\hat{k}_2-\hat{\ell}_2|)^2} 
\bar{C_1}(\hat{\ell}) \bar{C}_2(\hat{k}-\hat{\ell}) d\hat{\ell} d\hat{k}.
 \end{align*} 
 It is here that our assumption $\alpha>\frac14$ becomes relevant to assure the convergence of the $d\hat{\ell}$ integral.
We get
 \begin{align*}
& \int_{|(\hat{\ell}-\hat{k})_1| + |(\hat{\ell}-\hat{k})_2|> 1}  \frac{1}{((\hat{k}-\hat{\ell})_1^2 + |\hat{k}_2-\hat{\ell}_2|)^2} 
\bar{C_1}(\hat{\ell}) \bar{C}_2(\hat{k}-\hat{\ell}) d\hat{\ell},
 \end{align*} 
  which converges for $|\hat{\ell}| \to \infty$ because of $4+ 2(\lambda_1 +\lambda_2)= 2+4\alpha'$ which is larger than the parabolic 
  dimension $3$ due to $\alpha'>\frac14$ and because $\kappa_1,\kappa_2\ll1$.  The
   convergence of the resulting $d\hat{k}$ integral near the origin is guaranteed by $2(\lambda_1+\lambda_2) -3 = 4\alpha' -5<3$
   and   
   for $|\hat{k}| \to \infty$ by the exponential
  decay of $\hat{\psi}_1^2(\hat{k})$.

\subsection{Proof of Corollary~\ref{cch2}}
%\addcontentsline{toc}{section}{\protect\numberline{}Corollary~\ref{cch2}}
%{\sc Proof of Corollary~\ref{cch2}}.

%We aim at giving a meaning to  the products $v(\cdot, a_0)\diam f$, $v(\cdot, a_0)\diam \partial_1^2 v(\cdot, a_0')$ and obtaining bounds for the families of 
%commutators $[v ( \cdot,a_0), (\cdot)_T ]\diam  f$, $[v (\cdot,a_0), ( \cdot)_T]\diam \partial_1^2 v( \cdot, a_0')$ derived from them. 

The quantity $\partial_1^2 v(\cdot, a_0)$ is obtained from $f$ through a regularity-preserving transformation, as can be expressed in terms of the Fourier transform
\begin{align*}
\widehat{\partial_1^2 v}(k, a_0) =  \frac{ k_1^2}{a_0k_1^2 - i k_2} \hat{f}(k).
\end{align*}
%
%and noting that $ \frac{ k_1^2}{a_0 k_1^2 - i k_2}$ is a bounded symbol (see also Lemma \ref{LA1}). Therefore, the proofs for  $v(\cdot, a_0) \diam f$ and $v(\cdot, a_0) \diam \partial_1^2 v(\cdot, a_0')$ are essentially identical. 

Derivatives with respect to $a_0$ and $a_0'$  do not change the regularity either as can be seen from
% for example we have for any $n \geq 1$
%
\begin{align}\label{BoundedFourier}
\Big(\frac{\partial }{\partial a_0} \Big)^n\hat{G}(k, a_0)   =
 \frac{(-1)^n n! k_1^{2n}}{(a_0 k_1^2 - i k_2)^n} \hat{G}(k, a_0)  \qquad n \geq 1
\end{align}
and for every $n$ the symbol $\frac{(-1)^n n! k_1^{2n}}{(a_0 k_1^2 - i k_2)^{n}}$ is also bounded. Therefore,
the estimate \eqref{mn02} follows immediately from \eqref{mn01} either with $f_\eps$ in the role of $f'$ (i.e. $\hat{M}_1 = \hat{\ph}_{\eps_1}$) 
or
$ \big(\frac{\partial}{\partial a_0'}\big)^{m}$ $ \partial_1^2 v_\eps(\cdot, a_0')$ in the role of $f'$ which amounts to 
\begin{align*}
\hat{M}_1(k) =   \frac{(-1)^m m! k_1^{2m}}{(a_0 k_1^2 - i k_2)^m}  \frac{-k_1^2}{a_0' k_1^2- i k_2} \hat{\ph}_{\eps_1}(k)
\end{align*}
and with $ \big(\frac{\partial}{\partial a_0}\big)^{n}$ $  v_{\eps_0}(\cdot, a_0)$ in the role of $v'$ i.e. 
\begin{align*}
\hat{M}_2(k) =   \frac{(-1)^n n! k_1^{2n}}{(a_0' k_1^2 - i k_2)^n}  \hat{\ph}_{\eps_0}(k).
\end{align*}
\medskip

For the derivatives with respect to $\eps_i$ 
%we observe that the product rule applies in the form
%%
%\begin{align*}
%&\frac{\partial}{\partial \eps} 
%\Big(
%	\Big[ 
%		\big(\frac{\partial}{\partial a_0}\big)^n v_\eps( \cdot, a_0),(\cdot)_T 
%	\Big] 
%	\diam \{  f_\eps, \big(\frac{\partial}{\partial (a_0')\big)^{n'}} \partial_1^2 v_\eps \}
%\Big) \\
%&= 
% \Big[ 
%	\frac{\partial}{\partial \eps}  \frac{\partial}{\partial a_0^n}   v_\eps( \cdot, a_0),(\cdot)_T 
%\Big] 
%	\diam \{  f_\eps, \frac{\partial}{\partial (a_0')^{n'}} \partial_1^2 v_\eps \}
%\\
% & \qquad +  
% \Big[ 
% 	\frac{\partial}{\partial a_0^n} v_\eps( \cdot, a_0),(\cdot)_T 
%\Big] 
%\diam \frac{\partial}{\partial \eps} 
%	\{  
%		f_\eps, \frac{\partial}{\partial (a_0')^{n'}} \partial_1^2 v_\eps 
%	\}. 
%\end{align*}
%%
%We then apply \eqref{mn02} to each of the terms on the right hand side separately. 
the multipliers $\hat{M}_1$, $\hat{M}_2$ are the same as above only with $\hat{\ph}_{\eps_i}$ replaced by 
$|\eps_i \frac{\partial}{\partial \eps_i} \hat{\ph}_{\eps_i} |\lesssim  \big( (k_1^4 + k_2^2) \eps \big)^\frac{\kappa_i}{4} $ 
in $\hat{M}_2$ for $i=0$ and  in $\hat{M}_1$  if $i=1$.

%\bigskip
%
%
%
\subsection{Proof of Proposition~\ref{P2}}
%\addcontentsline{toc}{section}{\protect\numberline{}Proposition \ref{P2}}
%{\sc Proof of Proposition~\ref{P2}}.
%
%
%
\newcounter{P2St} % proofstep = 0
\refstepcounter{P2St} % increases value by 1

{\sc Step} \arabic{P2St}\label{P2St1}\refstepcounter{P2St}. Bound on the supremum over $x$ and $T$.
Our first claim is that for all $a_0, a_0' \in [\lambda,\frac{1}{\lambda}]$, $\eps_0,\eps_1 \in (0,1]$, for $\kappa \ll 1$ and for all $n, m\geq 1$ and $i=0,1$ we have
\begin{align}
\notag
 & \Big\langle 
 	\Big( 
		\sup_{T \leq 1}   (T^{\frac14})^{2-2\alpha} 
		\Big\| 
			\Big(\frac{\partial}{\partial a_0}\Big)^n 
			\Big(\frac{\partial}{\partial {a_0'}}\Big)^m 
			[ v_{\eps_0}( \cdot, a_0),(\cdot)_T ]\\
 \label{yoo1}
  & \qquad \qquad  \qquad
  			\diam \{ f_{\eps_1}, \partial_1^2 v_{\eps_1}(\cdot, a_0')  \} 
		\Big\| 
	\Big)^p 
\Big\rangle^{\frac{1}{p}} 
\lesssim 1,\\
\notag
 & \Big\langle 
 	\Big( 
	\sup_{T \leq 1}   (T^{\frac14})^{2-2\alpha + \kappa} 
	\Big\| 
		\eps_i \frac{\partial}{\partial \eps_i} 
		\Big(\frac{\partial}{\partial a_0}\Big)^n 
		\Big(\frac{\partial}{\partial a_0'}\Big)^m 
		[ v_{\eps_0}( \cdot, a_0),(\cdot)_T ]\\
 \label{yoo2}
  & \qquad \qquad  \qquad 
  		\diam \{ f_{\eps_1}, \partial_1^2 v_{\eps_1}(\cdot, a_0')  \} 
	\Big\| 
	\Big)^p 
\Big\rangle^{\frac{1}{p}} 
\lesssim \eps_i^{\frac{\kappa}{4}} \qquad \text{for all $p<\infty$}.
\end{align}
To keep the notation concise, for the moment we restrict ourselves to the bound for 
$[v_{\eps_0}, (\cdot)_T ] \diam \partial_1^2 v_{\eps_1}$ 
without the derivatives with respect to $a_0$, $a_0'$, $\eps_i$. The general case of \eqref{yoo1} follows in the identical way and so does \eqref{yoo2} if in the proof \eqref{mn02} is replaced by \eqref{mn03}. To simplify the notation further we drop the subscript $\eps_i$
as well as the dependence on $a_0, a_0'$ for the moment.

\medskip
First of all $[v, (\cdot)_T] \diam \partial_1^2 v$ 
is a random variable in the second Wiener chaos over the Gaussian field $f$ such that by equivalence of moments (see e.g. \cite[Chapter 1]{Nualart}, \cite[Section 1.6]{Bogachev}, or \cite[Section 3]{MourratWeberXu})   for random variables in the second Wiener chaos and by stationarity, the bound \eqref{mn02}  can be upgraded to 
\begin{align}\label{P1-1}
\langle  | [v, (\cdot)_T] \diam \partial_1^2 v|^{p} \rangle^{\frac{1}{p}} 
\lesssim (T^{\frac14})^{2\alpha'-2}  \qquad \text{for all $p<\infty$}.
\end{align}
We now aim to upgrade this $L^p$ bound to an $L^\infty$ bound over $x$. At the same time, we want to show that the supremum over all $T \leq 1$ can be reduced to a supremum over all dyadic $T$. For any given $T \leq 1$  there is a unique a dyadic $T' \leq \frac12$ such that $T = 2 T' + t$  with $2T' \leq T < 4T'$ and we refer to this choice when we write $T'$ in the sequel. 

\medskip
We make use of the commutator identity \eqref{v36} in the form of
\begin{align}\label{P1-2}
[v,(\cdot)_{T}] \diam \partial_1^2 v  
= \big( [v, (\cdot)_{T'}] \diam \partial_1^2 v \big)_{T'+t} 
	+ [v, (\cdot )_{T'+t} ] (\partial_1^2 v)_{T'}.
\end{align}
The second term on the right hand side of \eqref{P1-2} can be bounded directly by making the convolution with $\psi_{T'+t}$ explicit
%and 
% using Lemma~\ref{LA1} and \eqref{1.13} to get
%
\begin{align*}
&\Big|\big([v, (\cdot )_{T'+t} ] (\partial_1^2 v)_{T'}\big)(x) \Big| \\
&= \Big| \int (v(x) - v(y) ) \psi_{T'+t}(x-y) (\partial_1^2 v)_{T'}(y) dy \Big|\\
&\leq [ v]_{{\alpha}} \| (\partial_1^2 v)_{T'}\| \int d^{\alpha} (x,y) |\psi_{T'+t}(x-y)| dy\\
&\overset{\eqref{1.13}}{\lesssim}  ({T'}^\frac14)^{{\alpha}-2} ((T'+t)^\frac14)^{{\alpha}} [ v]_{\alpha}^2   
\overset{\eqref{w31}}{\lesssim} (T^\frac14)^{2{\alpha}-2}     \|f\|_{\alpha-2}^2.
\end{align*}
%
%(Note that here we have used Lemma~\ref{LA1}  to bound the semi-norm $[v]_\alpha$ by the $C^{\alpha - 2}$ 
%norm of $f$. The same argument applies to $\frac{\partial}{\partial a_0^n} v( \cdot, a_0)$ using an identity
%like \eqref{x90} and iterating Lemma~\ref{LA1}. Similarly the $C^{\alpha -2}$ norm of $\frac{\partial}{\partial a_0^n} \partial_1^2 v( \cdot, a_0)$ is controlled by the $C^{\alpha -2}$ norm of $f$.)
Derivatives with respect to $a_0$, $a_0'$ can be dealt with as in Step~\ref{C1S0} of the proof of Corollary~\ref{C1}.

Taking the sup over $x$ and $T$ and then the $p$-th moment in the expectation we get from Lemma~\ref{lem5}
\begin{align}
\notag
\big \langle \big( \sup_{T \leq 1 }(T^{\frac14})^{2-2\alpha} \| [v, (\cdot )_{T'+t} ] (\partial_1^2 v)_{T'} \| \big)^p \big\rangle^{\frac{1}{p}} 
\notag
&\lesssim   \big \langle  \|f\|_{\alpha-2}^{2p}   \big\rangle^{\frac{1}{p}}
\lesssim 1 .
\end{align}

\medskip

 To bound the first term on the right hand side of \eqref{P1-2}  we use Young's inequality (on the torus)
 in the form
 \begin{align*}
\|\big( [v, (\cdot)_{T'}] \diam \partial_1^2 v \big)_{T'+t} \| 
\lesssim
	\|  [v, (\cdot)_{T'}] \diam \partial_1^2 v  \|_{L^p} 
	\| \psi_{T'+t, \text{per}} \|_{L^{p'}} ,
 \end{align*}
where we use the notation of Step~\ref{Llem55S1} in the proof of Lemma~\ref{lem5}, resulting in
%
% where $L^p$ denotes the $L^p$ norm on the torus $[0,1)^2$, $p' = \frac{p}{p-1}$ is the dual exponent of $p$ and 
% $\psi_{T, \text{per}}(x) = \sum_{k \in \Z^2} \psi_T(x+k)$ is the periodisation of $\psi_T$.  By observing that for for small $T$ the difference 
% $\big| \| \psi_{T, \text{per}} \|_{L^{p'}} - \big(\int_{\R^2} \psi_{T}^{p'}dx \big)^{\frac{1}{p'}}\big|$ 
% stays bounded and scaling we get 
% $ \| \psi_{T, \text{per}} \|_{L^{p'}} \lesssim (T^{\frac14})^{-\frac{3}{p}}$ 
% such that finally
%
\begin{align*}
\|\big( 
	[v, (\cdot)_{T'}] \diam \partial_1^2 v 
\big)_{T'+t} \| 
\lesssim  ((T'+t)^{\frac14})^{-\frac{3}{{p}}}  
	\|  
		 [v, (\cdot)_{T'}] \diam \partial_1^2 v  
	\|_{L^{p}}.
\end{align*}
Taking the supremum over $T$ we get for any $p$
\begin{align*}
&\Big( \sup_{T \leq 1} (T^{\frac14})^{2 - 2\alpha}   
		\| \big( [v, (\cdot)_{T'}] \diam \partial_1^2 v \big)_{T'+t} \| \Big)^p \\
&\lesssim 
\sum_{T' \leq \frac12, \mathrm{dyadic}}
	  ({T'}^{\frac14})^{p(2 - 2\alpha)}   ((T')^{\frac14})^{-3}  
		\|  
			 [v, (\cdot)_{T'}] \diam \partial_1^2 v  
		\|_{L^{p}}^p.
\end{align*}
Finally, we take the expectation of this estimate  and use \eqref{P1-1} and the stationarity to get
\begin{align*}
&\Big \langle  
	\Big(
		\sup_{T \leq 1} (T^{\frac14})^{2 - 2\alpha}   
		\| \big( [v, (\cdot)_{T'}] \diam \partial_1^2 v \big)_{T'+t} \| 
	\Big)^{p} 
\Big\rangle \\
&\lesssim 
\sum_{T' \leq \frac12, \mathrm{dyadic}}({T'}^{\frac14})^{p(2 - 2\alpha)} ({T'}^{\frac14})^{-3}  
	\big\langle   
		\| [v, (\cdot)_{T'}] \diam \partial_1^2 v \|_{L^p}^p
	\big\rangle \\
&\lesssim   
\sum_{T' \leq \frac12, \mathrm{dyadic}}({T'}^{\frac14})^{p(2\alpha' - 2\alpha)} ({T'}^{\frac14})^{-3}  .
\end{align*}
Estimate \eqref{yoo1} for $p>\frac{3}{2(\alpha'-\alpha)}$ then follows by summing this geometric series. 
The same bound for smaller $p$ can be derived from the bound for large $p$ and Jensen's inequality.
\medskip

{\sc Step} \arabic{P2St}\label{P2St2}\refstepcounter{P2St}. 
Bounding the supremum over $a_0$, $a_0'$. 
In the following steps we use  the abbreviation
\begin{align}\label{yoo6}
A(\cdot,T,a_0,a_0',\eps_0,\eps_1) 
= 
	\Big(\frac{\partial}{\partial a_0}\Big)^n 
	\Big(\frac{\partial}{\partial a_0'}\Big)^m 
	[ v_{\eps_0}( \cdot, a_0),(\cdot)_T] \diam 
	\{ f_{\eps_1}, \partial_1^2 v_{\eps_1}(\cdot, a_0')  \} .
\end{align}
In this step we show that for  $\eps_0,\eps_1 \in (0,1]$ and $\kappa \ll 1$
\begin{align}
 \label{FinalStoch1A}
 & \Big\langle 
 	\Big( 
		\sup_{a_0, a_0' \in [\lambda,\frac{1}{\lambda}]} 
		\sup_{T \leq 1} (T^{\frac14})^{2-2\alpha} 
		\|A   \| 
	\Big)^p
\Big\rangle^{\frac{1}{p}} 
\lesssim 1,
  \\
 \label{FinalStoch1AA}
 & \Big\langle 
 	\Big( 
		\sup_{a_0, a_0' \in [\lambda,\frac{1}{\lambda}]} 
		\sup_{T \leq 1} (T^{\frac14})^{2-2\alpha} 
		\Big\| \eps_i\frac{\partial}{\partial \eps_i} A   \Big\| 
	\Big)^p
\Big\rangle^{\frac{1}{p}} 
\lesssim \eps_i^{\frac{\kappa}{4}} \qquad \text{for all $p<\infty$}.
\end{align}
For \eqref{FinalStoch1A} we use the Sobolev inequality 
\begin{align}
\notag
\sup_{a_0,a_0' \in [\lambda,\frac{1}{\lambda}]} \|A \|^p %\\
%\notag
&\lesssim   \int_{[\lambda,\frac{1}{\lambda}]} \int_{[\lambda,\frac{1}{\lambda}]} \Big\|\Big \{1, \frac{\partial}{\partial a_0},\frac{\partial}{\partial a_0'}  \Big\}A \Big\|^p da_0 \, da_0'   
%+   \int_{[\lambda,\frac{1}{\lambda}]}\int_{[\lambda,\frac{1}{\lambda}]} \| \frac{\partial}{\partial a_0}A \|^p da_0 \, da_0' \\
%\notag
%&\qquad 
%+  \int_{[\lambda,\frac{1}{\lambda}]}\int_{[\lambda,\frac{1}{\lambda}]} \| \frac{\partial}{\partial a_0'}A \|^p da_0 \, da_0' 
\end{align}
which holds for  $p >2$. Taking the supremum over $T$, then the expectation 
and invoking Fubini's theorem and \eqref{yoo1} yields
\begin{align*}
\notag
&\Big\langle
	\Big(
	 	 \sup_{a_0 ,a_0'\in [\lambda,\frac{1}{\lambda}]} 
%		 \sup_{a_0' \in [\lambda,\frac{1}{\lambda}]}
%		 \sup_{\eps \in (0,1]} 
		 \sup_{T \leq 1} (T^{\frac14})^{2 -2 \alpha }
		  \|  A \| 
	\Big)^p 
\Big\rangle\\
\notag
&\lesssim  
\int_{[\lambda,\frac{1}{\lambda}]}\int_{[\lambda,\frac{1}{\lambda}]} 
	\Big\langle
		\Big(
		 	 \sup_{T \leq 1} (T^{\frac14})^{2 -2 \alpha} 
			 \Big\|   \Big\{1,\frac{\partial}{\partial a_0},\frac{\partial}{\partial a_0'}  \Big\}A \Big\|
		\Big)^p 
	\Big\rangle 
da_0 \, da_0' 
%\\
%\notag
%&\qquad + 
%\int_{[\lambda,\frac{1}{\lambda}]}\int_{[\lambda,\frac{1}{\lambda}]} 
%	\Big\langle 
%		\Big( 
%			\sup_{T \leq 1} (T^{\frac14})^{2- 2 \alpha}  
%			\| \frac{\partial}{\partial a_0}A \|
%		\Big)^p 
%	\Big\rangle
%da_0 \, da_0' \\
%& \qquad +
%\int_{[\lambda,\frac{1}{\lambda}]} \int_{[\lambda,\frac{1}{\lambda}]}  
%	\Big\langle
%		\Big(
%			\sup_{T \leq 1} (T^{\frac14})^{2 -2 \alpha}  
%			\| \frac{\partial}{\partial a_0'} A \|
%		\Big)^p 
%	\Big\rangle 
%da_0 \, da_0'  
\lesssim 1,
\end{align*}
so \eqref{FinalStoch1A} follows. For \eqref{FinalStoch1AA} we repeat the same calculation with $A$ replaced by $\eps_i\frac{\partial}{\partial \eps_i} A$ and \eqref{yoo1} replaced by \eqref{yoo2}.

\medskip

{\sc Step} \arabic{P2St}\label{P2St3}\refstepcounter{P2St}. Bounding the supremum over $\eps_i$.
Let $A$ be defined as in \eqref{yoo6} above. In this step we upgrade \eqref{FinalStoch1A}and \eqref{FinalStoch1AA} to
\begin{align}
 \label{FinalStoch1AAAA}
 & \Big\langle 
 	\Big( 
		\sup_{\eps_0,\eps_1 \in (0,1]}
		\sup_{a_0, a_0' \in [\lambda,\frac{1}{\lambda}]} 
		\sup_{T \leq 1} (T^{\frac14})^{2-2\alpha} 
		\|A   \| 
	\Big)^p
\Big\rangle^{\frac{1}{p}} 
\lesssim 1
\end{align}
valid for $\alpha < \alpha'$.
As in the previous step, we use the Sobolev inequality
\begin{align*}
\sup_{\eps_0,\eps_1 \in (0,1]} |A(\eps)|^p 
\lesssim \int_{[0,1]} \int_{[0,1]} \Big| \Big\{1,\frac{\partial}{\partial \eps_0} ,\frac{\partial}{\partial \eps_1} \Big\} A(\eps) \Big|^p d \eps_0 d \eps_1 ,
%+ \int_0^1  | \frac{\partial}{\partial \eps} A(\eps)| d \eps .
\end{align*}
valid for $p>2$.
We now multiply with $(T^{\frac14})^{2 - \alpha +\kappa}$ for some $\alpha < \alpha'$ and $0<\kappa \ll 1$,  take the supremum over $x$, $T$, $a_0$, $a_0'$ of this estimate
and finally take the expectation to arrive at
\begin{align*}
\notag
&\Big\langle
	\Big(
		\sup_{\eps_0,\eps_1 \in (0,1]} 
	 	\sup_{a_0 ,a_0'\in [\lambda,\frac{1}{\lambda}]} 
%		\sup_{a_0' \in [\lambda,\frac{1}{\lambda}]}
		\sup_{T \leq 1} (T^{\frac14})^{2 -2 \alpha +\kappa }
		\|  A \| 
	\Big)^p 
\Big\rangle^{\frac{1}{p}}\\
\notag
&\lesssim  
\notag
\int_{[0,1]} \int_{[0,1]}
	\Big\langle
		\Big(
	 		 \sup_{a_0 ,a_0'\in [\lambda,\frac{1}{\lambda}]} 
%			 \sup_{a_0' \in [\lambda,\frac{1}{\lambda}]}
			 \sup_{T \leq 1} (T^{\frac14})^{2 -2 \alpha +\kappa }
		  	\Big\| \Big\{1,\frac{\partial}{\partial \eps_0}, \frac{\partial}{\partial \eps_1}   \Big\}  A \Big\| 
		\Big)^p 
	\Big\rangle^{\frac{1}{p}}
d \eps_0 d \eps_1\\
%\notag
%& \qquad +
%\int_0^1
%	\Big\langle
%		\Big(
%	 		 \sup_{a_0 ,a_0'\in [\lambda,\frac{1}{\lambda}]} 
%%			 \sup_{a_0' \in [\lambda,\frac{1}{\lambda}]}
%			 \sup_{T \leq 1} (T^{\frac14})^{2 -2 \alpha +\kappa }
%		  	\|  \frac{\partial}{\partial \eps} A \| 
%		\Big)^p 
%	\Big\rangle^{\frac{1}{p}}
%d \eps \\
&\overset{\eqref{yoo1},\eqref{yoo2}}{\lesssim} 
\int_{[0,1]} \int_{[0,1]} \big\{1, \eps_0^{\frac{\kappa}{4}-1}, \eps_1^{\frac{\kappa}{4}-1} \big\} d \eps_0 d\eps_1 \lesssim 1.
\end{align*}
Now \eqref{FinalStoch1AAAA} follows by relabelling $-2 \alpha + \kappa$ as $-2 \alpha$.

\medskip

{\sc Step} \arabic{P2St}\label{P2St4}\refstepcounter{P2St}. Bounding  $\eps$ differences.
In this step we only consider the diagonal where $\eps_0=\eps_1=\eps$ in $A$ defined 
in \eqref{yoo6} 
and simply write $A(\eps)$ instead of $A(\eps,\eps)$. Note that with this notation
\[
\eps \frac{\partial}{\partial \eps} A(\eps) = \eps_0\frac{\partial}{\partial \eps_0} A(\eps_0,\eps_1)\Big|_{\eps_0=\eps_1=\eps} + \frac{\partial}{\partial \eps_1} A(\eps_0,\eps_1)\Big|_{\eps_0=\eps_1=\eps}.
\]
We claim that for $ \kappa \ll 1$
and all $p< \infty$ and $\alpha < \alpha'$ 
\begin{align}
%\notag
 & \Big\langle 
 	\Big( 
		\sup_{a_0, a_0' \in [\lambda,\frac{1}{\lambda}]} 
		\sup_{\eps_0 \neq \eps_1 \in (0,1]}   
		\sup_{T \leq 1} (T^{\frac14})^{2-2\alpha+\kappa} 
		|\eps_1-\eps_0|^{-\frac{\kappa}{4}}
%		 \\
% & \qquad 
 \|A(\eps_1) -  A(\eps_0)  \| 
	\Big)^p
\Big\rangle^{\frac{1}{p}} 
\lesssim 1.
  \label{FinalStoch1C}
\end{align}
We start the argument with Sobolev's inequality  
\begin{align*}
\sup_{\eps_0 \neq \eps_1 \in (0,1]} \frac{ |A(\eps) -  A(\bar{\eps})| }{|\eps_1 - \eps_0|^{\frac{\kappa}{4}} }
% &= 
% \Big| 
% 	\int_{\eps_1}^{\eps_2}  
%		\frac{\partial}{ \partial \eps} A(\eps) 
%	d \eps 
%\Big| \\
& \leq 
%|\eps_2 - \eps_1|^{\frac{\kappa}{4}} 
\Big(
	\int_{0}^{1} 
		\big| \frac{\partial}{ \partial \eps} A(\eps) \big|^{\frac{1}{1-\frac{\kappa}{4}}} 
	d \eps
\Big)^{1-\frac{\kappa}{4}} .
\end{align*}
Now, we multiply this estimate with $(T^{\frac14})^{2 - 2\alpha +  \kappa + \bar{\kappa}}$ for another $0<\bar{\kappa} \ll 1$, 
take the supremum over $x$, $T$, $a_0$ and $a_0'$, then $p$-th moments, and finally invoke Minkowski's inequality (for $p> \frac{1}{1-\frac{\kappa}{4}}$)
and  \eqref{yoo2}
 to get
\begin{align*}
\notag
 & \Big\langle 
 	\Big( 
		\sup_{a_0, a_0' \in [\lambda,\frac{1}{\lambda}]} 
		\sup_{\eps_0 \neq \eps_1 \in (0,1]}   
		\sup_{T \leq 1}  (T^{\frac14})^{2-2\alpha+\kappa + \bar\kappa} 
		|\eps_1-\eps_0|^{-\frac{\kappa}{4}}
%		 \\
% & \qquad \times 
	 	\|A(\eps_1) -  A(\eps_0)  \| 
	\Big)^p
\Big\rangle^{\frac{1}{p}} \\
& \lesssim 
\Big(
	\int_{0}^1
			\Big\langle 
				\Big( 
					\sup_{a_0, a_0' \in [\lambda,\frac{1}{\lambda}]}
					\sup_{T \leq 1} (T^{\frac14})^{2 -2 \alpha + \kappa + \bar\kappa} 
					\Big\| \frac{\partial}{\partial \eps}A \Big\|^p 
				\Big)
			\Big\rangle^{
					 \frac{1}{p} 
					 \frac{1}{1- \frac{\kappa}{4} }
					 } 
	d \eps
\Big)^{1 - \frac{\kappa}{4}}
	\\
\notag
& \lesssim 
\int_0^1
	\eps^{
		(\frac{\kappa +  \bar{\kappa}}{4}-1 )
		\frac{1}{1 - \frac{ \kappa}{4}}
		} 
d \eps
\lesssim 1,
\end{align*}
so \eqref{FinalStoch1C} follows by relabelling $-2 \alpha + \bar{\kappa}$ as $-2\alpha$.

 \medskip
 {\sc Step} \arabic{P2St}\label{P2St5}\refstepcounter{P2St}. Conclusion. To shorten notation, we only treat 
 the product $v_\eps \diam f_\eps$. Writing 
 \begin{align*}
( v_\eps  \diam f_\eps )_T = v_\eps (f_{\eps})_T - [v_\eps, (\cdot)_T] \diam f_\eps, 
 \end{align*}
and invoking \eqref{SB2} and \eqref{SB2A} for the first and \eqref{FinalStoch1C} for the second term 
imply that  $v_\eps \diam f_\eps$ converges almost surely with respect to the $C^{\alpha -2}$ norm to a limit
$v \diam f$. Furthermore, the estimates \eqref{FinalStoch1AAAA} and \eqref{FinalStoch1C} remain true if the supremum over $\eps \in (0,1]$ is extended to include the limit as $\eps \to 0$.

 \bigskip

%%%%%%%%%%%%%%%%%%%%%%%%%%%%%%%%%%%%%%%%%%
%%%%%%%%%%%%%%%%%%%%%%%%%%%%%%%%%%%%%%%%%%
%%%%%%%%%%%%%%%%%%%%%%%%%%%%%%%%%%%%%%%%%%

\subsection{Proof of Lemma~\ref{lem6}}
\newcounter{L6St} % proofstep = 0
\refstepcounter{L6St} % increases value by 1

{\sc Step} \arabic{L6St}\label{L6St1}\refstepcounter{L6St}. Proof of (i).
By stationarity and \eqref{SB1} we may write
\begin{align*}
&g_1(\eps,a_0)\overset{\eqref{Intro77}}{=}\langle v_\eps(0, a_0) f_\eps(0) \rangle\\
  &= 
 \Big\langle 
 	\int_{[0,1)^2} 
		v_\eps(x, a_0) f_\eps(x) 
	dx 
\Big\rangle 
= 
\sum_{k \in (2\pi \Z)^2} 
	\langle 
		\hat{v}_\eps(k, a_0) \hat{f}_\eps(-k) 
	\rangle\\
&=  
\sum_{k \in (2\pi \Z)^2} 
	\hat{G}( k, a_0) 
	\langle 
		f_\eps(k) f_\eps(-k) 
	\rangle 
=   
\sum_{k \in (2\pi \Z)^2} 
	\hat{G}( k, a_0) 
	\hat{C}(k) |\hat{\ph}_{\eps} (k)|^2,
\end{align*}
where $\hat{G}$ denotes the Fourier transform of the Greens function introduced in \eqref{e:Greens} above.
As the left hand side of this expression is real valued, the imaginary part of the sum of the right hand side also has to vanish.  As $\hat{C}$ is
real valued this means that we can replace $\hat{G}( \cdot, a_0)$ by its real part  (given in \eqref{e:Greens}) thereby yielding \eqref{const1}.

\medskip

The same calculation yields
\begin{align*}
g_2(\eps,a_0,a_0')\overset{\eqref{Intro77}}{=}	\langle  v_\eps (0, a_0)  \partial_1^2 v_\eps(0, a_0') 
\rangle  
%&= 
%\sum_{k \in (2\pi \Z)^2} 
%	\hat{G}( k, a_0) (-k_1^2)\hat{G}(-k,a_0') 
%	\langle 
%		f_\eps(k) f_\eps(-k) 
%	\rangle \\
 =    
 \sum_{k \in (2\pi \Z)^2} 
 	\hat{G}( k, a_0) (-k_1^2 )
	\hat{G}(k,a_0') 
	\hat{C}(k) 
	|\hat{\varphi}_{\eps} (k)|^2,
\end{align*}
and after calculating the real part of $\hat{G}( k, a_0) (-k_1^2)\hat{G}(k,a_0')$ we arrive at \eqref{const2}.

\medskip
{\sc Step} \arabic{L6St}\label{L6St2}\refstepcounter{L6St}. Proof of (ii).
By the condition $a_0 \in [\lambda,\frac{1}{\lambda}]$ we immediately see from \eqref{const1} that 
convergence of $g_1 (\eps,a_0)$ is equivalent to \eqref{A2}.
%\begin{align}\label{A2}
%\sum_{k \in (2 \pi \Z)^2 \setminus \{0\}}  \frac{k_1^2}{ k_1^4 + k_2^2} \hat{C}(k)  < \infty.
%\end{align}
%
Furthermore, given that the ratio of the kernels appearing in \eqref{const1} and \eqref{const2} is bounded 
\begin{align*}
 \Big| \frac{- a_0' k_1^4 + a_0^{-1}k_2^2}{(a_0')^2 k_1^4 + k_2^2}\Big| \leq \lambda^{-3},
\end{align*}
\eqref{A2} also implies the convergence of the $g_2(\eps, a_0,a_0')$ as $\eps$ goes to zero. 
 The condition \eqref{A2} also implies the convergence  for arbitrary derivatives of  $g_1$, $g_2$ with respect to $a_0,a_0'$.
%
%
%The quantity $\partial_1^2 v(\cdot, a_0)$ is obtained from $f$ through a regularity-preserving transformation, as can be seen in terms of the Fourier transform
%%
%\begin{align*}
%\widehat{\partial_1^2 v}(k, a_0) =  \frac{ k_1^2}{a_0k_1^2 - i k_2} \hat{f}(k)
%\end{align*}
%%
%and noting that $ \frac{ k_1^2}{a_0 k_1^2 - i k_2}$ is a bounded symbol (see also Lemma \ref{LA1}). Therefore, the proofs for  $v(\cdot, a_0) \diam f$ and $v(\cdot, a_0) \diam \partial_1^2 v(\cdot, a_0')$ are essentially identical. 
%The list of commutators needed for the deterministic analysis  also includes various derivatives with respect to $a_0$ and $a_0'$, but these derivatives do not change the regularity either. 
%
%For example we have for any $n \geq 1$
%%
%\begin{align}\label{BoundedFourier}
%\frac{\partial^n }{\partial a_0^n}\hat{v}(k, a_0)   = \frac{\partial^n}{\partial a_0^n} \hat{G}( k, a_0)\hat{f}(k) = \frac{(-1)^n n! k_1^{2n}}{(a_0 k_1^2 - i k_2)^n} \hat{v}(k, a_0),
%\end{align}
%%
%and for every $n$ the symbol $\frac{(-1)^n n! k_1^{2n}}{(a_0 k_1^2 - i k_2)^{n}}$ is also bounded.
%
%\medskip
%
%
 For example, recalling the fact that the term $\frac{a_0  k_1^2 }{ a_0^2 k_1^4 +  k_2^2}$ is nothing but the real part $\mathfrak R$ of $\hat{G}(k, a_0)$ we can write
\begin{align*}
\Big(\frac{\partial}{\partial a_0}\Big)^n 
g_1(\eps, a_0)  
&= 
\sum_{k \in (2 \pi \Z)^2 \setminus \{0\}} 
	\mathfrak{R} 
	\Big( 
		\Big(\frac{\partial}{\partial a_0}\Big)^n 
		\hat{G}(k, a_0)
	\Big)  
	\hat{C}(k)|\hat{\ph}_{\eps}(k)|^2\\
&\overset{\eqref{BoundedFourier}}{=} 
\sum_{k \in (2 \pi \Z)^2 \setminus \{0\}} 
	\mathfrak{R}
	\Big( 
		\frac{(-1)^n n! k_1^{2n}}{(a_0 k_1^2 - i k_2)^n} 
		\hat{G}(k, a_0)
	\Big) 
	\hat{C}(k) |\hat{\ph}_{\eps}(k)|^2.
\end{align*}
Given that for any $n \geq 1$ the modulus of the quantity under the real part $\mathfrak R$ is $\lesssim \frac{  k_1^2 }{  k_1^4 +  k_2^2}$ 
 the convergence  as $\eps \to 0$ under \eqref{A2} follows. A similar argument works for $g_2$.

\bigskip
%%%%%%%%%%%%%%%%%%%%%%%%%%%%%%%%%%%%%%%%%%%%
%%%%%%%%%%%%%%%%%%%%%%%%%%%%%%%%%%%%%%%%%%%%
%%%%%%%%%%%%%%%%%%%%%%%%%%%%%%%%%%%%%%%%%%%%

%%%%%%%%%%%%%%%%%%%%%%%%%%%%%%%%%%%%%%%%%%%%%%%%%%%%%%%%%%%%%%%%%%%%%%%%%%%%%%%%%%%%%%%%%%%
%%%%%%%%%%%%%%%%%%%%%%%%%%%%%%%%%%%%%%%%%%%%%%%%%%%%%%%%%%%%%%%%%%%%%%%%%%%%%%%%%%%%%%%%%%%
%%%%%%%%%%%%%%%%%%%%%%%%%%%%%%%%%%%%%%%%%%%%%%%%%%%%%%%%%%%%%%%%%%%%%%%%%%%%%%%%%%%%%%%%%%%
%%%%%%%%%%%%%%%%%%%%%%%%%%%%%%%%%%%%%%%%%%%%%%%%%%%%%%%%%%%%%%%%%%%%%%%%%%%%%%%%%%%%%%%%%%%
%%%%%%%%%%%%%%%%%%%%%%%%%%%%%%%%%%%%%%%%%%%%%%%%%%%%%%%%%%%%%%%%%%%%%%%%%%%%%%%%%%%%%%%%%%%

\section{Proofs of Theorems \ref{Theo:Introduction1} and \ref{Theo:Introduction2}}

According to Lemma~\ref{lem5}  under assumption~\eqref{A1}, we have
\[
\Big\langle  \sup_{\eps \in [0,1] }\| f\|_{\alpha-2}^p \Big\rangle^{\frac{1}{p}} <\infty
\]
for any $\alpha< \alpha'$ and $p<\infty$ and we have almost surely and in every stochastic 
$L^p$ space that $\| f - f_\eps \|_{\alpha-2} \to 0$.
Under the same assumption according to Proposition~\ref{P2} the renormalized products
$ v_\eps( \cdot, a_0) \diam f_\eps$ and $ v_\eps( \cdot, a_0) \diam \partial_1^2 v_\eps(\cdot, a_0')$
defined in \eqref{e:RenProducts}
%
%
%\begin{align*}
% v_\eps( \cdot, a_0) \diam f_\eps :=& v_\eps( \cdot, a_0) f_\eps - \langle v_\eps( \cdot, a_0) f_\eps \rangle, \\
% v_\eps( \cdot, a_0) \diam \partial_1^2 v_\eps(\cdot, a_0')   :=& v_\eps( \cdot, a_0) \partial_1^2 v_\eps(\cdot, a_0')  -\langle v_\eps( \cdot, a_0) \partial_1^2 v_\eps(\cdot, a_0')  \rangle  
%\end{align*}  
converge to limits denoted by $v\diam f$ and $v \diam \partial_1^2 v$  as $\eps$ goes to zero in the sense that  almost surely  the quantities
\begin{align*}	
		\big\|   
				[v_\eps,(\cdot) ] 					
				\diam  f_\eps 
				-[v,(\cdot) ] 
				\diam  f 
			\big\|_{2\alpha-2,2} \;, \quad
			\big\|			
				[v_\eps,(\cdot) ] 					
				\diam \partial_1^2 v_\eps    
				-[v,(\cdot) ] 
				\diam \partial_1^2 v     
			\big\|_{2\alpha-2,2,2} 
\end{align*}
converge to zero. Furthermore, we have the moment bounds
\begin{align}
%\notag
  \Big\langle
% 	 \Big(
%	 	 \sup_{a_0, a_0' \in [\lambda,\frac{1}{\lambda}]} 
		 \sup_{\eps_0,\eps_1 \in [0,1]}    
%		 \sup_{T \leq 1}  (T^{\frac14})^{2-2\alpha'} 
		 	\| 
%				\frac{\partial^n}{\partial a_0^n} 
%				\frac{\partial^{m}}{\partial (a_0')^{m}} 
				[ v_{\eps_0},(\cdot) ]
  				\diam f_{\eps_1}
			\|_{2\alpha-2,2} ^p
%	\Big)^p 
\Big\rangle^{\frac{1}{p}} , \quad
  \Big\langle
% 	 \Big(
%	 	 \sup_{a_0, a_0' \in [\lambda,\frac{1}{\lambda}]} 
		 \sup_{\eps_0,\eps_1 \in [0,1]}  
%		 \sup_{T \leq 1}  (T^{\frac14})^{2-2\alpha'} 
		 	\| 
%				\frac{\partial^n}{\partial a_0^n} 
%				\frac{\partial^{m}}{\partial (a_0')^{m}} 
				[ v_{\eps_0},(\cdot) ]
 \label{FinalStoch1bis}
  				\diam \partial_1^2 v_{\eps_1}  
			\|_{2\alpha-2,2,2} ^p
%	\Big)^p 
\Big\rangle^{\frac{1}{p}} 
<\infty ,
\end{align}
for all $p<\infty$.

\medskip
Let  $N_0 \ll 1$  be so small that Theorem~\ref{Theo} holds and set
\begin{align*}
 \eta_0^{-1} = \frac{1}{N_0}
\sup_{\eps,\eps_0,\eps_1 \in [0,1]}  \max\Big\{ 
	\| f_\eps \|_{\alpha-2},
	\|[v_{\eps_0},(\cdot)] \diam f_{\eps_1} \|_{2\alpha-2,2,2}^{\frac12},
	\| [v_{\eps_0},(\cdot)] \diam \partial_1^2 v_{\eps_1}\|_{2\alpha-2,2,2}^{\frac12}
	\Big\}.
\end{align*}
Then the moment bound \eqref{eta-moment-bound} holds, 
and for all $\eta\leq \eta_0$ the functions/distributions  $\eta f_\eps$, $\eta^2 [v_\eps, (\cdot)] f_\eps$ and $\eta^2 [v_\eps, (\cdot)] \partial_1^2v_\eps$  satisfy the smallness condition \eqref{x92}, \eqref{x93}, and \eqref{v14} uniformly in $\eps\in [0,1]$.  Thus Theorem~\ref{Theo} part (i) yields the existence and uniqueness of  a solution $u$ to \eqref{v95A}, \eqref{v94A}, \eqref{v96bisA}, as well as solutions $u_\eps$ to 
the corresponding regularized problems with $\delta $ in Theorem~\ref{Theo:Introduction1} being the implicit constant in \eqref{v96bis}.
 By Corollary~\ref{Coro}  the regularized problem takes the form of \eqref{glatteGleichung}
and part (ii) of Theorem~\ref{Theo}, more precisely estimate \eqref{wk54}, yields the convergence to zero of $\| u - u_\eps \|+ [u-u_\eps]_{\alpha}$.
%which according to takes the form of .
% Finally, it remains to argue that 
%in the regularized case  \eqref{v95A},\eqref{v94A}, \eqref{v96bisA} is equivalent to  and this is the content of .

%%%%%%%%%%%%%%%%%%%%%%%%%%%%%%%%%%%%%%%%%%%%%%%%%%%%%%%%%%%%%%%%%%%%%%%%%%%%%%%%%%%%%%%%%%%%%%%%
%%%%%%%%%%%%%%%%%%%%%%%%%%%%%%%%%%%%%%%%%%%%%%%%%%%%%%%%%%%%%%%%%%%%%%%%%%%%%%%%%%%%%%%%%%%%%%%%
%%%%%%%%%%%%%%%%%%%%%%%%%%%%%%%%%%%%%%%%%%%%%%%%%%%%%%%%%%%%%%%%%%%%%%%%%%%%%%%%%%%%%%%%%%%%%%%%
%%%%%%%%%%%%%%%%%%%%%%%%%%%%%%%%%%%%%%%%%%%%%%%%%%%%%%%%%%%%%%%%%%%%%%%%%%%%%%%%%%%%%%%%%%%%%%%%
%%%%%%%%%%%%%%%%%%%%%%%%%%%%%%%%%%%%%%%%%%%%%%%%%%%%%%%%%%%%%%%%%%%%%%%%%%%%%%%%%%%%%%%%%%%%%%%%
%%%%%%%%%%%%%%%%%%%%%%%%%%%%%%%%%%%%%%%%%%%%%%%%%%%%%%%%%%%%%%%%%%%%%%%%%%%%%%%%%%%%%%%%%%%%%%%%
%%%%%%%%%%%%%%%%%%%%%%%%%%%%%%%%%%%%%%%%%%%%%%%%%%%%%%%%%%%%%%%%%%%%%%%%%%%%%%%%%%%%%%%%%%%%%%%%
%%%%%%%%%%%%%%%%%%%%%%%%%%%%%%%%%%%%%%%%%%%%%%%%%%%%%%%%%%%%%%%%%%%%%%%%%%%%%%%%%%%%%%%%%%%%%%%%
%%%%%%%%%%%%%%%%%%%%%%%%%%%%%%%%%%%%%%%%%%%%%%%%%%%%%%%%%%%%%%%%%%%%%%%%%%%%%%%%%%%%%%%%%%%%%%%%
%%%%%%%%%%%%%%%%%%%%%%%%%%%%%%%%%%%%%%%%%%%%%%%%%%%%%%%%%%%%%%%%%%%%%%%%%%%%%%%%%%%%%%%%%%%%%%%%
%%%%%%%%%%%%%%%%%%%%%%%%%%%%%%%%%%%%%%%%%%%%%%%%%%%%%%%%%%%%%%%%%%%%%%%%%%%%%%%%%%%%%%%%%%%%%%%%
%%%%%%%%%%%%%%%%%%%%%%%%%%%%%%%%%%%%%%%%%%%%%%%%%%%%%%%%%%%%%%%%%%%%%%%%%%%%%%%%%%%%%%%%%%%%%%%%

\appendix

\section{Some additional lemmas}\label{s:Appendix}

\addtocontents{toc}{\protect\setcounter{tocdepth}{1}}

\begin{lemma}\label{LA1}
The (mean-free) solution of (\ref{v53}) satisfies the estimate
\begin{align}\label{w31}
\sup_{a_0}[v(\cdot,a_0)]_\alpha\lesssim \|f\|_{\alpha-2}.
\end{align}
\end{lemma}

\subsection*{Proof of Lemma \ref{LA1}}
%{\sc Proof of Lemma \ref{LA1}}.

All functions are space-time period if not stated otherwise.

\newcounter{LA1S} % proofstep = 0
\refstepcounter{LA1S} % increases value by 1

{\sc Step} \arabic{LA1S}.\label{LA1S1}\refstepcounter{LA1S}
Reduction. We claim that it is enough to show
\begin{align}\label{w30}
 \sup_{T\le 1}(T^\frac{1}{4})^{2-\alpha}\|f_T\| 
\sim\inf\Big\{[f_1]_\alpha+[f_2]_\alpha+|c|\,\Big|\,f=\partial_1^2f_1+\partial_2f_2+c\;\Big\},
\end{align}
where the infimum is over all triplets $(f_1,f_2,c)$ of two functions and a constant.
Incidentally, the equivalence confirms that the left hand side indeed
defines the (parabolic) $C^{\alpha-2}$-norm. Let the decomposition
$f=\partial_1^2f_1+\partial_2f_2+c$ be near-optimal in the right hand side of (\ref{w30}),
that is,
\begin{align}\label{w37}
[f_1]_\alpha+[f_2]_\alpha\leq 2\sup_{T\leq 1}\|f_T\|.
\end{align}
By uniqueness of the mean-free
solution of (\ref{v53}) this induces $v(\cdot,a_0)$ $=\partial_1^2v_1$ $+\partial_2v_2$ where $v_i$, $i=1,2$,
denote the mean-free solutions of $(\partial_2-a_0\partial_1^2)v_i=f_i$.
%%
%\begin{align*}
%(\partial_2-a_0\partial_1^2)v_1=f_1\quad\mbox{and}\quad
%(\partial_2-a_0\partial_1^2)v_2=f_2.
%\end{align*}
%%
By classical $C^{\alpha+2}$-Schauder theory \cite[Theorem 8.6.1]{Krylov} we have 
$[\partial_1^2v_i]_\alpha$ $+[\partial_2v_i]_\alpha$ $\lesssim[f_i]_\alpha$,
so that (\ref{w31}) follows from (\ref{w37}).

\medskip

{\sc Step} \arabic{LA1S}.\label{LA1S2}\refstepcounter{LA1S}
For the solution of
\begin{align}\label{w32}
(\partial_2-\partial_1^2)v=Pf
\end{align}
we claim
\begin{align}\label{w33}
\|v_T-v\|\lesssim N_0\max\{(T^\frac{1}{4})^{\alpha},(T^\frac{1}{4})^{2}\}\quad\mbox{for all}\;T>0,
\end{align}
where we have set for abbreviation 
\begin{equation}\label{pppi}
N_0:=\sup_{T\le 1}(T^\frac{1}{4})^{2-\alpha}\|f_T\|.
\end{equation}
We start by noting that the definition of $N_0$ may be extended to
the control of $T\ge 1$ by the semi-group property (\ref{1.10}) in form of $f_T=(f_1)_{T-1}$
and (\ref{1.13}) in form of $\|f_T\|\lesssim\|f_1\|$. We thus have
\begin{align}\label{w34}
\|f_T\|\lesssim N_0\max\{T^{\alpha-2},1\}.
\end{align}
By approximation through (standard) convolution, which preserves
(\ref{w32}) and does increase $N_0$, we may assume that $f$ and $v$ are smooth.
By definition of the convolution $(\cdot)_t$ we have
\begin{align*}
\partial_tv_t&=-(\partial_1^4-\partial_2^2)v_t=(-\partial_1^2-\partial_2)(\partial_2-\partial_1^2)v_t
\stackrel{(\ref{w32})}{=}(-\partial_1^2-\partial_2)Pf_t\nonumber\\
&\stackrel{(\ref{1.10})}{=}(-\partial_1^2-\partial_2)(f_{\frac{t}{2}})_{\frac{t}{2}}.
\end{align*}
Hence we obtain by (\ref{1.13}) for all $T\le 1$
\begin{align*}
\|\partial_tv_t\|\lesssim(t^\frac{1}{4})^{-2}\|f_{\frac{t}{2}}\|
\stackrel{(\ref{w34})}{\lesssim} N_0\max\{(t^\frac{1}{4})^{\alpha-4},(t^\frac{1}{4})^{-2}\}.
\end{align*}
Integrating over $t\in(0,T)$ we obtain (\ref{w33}) by the triangle inequality. 

\medskip

{\sc Step} \arabic{LA1S}.\label{LA1S3}\refstepcounter{LA1S}
For $v$ defined through (\ref{w32}) we have
\begin{align}\label{w36}
[v]_\alpha\lesssim N_0,
\end{align}
where $N_0$ is as in \eqref{pppi}.
As in Step \ref{LA1S2} we may assume that $f$ and $v$ are smooth so that $[v]_\alpha$ is finite.
Because of periodicity, it is sufficient to probe H\"older continuity
for pairs $(x,y)$ of points with $d(y,x)\le 4$.
For any $T>0$ we have the identity
\begin{align*}
\lefteqn{v(y)-v(x)=(v_T-v)(y)-(v_T-v)(x)}\nonumber\\
&-\int_0^1\partial_1v_T(sy+(1-s)x)(y-x)_1+\partial_2v_T(sy+(1-s)x)(y-x)_2ds,
\end{align*}
from which we obtain the inequality
\begin{align*}
|v(y)-v(x)|\le 2\|v_T-v\|+\|\partial_1v_T\|d(y,x)+\|\partial_2v_T\|d^2(y,x).
\end{align*}
From Step \ref{LA1S2} and (\ref{1.13}) we obtain the estimate
\begin{align*}
\lefteqn{|v(y)-v(x)|}\nonumber\\
&\lesssim N_0\max\{(T^\frac{1}{4})^\alpha,(T^\frac{1}{4})^2\}
+[v]_\alpha\big((T^\frac{1}{4})^{\alpha-1}d(y,x)+(T^\frac{1}{4})^{\alpha-2}d^2(y,x)\big).
\end{align*}
With the ansatz $T^\frac{1}{4}=\frac{1}{\epsilon}d(y,x)$ for some $\epsilon\le 1$
and making use of $d(y,x)\leq 1$ we obtain
\begin{align*}
|v(y)-v(x)|\lesssim \big(\epsilon^{-2}N_0+[v]_\alpha(\epsilon^{1-\alpha}+\epsilon^{2-\alpha})\big)d^\alpha(y,x).
\end{align*}
Fixing an $\epsilon$ sufficiently small to absorb the last right-hand-side term into the left hand side
we infer (\ref{w36}).

\medskip

{\sc Step} \arabic{LA1S}.\label{LA1S4}\refstepcounter{LA1S}
We finally establish the equivalence of norms (\ref{w30}). The direction
$\lesssim$ follows immediately from (\ref{1.13}). The direction $\gtrsim$ follows
from Step \ref{LA1S3} with $f_1=v$, $f_2=-v$, and $c=\int_{[0,1)^2}f$.

\bigskip

\begin{lemma}\label{LA2}
\begin{align}\label{LA2.2}
\|[x_1,(\cdot)]f\|_{\alpha-1}
\lesssim \|f\|_{\alpha-2}.
\end{align}
\end{lemma}

\subsection*{Proof of Lemma \ref{LA2}}
%{\sc Proof of Lemma \ref{LA2}}.

Introducing the kernel $\tilde\psi_T(x):=x_1\psi_T(x)$ we start by claiming the representation
\begin{align}\label{LA2.1}
[x_1,(\cdot)_T]f=2\tilde\psi_{\frac{T}{2}}*f_{\frac{T}{2}}.
\end{align}
Indeed, by definition of the commutator and $\tilde\psi_T$ we have
$[x_1,(\cdot)_T]f=\tilde\psi_T*f$, so that the above representation follows
from the formula 
\begin{align}\label{LA2.3}
\tilde\psi_T=2\tilde\psi_{\frac{T}{2}}*\psi_{\frac{T}{2}}.
\end{align}
The argument for (\ref{LA2.3}) relies on the fact that convolution 
is commutative in form of $\tilde\psi_{\frac{T}{2}}*\psi_{\frac{T}{2}}=
\psi_{\frac{T}{2}}*\tilde\psi_{\frac{T}{2}}$, which spelled out means
$\int dy(x_1-y_1)\psi_{\frac{T}{2}}(x-y)\psi_{\frac{T}{2}}(y)$
$=\int dy\psi_{\frac{T}{2}}(x-y)y_1\psi_{\frac{T}{2}}(y)$, and thus implies
$2\int dy(x_1-y_1)\psi_{\frac{T}{2}}(x-y)\psi_{\frac{T}{2}}(y)$
$=x_1\int dy\psi_{\frac{T}{2}}(x-y)\psi_{\frac{T}{2}}(y)$, that is
$2(\tilde\psi_{\frac{T}{2}}*\psi_{\frac{T}{2}})(x)$
$=x_1(\psi_{\frac{T}{2}}*\psi_{\frac{T}{2}})(x)$.
Together with the semi-group property (\ref{1.10}) in form of $\psi_{\frac{T}{2}}*\psi_{\frac{T}{2}}=\psi_T$
this yields (\ref{LA2.3}).

\medskip

From the representation (\ref{LA2.1}) we obtain the estimate
\begin{align*}
\|[x_1,(\cdot)_T]f\|\le2\int dx|x_1 \psi_{\frac{T}{2}}(x)|\|f_{\frac{T}{2}}\|
\stackrel{(\ref{1.13})}{\lesssim}T^\frac{1}{4}\|f_{\frac{T}{2}}\|,
\end{align*}
which yields the desired (\ref{LA2.2}).

\bigskip
The following lemma shows that the definitions \eqref{c.1}, \eqref{e:2alphadef} and \eqref{e:difference3alpha} are independent of the choice of convolution kernel.
\begin{lemma}\label{L:AMMM1}
Let  $\psi$ and $\psi'$ be  Schwartz functions over $\R^2$ with $\int \psi = \int \psi' = 1$. For $T>0$ define
\begin{align}\label{AMM1}
\psi_T(x_1,x_2) = T^{-\frac{3}{4}} \psi\Big(\frac{x_1}{T^{\frac14}}, \frac{x_2}{T^{\frac12}} \Big), \qquad  \psi'_T(x_1,x_2) = T^{-\frac{3}{4}} \psi'\Big(\frac{x_1}{T^{\frac14}}, \frac{x_2}{T^{\frac12}} \Big).
\end{align}
and for an arbitrary Schwartz distribution $f \in \mathcal{S}'(\R^2)$ set 
\begin{align}\label{AMM2}
(f)_T = f \ast \psi_T \qquad \text{and} \qquad  (f)'_T = f \ast \psi'_T.
\end{align}
i) For any $\gamma<0$ we have 
\begin{align}\label{AMM2a}
\sup_{T \leq 1}  (T^{\frac14})^{-\gamma}  \| (f)_T \|\lesssim \sup_{T \leq 1} (T^{\frac14})^{-\gamma}\| (f)'_T \| ,
\end{align}
where $\lesssim$ only refers to $\psi$, $\psi'$ and $\gamma$.
\medskip

ii) Let  $\alpha>0$ and $\gamma<0$. Let $u$ be a function of class $C^{\alpha}$
and $f$ a distribution of class $C^{\gamma}$. Furthermore, let $u \diam f$ be 
an arbitrary distribution of class $C^{\gamma}$ and define the generalized commutators
$[u,( \cdot)_T ] \diam f := u (f)_T - (u \diam f)_T $ and $[u,( \cdot)_T' ] \diam f := u (f)'_T - (u \diam f)'_T $. Then for $\bar \gamma = \gamma + \alpha$ we have
\begin{align}
\notag
 \sup_{T \leq 1}(T^{\frac14})^{-\bar \gamma} \| [ u, (\cdot)_T ]\diam f \|  
 &\lesssim 
 \sup_{T\leq 1} (T^{\frac14})^{-\bar \gamma}  \| [ u, (\cdot)'_T ] \diam f \|\\
 \label{AMM309}
 & \qquad  
 + [u]_\alpha \sup_{T \leq 1}(T^{\frac14})^{-\gamma} \| (f)'_T \| ,
\end{align}
where $\lesssim$ depends on $\alpha$, $\gamma$ as well as $\psi$ and $\psi'$.
\end{lemma}

\subsection*{Proof of Lemma \ref{L:AMMM1}}
%{\sc Proof of Lemma \ref{L:AMMM1}}.

%
%\begin{proof}[Proof of Lemma~\ref{L:AMMM1}]

\newcounter{ALemmS} % proofstep = 0
\refstepcounter{ALemmS} % increases value by 1

{\sc Step} \arabic{ALemmS}.\label{Representation}  \refstepcounter{ALemmS}
The proof relies on a variant of a construction from \cite{GloriaOtto} which we recall in this step. For the reader's convenience we give self-contained
 proofs of the identities  in Step~\ref{ReprentationProof} below.
First of all,  for any $p>0$ there exists a Schwartz function $\omega^0$ such that 
$\varphi' = \omega^0 \ast \psi'$ satisfies
\begin{align}\label{AMM3}
\int x^{n} \varphi'(x) dx = 
\begin{cases}
1 \qquad \text{for } \alpha =0\\
0 \qquad \text{for } 0< \|n\|_{\text{par}} < p,
\end{cases}
\end{align}
where for $n = (n_1, n_2)$ and $x=(x_1,x_2)$ we write $x^n = x_1^{n_1} x_2^{n_2}$ and use the parabolic norm $ \|n\|_{\text{par}}  = |n_1| + 2| n_2|$.  Furthermore, it is shown that for any $p$ and any $\varphi'$ satisfying \eqref{AMM3} as well as $\theta \ll 1$ (depending on $\varphi, \psi,p$), the function $\psi$ can be represented as
\begin{align}\label{AMM4}
\psi = \sum_{k=0}^\infty \omega^{(k)} \ast \varphi'_{\theta^k},
\end{align}
where $\varphi'_{\theta^k}$ is the rescaled version of $\varphi'$ defined as in \eqref{AMM1} for $T = \theta^k$, and the $\omega^{(k)}$ are Schwartz functions satisfying
\begin{align}\label{AMM5}
\int |\omega^{(k)}|  \lesssim (C_0 \theta^{\frac{p}{4}})^k,
\end{align}
where $C_0 = C_0( \varphi', \psi,p)$.
The convergence of the sum in \eqref{AMM4} holds  in $L^1(\R^2)$ . 
Additionally, we will make use of the bounds
\begin{align}\label{newbound111}
 \int d^{\alpha}(0,x) |\omega^{(k)}  (x) | dx \lesssim   (C_0 \theta^{\frac{p}{4}})^k.
\end{align}

\medskip

We summarize this as
$
\psi = \sum_{k=0}^\infty \omega^{(k)} \ast \omega^0_{\theta^k} \ast \psi'_{\theta^k},
$
which can be rescaled as
\begin{align}\label{AMM6}
\psi_{T} = \sum_{k=0}^\infty \omega^{(k)}_{ T} \ast \omega_{\theta^k  T}^0 \ast \psi'_{\theta^k T},
\end{align}
where as before the index $T$ expresses that a function is rescaled by $T$ as in \eqref{AMM1}. 
%Note for further use below that this scaling 
%preserves the $L^1$ norm such that by  \eqref{AMM5}
%\begin{align}\label{AMM7}
% \int |\omega^{(k)}_{T}| =  \int |\omega^{(k)} |  \lesssim  (C_0 \theta^{\frac{p}{4}})^k
%\end{align}
% with constant independent of $T$.

\medskip

{\sc Step} \arabic{ALemmS}.\label{ALemmSProofi}  \refstepcounter{ALemmS}
Equipped with these results  we now proceed to prove \eqref{AMM2a}. 
%
%First, we claim that we can assume wlog that $f$ is smooth. Indeed, using the 
%lower semicontinuity of $\|  \cdot \|$ with respect to removing convolutions in the first step and 
% assuming that \eqref{AMM2} holds for the smooth functions $f_t$ in the second step, we can write
%\begin{align}
%\notag
% \| (f)_T \| &\leq \limsup_{t \to 0} \| (f_t)_{T} \| \lesssim \limsup_{t \to 0} ((T+t)^{\frac14})^{\gamma}  \sup_{\bar{T} \leq 1} (\bar{ T}^{\frac14})^{-\gamma}\| (f)'_{\bar{T}} \| \\  
%\label{poi0}
% &\leq (T^{\frac14})^{\gamma}\sup_{\bar{T} \leq 1} (\bar{ T}^{\frac14})^{-\gamma}\| (f)'_{\bar{T}} \|  ,
%\end{align}
%which implies \eqref{AMM2} also in the general case.
  Set $N_0:=\sup_{T\leq 1} (T^{\frac14})^{-\gamma}  \|(f)'_T\| $
and write
\begin{align*}
\| (f)_{ T} \| & \overset{\eqref{AMM6}}{=}  \| \sum_{k=0}^{\infty}  (\omega^{(k)}_{T} \ast \omega_{\theta^k T}^0 ) \ast (f)'_{\theta^k T } \|  
%&= \|  \sum_{k=0}^{\infty} (f)'_{\theta^k T } \ast ( \omega^{(k)}_{T} \ast \omega_{\theta^k T}^0 ) \| 
 \leq  \sum_{k=0}^{\infty} \int  | \omega^{(k)}_{T} | \int|  \omega_{\theta^k T}^0 |    \| (f)'_{\theta^k T } \|  \\
& \overset{\eqref{AMM5}}{\lesssim} N_0 \sum_{k=0}^{\infty} (\theta ^\frac{k}{4} T^{\frac14} )^{ \gamma} (C_0 \theta^{\frac{p}{4}})^k \int |\omega^0 |.
\end{align*} 
%where in the last line, we have used~ (we have also used the smoothness of $f$ in the second identity to ensure that we can exchange the convolution and the sum which converges in $L^1$). 
Then \eqref{AMM2a} follows by choosing first $p>|\gamma|$ and then $\theta^{\frac14} \leq \frac{1}{2C_0}$  and then summing the geometric series
over $k$.

\medskip

{\sc Step} \arabic{ALemmS}.\label{ALemmSProofii}  \refstepcounter{ALemmS}
%
%As in the previous step we assume wlog that $f$, $u \diam f$ are smooth (however, 
%we do not assume that $u \diam f = u f$). Indeed, else consider $[u, ( \cdot)_{T+t} ] \diam f = 
%u (f_t)_T - ((u \diam f)_t)_T$ and argue as in \eqref{poi0}.
%
We set   $ N_1:=$  $\sup_{T \leq 1}$ $ (T^{\frac14})^{-\bar \gamma}   \| [ u, (\cdot)'_{T} ] \diam f \| $  and   $N_0':=\sup_{T\leq 1} (T^{\frac14})^{-\gamma}  \|(f)'_T\| $ as before. 
Again, we make use of the representation \eqref{AMM6} of $\psi_{T}$ to write
\begin{align*}
[u, ( \cdot)_{T}] \diam f =   \sum_{k=0}^{\infty} [u,  \omega^{(k)}_{T} \ast \omega_{\theta^k T}^0  \ast   \psi'_{\theta^k T} \ast ] \diam f.
\end{align*}
We apply the commutator relation $[A, BC] = [A,B] C + B [A,C]$ twice,  to 
rewrite each term in this sum as
\begin{align}
\notag
&[u,  \omega^{(k)}_{T} \ast \omega_{\theta^k T}^0  \ast   \varphi_{\theta^k T} \ast ] \diam  f \\
\notag
&= [u,  \omega^{(k)}_{T} \ast \omega_{\theta^k T}^0 \ast  ] (f)'_{\theta^k T} + \omega^{(k)}_{T} \ast \omega_{\theta^k T}^0  \ast [u, ( \cdot)'_{ \theta^k T }] \diam f \\
\notag
& =  [u,   \omega^{(k)}_{T} \ast   ] (\omega_{\theta^k T}^0  \ast (f)'_{\theta^k T}) +  \omega^{(k)}_{T} \ast ( [u,  \omega_{\theta^k T}^0 \ast   ] (f)'_{\theta^k T} )\\
& \qquad   +  \omega^{(k)}_{T} \ast \omega_{\theta^k T}^0  \ast [u, ( \cdot)'_{ \theta^k T }] \diam f  .\label{AMM21}
\end{align}
Note that only the last commutator on the rhs requires the definition of $u \diam f$ and 
all the other commutators are defined classically.
We bound the terms on the right hand side of \eqref{AMM21} one by one, starting with the last.
This expression can be directly bounded 
\begin{align*}
 \| (\omega^{(k)}_{T} \ast \omega_{\theta^k T}^0 ) \ast [u, ( \cdot)_{\theta^k T}']  \diam f \| 
% &\leq  \int |\omega^{(k)}_{T}| \, \int  |\omega_{\theta^k T}^0 |  \, (\theta^{\frac{k}4} T^{\frac14} )^{\bar \gamma} N_1\\
 &= \int |\omega^{(k)}| \, \int  |\omega^0 |  \, (\theta^{\frac{k}4} T^{\frac14} )^{\bar \gamma} N_1.
 \end{align*}
Therefore, the sum in $k$ over this term is controlled by  invoking \eqref{AMM5} for  $p$ large enough, then choosing $\theta$ small enough,  resulting with a geometric series as in Step~\ref{ALemmSProofi}.
\medskip

By Young's inequality, the second term on the right hand side of \eqref{AMM21} is bounded:
\begin{align*}
  \| \omega^{(k)}_{T} \ast ( [u,  \omega_{\theta^k T}^0 \ast   ] (f)'_{\theta^k T}) \|  \leq \int |\omega^{(k)} | \,\, \|  [u,  \omega_{\theta^k T}^0 \ast   ] (f)'_{\theta^k T}\| .
\end{align*}
According to  \eqref{AMM5} the first factor on the rhs is bounded by $ \lesssim (C_0 \theta^{\frac{p}{4}})^k$, while the second factor can be bounded as 
\begin{align*}
 &\|  [u,  \omega_{\theta^k T}^0 \ast  ] (f)'_{\theta^k T}\|\\
 & = \sup_x  \Big| \int (u(x) - u(y)) \omega_{\theta^k T}^0 (y-x)  (f)'_{\theta^k T}(y) dy \Big|\\
 &\leq [u]_{\alpha} N_0 (\theta^{\frac{k}{4}}T^{\frac14} )^{\gamma} \sup_x\int d^\alpha(x,y) |\omega_{\theta^k T}^0  (y-x) |dy \\
 &  = [u]_{\alpha} N_0 (\theta^{\frac{k}{4}}T^{\frac14})^{\bar{\gamma}}   %(\theta^{\frac{k}{4}} T^{\frac14})^{\alpha}  
 \int d^{\alpha}(0,z) |\omega^{0} (z)| dz,
\end{align*}
so that summing these terms over $k$ also yields the required bound as above. 

\medskip
It remains to bound the first term on the rhs of \eqref{AMM21} and for this we write
\begin{align*}
&\|  [u,   \omega^{(k)}_{T} \ast   ] ( \omega_{\theta^k T}^0 \ast  (f)'_{\theta^k T})  \| \\
&\leq \sup_{x} \int | u(y) - u(x)|\, | \omega^{(k)}_{T} (y-x) |dy \, \Big(  \int | \omega_{\theta^k T}^0 |  \Big)\| (f)'_{\theta^k T} \|    \\
%
%& \leq  [u]_{\alpha} \sup_{x} \int  d^\alpha(x,y) |\omega^{(k)}_{T} (y-x) | dy \, \Big(  \int | \omega^0 |  \Big)   N_0    ( \theta^{\frac{k}{4}}   T^{\frac14})^{\gamma}\\
&\leq [u]_{\alpha}    (T^{\frac14})^{\alpha} N_0  ( \theta^{\frac{k}{4}} T^{\frac14})^{\gamma}  \int  d^\alpha (0,z) |\omega^{(k)}  (z) | dz  \, \Big(  \int | \omega^0 |  \Big) .
\end{align*}
The first integral on the rhs is bounded $\lesssim (C_0 \theta^{\frac{p}{4}})^k$ in \eqref{newbound111}, so that finally \eqref{AMM309}  follows once more by choosing $p$ large enough and $\theta$ small enough and summing over $k$.

\medskip

{\sc Step} \arabic{ALemmS}.\label{ReprentationProof}  \refstepcounter{ALemmS}
It remains to give the argument for \eqref{AMM3}, \eqref{AMM4} and \eqref{newbound111} following \cite{GloriaOtto}. The construction of $\omega^0$ is based on the identity
\begin{align*}
A_{ n, m}&:= \int x^{n} \partial^m \psi'(x) dx \\
&= 
\left\{
\begin{array}{ll}
0 \qquad &\text{if } \|n \|_{\text{par}} \leq \|m \|_{\text{par}}, \, n \neq m\\
(-1)^{|m_1|+|m_2|} m_1 ! m_2! \qquad &\text{if } n = m
\end{array}
\right\}.
\end{align*}
This trigonal structure implies  that for any fixed $p$ the linear map $$(a_m)_{\|m\|_{\text{par}}<p} \mapsto (\sum_{ \|m \|_{\text{par}}<p}  A_{ n,m} a_m)_{\|n\|_{\text{par}}<p}$$ is invertible. Furthermore, for each 
$n, m $ the numbers $A_{ n,m}^r:= \int x^{n} \partial^m (\psi'_r \ast \psi')(x) dx$ converge to $A_{n,m}$ as 
$r \to 0$ and  for $r>0$ small enough the linear map associated to $(A_{n, m}^r)_{ \|n\|_{\text{par}}, \|m \|_{\text{par}} < p}$ is still invertible.
This implies in particular the existence of coefficients $(a_{m})$ such that 
\begin{align*}
\sum_{\|m\|_{\text{par}} < p} A^r_{ n,m} a_{m} &= \sum_{\|m\|_{\text{par}} < p} a_m    \int x^{n} \partial^m (\psi'_r \ast \psi')(x) dx  \\
&= \left\{\begin{array}{ll}
1 \qquad \text{if } n =0 \\
0 \qquad \text{else} 
\end{array}
\right\}.
\end{align*}
The identity \eqref{AMM3} thus follows for $\omega^0 = \sum_{\|m \|_{\text{par}} < p} a_{m}   \partial^m \psi'_r$. 

\medskip

The key ingredient for the proof of \eqref{AMM4} and \eqref{newbound111} are the following 
estimates \eqref{yyy35}--\eqref{yyy36}. We claim that for an arbitrary Schwartz function $\omega$ and any multi-index $m = (m_1,m_2)$ with $\| m \|_{\text{par}} \leq p+1$
we have for any $T>0 $
\begin{align}
\label{yyy35}
\int |\partial^m (\omega - \varphi_T' \ast \omega) | &\leq C_0   \int |\partial^m \omega|,\\
\notag
\int d^\alpha(0,x)  |\partial^m (\omega & - \varphi_T' \ast \omega) | dx \\
\label{yyy37}
&\leq C_0  \Big( \int d^\alpha(0,x) |\partial^m \omega| dx + (T^{\frac14})^\alpha\int  |\partial^m \omega| dx \Big) .
\end{align}
Furthermore,  for $T \leq 1$
\begin{align}
\label{yyy34}
\int |\omega - \varphi_T' \ast \omega | &\leq C_0 (T^{\frac14})^p  \sum_{\| m\|_{\text{par}} = p, p+1} \int |\partial^m \omega|\\
\notag
\int d^\alpha(0,x)  | \omega - \varphi_T' \ast \omega | &\\
\label{yyy36}
\leq C_0 (T^{\frac14})^p  \sum_{\| m\|_{\text{par}} = p, p+1}   \Big( &\int d^\alpha(0,x)  |\partial^m \omega|  + (T^{\frac14})^\alpha  \int |\partial^m \omega| \Big),
\end{align}
where we have $C_0 = C_0(p,\varphi')$ in \eqref{yyy35} -- \eqref{yyy36}.
The estimates \eqref{yyy34} and \eqref{yyy36} rely on the Assumption \eqref{AMM3} that $\varphi'$ integrates to zero against monomials of degree $0 < \| n \|_{\text{par}}<p$.
Once these bounds are established, the representation \eqref{AMM4} follows if we define the $\omega^{(k)}$ recursively by
\begin{align*}
\omega^{(0)} = \psi \qquad \text{and } \omega^{(k+1)} = \omega^{(k)} - \varphi_{\theta^{k}}' \ast \omega^{(k)}
\end{align*}
for a $\theta >0$ small enough.
Indeed, iterating \eqref{yyy35} and \eqref{yyy37} yields
\begin{align*}
&\sum_{\| m \|_{\text{par}} = p, p+1} \int (1 + d^\alpha(0,x)) |\partial^ m\omega^{(k)}| dx \\
&\leq (2C_0)^k \sum_{\| m \|_{\text{par}} = p, p+1} \int (1 + d^\alpha(0,x)) |\partial^ m \psi| dx,
\end{align*}
which can then be plugged into \eqref{yyy34} and \eqref{yyy36} to yield
\begin{align*}
&\int (1 + d^\alpha(0,x))   |\omega^{(k+1)}| dx \\
& \leq (2 C_0)^{k+1} (\theta^{\frac{k}{4}})^p \sum_{\| m \|_{\text{par}} = p, p+1} \int (1 + d^\alpha(0,x)) |\partial^ m \psi| dx,
\end{align*}
which in turn yields \eqref{AMM5} and \eqref{newbound111}. The representation then follows by observing 
\begin{align*}
\psi = \omega^{(0)} = \omega^{(0)} \ast \varphi' + \omega^{(1)} = \omega^{(0)} \ast \varphi' + \omega^{(1)} \ast \varphi_{\theta}' + \omega^{(2)} = \ldots.
\end{align*}
which together with \eqref{AMM5} implies that the convergence holds in $L^1$. 
\medskip

The bounds \eqref{yyy35} and \eqref{yyy34}   are provided in the discussion following equation (295) in \cite{GloriaOtto} (up to the parabolic 
scaling which can be included in the same way as in the following argument). Here 
we only present the proofs for  \eqref{yyy37} and  \eqref{yyy36} which follow along similar lines. First of all,
in order to bound $\int d^\alpha(0,x) | \partial^m  \omega - \varphi'_T \ast  \partial^m   \omega| dx$ we  make use of the triangle inequality in the form 
$| \partial^m  \omega - \varphi'_T \ast \partial^m \omega   | \leq | \partial^m \omega| + |\varphi'_T \ast \partial^m \omega   |$. The integral resulting from the first 
term then already has the desired form. For the second term, we write $| \varphi'_T \ast  \partial^m   \omega(x) | \leq  \int  |\varphi'_T (x-y) \partial^m   \omega(y)| dy$
and use the triangle inequality once more, this time in the form $d^\alpha(0,x) \leq d^\alpha(0,x-y) + d^\alpha(0,y)$. 
It hence remains to bound the two integrals
\begin{align*}
&\int \int  d^\alpha(0,x-y) |\varphi'_T (x-y) | \; | \partial^m    \omega(y)|  dx dy \\
%&\qquad \qquad  = \int  d^\alpha(0,z) |\varphi'_T (z) | dz \;   \int | \partial^m    \omega(y)|   dy ,   \\
&\qquad \qquad  = (T^{\frac14})^{\alpha} \int  d^\alpha(0,\hat{z}) |\varphi' (\hat{z}) | d \hat{z} \;   \int | \partial^m    \omega(y)|   dy ,   \\
%\end{align*}
%\begin{align*}
&\int \int  d^\alpha(0,y) |\varphi'_T (x-y) | \; | \partial^m    \omega(y)|  dx dy \\
&  \qquad \qquad  \leq \int   |\varphi'_T (z) | dz \;   \int d^\alpha(0,y) | \partial^m    \omega(y)|   dy, 
\end{align*}
and estimate \eqref{yyy37} follows. 

\medskip

To obtain \eqref{yyy36}, similar to \cite{GloriaOtto} we obtain the pointwise bound
\begin{align}
 \notag
 &|\varphi'_T \ast \omega- \omega |(x) \\
 \label{zzz89}
 & \leq 2 \sum_{\|m \|_{\text{par}}= p, p+1  } \int_0^1  \int d^{\| m \|_{\text{par}}}(0,z)  |\varphi_T'(-z) | \; |\partial^m \omega(x+sz)| dz ds.
\end{align}
We recall the argument from~\cite{GloriaOtto} (adjusted to the case of parabolic scaling): 
First, according to  \eqref{AMM3} $\varphi'$ integrates non-constant  monomials  of (parabolic) degree $<p$ to zero
which permits us to write
$( \varphi'_T \ast \omega- \omega)(x)  = \int \Big( \omega(x+z) - \sum_{\|m\|_{\text{par}} < p} \frac{1}{m_1!m_2!}\partial^m \omega(x) z^m \Big)$ $ \varphi'_T(-z)dz$.
At this point we seek to apply Taylor's formula, but unlike \cite{GloriaOtto} we need an anisotropic version of the error term. In order to formulate this  we define for $m = (m_1, m_2)$ 
\begin{align*}
F^m := \frac{\partial^m \omega(x) z^m}{ (m_1 + m_2)!} \qquad E^m := \int_0^1 \frac{(1-s)^{m_1+m_2-1}}{(m_1+m_2-1)!} z^m \partial^m \omega(x+s z) d s,
\end{align*}
and observe the elementary identities $\omega(x+z) - \omega(x) = E^{(1,0)} + E^{(0,1)}$ as well as $E^m = F^m + E^{(m_1+1, m_2)} + E^{(m_1,m_2+1)}$ which permit to recursively obtain
\begin{align*}
&\Big|  \omega(x+z) - \sum_{\|m\|_{\text{par}} < p} \frac{1}{m_1! m_2!}\partial^m \omega(x) z^m \Big|  \\
& = \Big| \sum_{\| m\|_{\text{par}} = p}    {{m_1+m_2}\choose{m_1}}  E^{(m_1, m_2)} + \!\!\!\!  \sum_{\| m\|_{\text{par}} = p-1}   {{m_1+m_2}\choose{m_1}}  E^{(m_1, m_2+1)} \Big|\\
& \leq \sum_{\| m\|_{\text{par}} = p, p+1} {{m_1+m_2}\choose{m_1}}  \big| E^{(m_1, m_2)} \big|. 
\end{align*}
Then bounding $|z^m| \leq d^{\| m\|_{\text{par}}}(0,z)$ and observing that  the combinatorial pre-factor satisfies $\frac{1}{(m_1 + m_2-1)!}{{m_1+m_2}\choose{m_1}} \leq 2$ and dropping $(1-s)^{m_1 + m_2-1} \leq 1$ 
the claimed inequality \eqref{zzz89} follows. 

\medskip
To bound $\int d^\alpha(0, x)  |\varphi'_T \ast \omega- \omega |(x) dx$ 
we then use the triangle inequality in the form $ d^\alpha(0, x) \leq d^\alpha(0,z)+ d^\alpha(0,x+sz) $ which prompts to bound the two integrals 
\begin{align*}
&\int  \int_0^1  \int  d^{\alpha+\|m \|_{\text{par}}}(0,z) | \varphi_T'(-z) | \; |\partial^m \omega(x+sz)| dz ds dx \\
&=   \Big(\int d^{\alpha+\|m \|_{\text{par}}}(0,z)  | \varphi_T'(-z) | dz \Big) \; \Big(  \int | \partial^m \omega(x)|dx  \Big) ,  \\
&\int   \int_0^1  \int d^\alpha(0,x+sz) d^{\|m \|_{\text{par}}}(0,z)  |\varphi_T'(-z) |  \; |\partial^m \omega(x+sz)| dz ds dx \\
& =   \Big(  \int d^{\|m \|_{\text{par}}}(0,z)  | \varphi_T'(-z) | dz \Big)\;  \Big( \int  d^\alpha(0,y) |\partial^m \omega(y)|dy \Big) ,
\end{align*}
both of which are bounded as claimed  in \eqref{yyy36}.

%\bigskip
\subsection*{Acknowledgements}
HW is supported by the Royal Society
through the University Research Fellowship UF140187.

\subsection*{Conflict of interest} The authors declare that they have no conflict of interest.

\end{document}